\documentclass[12pt]{amsart}
\usepackage{amsfonts,amssymb,amscd}
\newcommand{\R}{\mathbb R}
\newcommand{\Y}{\mathbb Y}
\newcommand{\T}{\mathbb T}
\newcommand{\A}{\mathbb A}
\newcommand{\B}{\mathbb B}
\newcommand{\Z}{\mathbb Z}
\newcommand{\N}{\mathbb N}

\newcommand{\te}{\theta}
\newcommand{\ti}{\tilde}
\newcommand{\om}{\omega}
\newcommand{\sm}{\sigma}
\newcommand{\Sm}{\Sigma}
\newcommand{\si}{\sigma_k^{-1}}
\newcommand{\la}{\lambda}

\newcommand{\s}{\sim}
\newcommand{\es}{\emptyset}
\newcommand{\bd}{\partial}
\newcommand{\st}{\stackrel}

\newcommand{\ap}{\approx}
\newcommand{\al}{\alpha}
\newcommand{\be}{\beta}
\newcommand{\ga}{\gamma}
\newcommand{\ph}{\varphi}

\newcommand{\Ga}{\Gamma}
\newcommand{\bu}{\bullet}
\newcommand{\vl}{\; | \;}
\newcommand{\ed}{\hfill $\blacksquare$}
\newcommand{\wt}[1]{\widetilde{#1}}
\newcommand{\ab}[1]{\langle #1 \rangle}

\newcommand{\lra}{\leftrightarrow}

\newcommand{\e}{\varepsilon}

\newtheorem{lemma}{Lemma}
\newtheorem{proposition}{Proposition}
\newtheorem{claim}{Claim}

\theoremstyle{definition}

\title{ Three-page encoding and complexity theory
 for spatial graphs}
\author[kurlin]{V.~Kurlin$^*$}
\address{Institut de Math\'ematiques de Bourgogne,
 BP 47870, 21078 Dijon cedex, France}
\email{ kourline@topolog.u-bourgogne.fr, vak26@yandex.ru }
\subjclass[2000]{57M25, 57M15, 57M05}
\keywords{Spatial graph, ambient isotopy,
 three-page embedding, encoding semigroup, graph tangle, vertex sum,
 edge sum, loop sum, three-page complexity, three-letters complexity}
\thanks{$^*$ The author was supported in part by grant INTAS YS 2001/2-30.}

\begin{document}
\vspace*{-12mm}
\begin{center}
\it This is a preprint.
Comments and suggestions are welcome.
\end{center}
\vspace{5mm}

\maketitle
\vspace{-3mm}

\begin{abstract}
For each $n\geq 2$, we construct a finitely presented semigroup $RSG_n$.
The center of $RSG_n$ encodes uniquely
 up to rigid ambient isotopy in $\R^3$
 all non-oriented spatial graphs with
 vertices of degree $\leq n$.
This encoding is obtained by using
 three-page embeddings of graphs
 into the product $\Y=T\times I$, 
 where $T$ is the cone on three points, and $I\ap [0,1]$ is the unit segment.
By exploiting three-page embeddings we introduce the notion of
 the three-page complexity for spatial graphs.
This complexity satisfies the properties of finiteness and
 additivity under natural operations.
\end{abstract}


\section{ Introduction}

\subsection{ Motivations and previous results.}
The notion of a spatial graph is motivated both theoretically and practically.
Firstly, the classification problem of spatial graphs
 up to ambient isotopy in $\R^3$ is a special case of
 the general topological classification of the embeddings into $\R^m$.
The theory of spatial graphs is a natural extension of 
 the classical knot theory to more
 complicated one-dimensional objects.
Many invariants of ordinary links, 
 inlcuding the Alexander polynomial, 
 the Jones polynomial, and 
 Vassiliev finite-type invariants,
 could be generalized to graphs \cite{Li, St, Ya}.
For other equivalence relations on spatial graphs, see \cite{Ta}.
Secondly, spatial graphs are useful mathematical models
 for long protein molecules in molecular biology.
For instance, it is of importance whether 
 a molecule can take the shape  reflection symmetric
 to the original shape or not \cite{Si}.
The isotopy classification problem of spatial graphs
 was intensively studied in \cite{JM, Kau}.
\smallskip

An embedding of a link in a structure, which looks like an open
 book with finitely many pages, was probably considered for
 the first time by Brunn in 1898 \cite{Br}.
More exactly, Brunn proved that each link could be projected into
 the plane $\R^2$ with a unique singular point.
Later an exploiting of such embeddings gave a new link invariant,
 namely \emph{the arc index} \cite{CN}.
It turned out that every link could be embedded into the product $\Y=T\times I$
 (\emph{the book with three pages}),
 where $T$ is the cone on 3 points and $I\ap [0,1]$.

In 1999 with help of \emph{three-page diagrams}
 Dynnikov reduced the isotopy classification of non-oriented links
 in $\R^3$ to a word problem in a finitely presented semigroup \cite{Dy1}.
To be more precise these diagrams will be called 
 \emph{three-page embeddings},
 see the formal definition in Subsection 3.1.
Dynnikov constructed the semigroup $DS$ such that its center encodes
 all non-oriented links up to ambient isotopy in $\R^3$.
In 2002 Dynnikov described an algorithm that allows
 to recognize the unknot using arc presentations
 related closely to three-page embeddings \cite{Dy2}.
Vershinin and the author have already extended the three-page approach to
 spatial \emph{3-graphs} (graphs with vertices of
 degree only 3) \cite{Ku1} and
 to \emph{singular knots} (links with finitely many double intersections in
 general position) \cite{VK}.
In 1990 Matveev introduced a complexity of 3-dimensional manifolds, which
 satisfies the properties of finiteness and additivity
 under connected sum \cite{Mat}.


\subsection{ Basic definitions.}

Any finite 1-dimensional CW-complex $G$ is called
 a \emph{finite graph}.
Every 0-dimensional (respectively, 1-dimensional) cell of $G$ 
 is said  to  be \emph{a vertex} (respectively, \emph{an edge})
 of the graph $G$.
Since \emph{hanging edges} having an endpoint of degree 1 could
 not be knotted, they are excluded.
All graphs are considered up to homeomorphism.
\smallskip

{\bf Definition 1} (\emph{$k$-vertices, $n$-graphs, $J$-graphs}).
A vertex $A\in G$ is called \emph{a $k$-vertex} 
 (i.e. $A$ has \emph{the degree} $k$),
 if the graph $G$ has exactly $k$ edges with the  endpoint $A$.
Fix an integer $n\geq 2$.
If a graph $G$ has $k$-vertices for $k=2,\ldots,n$ only,
 then $G$ is said to be \emph{an $n$-graph}.
Let $J=\{j_1,\ldots,j_k\}$ be any set of integers $j_l\geq 3$.
If a graph $G$ has $k$-vertices, where
 either $k=2$ or $k\in J$, then 
 $G$ will be called \emph{a $J$-graph}.
\ed
\smallskip

We consider only non-oriented graphs, possibly non-connected.
Self-loops and multiple edges are allowed. 
For example, any 2-graph is a disjoint union of several circles.
\smallskip

{\bf Definition 2} (\emph{spatial graphs, isotopies}).
Let $G$ be a finite graph.
We work in the PL-category, i.e. the edges of $G$
 under embeddings in $\R^3$ become finite polygonal lines.

(a) \emph{A spatial} (or \emph{knotted}) graph is
 a subset $G\subset\R^3$, homeomorphic to $G$.
We also assume that a neighbourhood of any vertex of $G$ lies in
 a plane.

(b) \emph{An ambient PL-isotopy} between two spatial graphs
 $G,H\subset \R^3$ is a continuous family of PL-homeomorphisms
 $\phi_t : \R^3 \to \R^3$, $t\in [0,1]$, such that
 $\phi_0=$ id and $\phi_1(G)=H$.

(c) If in addition, at each moment $t\in [0,1]$ of the isotopy $\phi_t$,
 a neighbourhood of every vertex of the graph $\phi_t(G)$
 lies in a plane depending on $t$,
 then $\phi_t$ is called \emph{a rigid isotopy}.
Otherwise $\phi_t$ is said to be \emph{a non-rigid isotopy}.
\ed
\smallskip

For example, \emph{singular knots} are spatial $\{4\}$-graphs
 considered up to rigid isotopy \cite{Ge}.
See a singular knot and a 3-graph in Fig.~4a, 4b, respectively (Section~3).
For spatial 3-graphs, any non-rigid isotopy could be transformed
 into a rigid one\footnote{
 Indeed, one should keep three arcs at each 3-vertex of a given
  graph in a (non-constant) plane.}.
For arbitrary spatial $n$-graphs with $n>3$,
 a non-rigid isotopy can permute edges at any vertex.
Except Subsection~5.4 we shall consider 
 spatial $n$-graphs up to rigid isotopy only.
\smallskip

{\bf Definition 3} (\emph{the encoding alphabet $\A_n$}).
For each $n\geq 2$, let us consider 
 the following \emph{encoding alphabet}:
 $\A_n=\{ \; a_i, \; b_i, \; c_i, \; d_i, \; x_{m,i} \vl
 i\in \Z_3, \; 3\leq m\leq n \}.$
Always the index $i$ belongs to the group $\Z_3=\{0,1,2\}$.
In particular, for $n=2$, we get the Dynnikov alphabet from \cite{Dy1}:
  $\A_2=\{\; a_0, \; a_1, \; a_2, \; b_0, \; b_1, \; b_2, \;
             c_0, \; c_1, \; c_2, \; d_0, \; d_1, \; d_2 \; \}.$
The alphabet $\A_n$ contains exactly $3(n+2)$ letters.
\ed
\smallskip

{\bf Definition 4} (\emph{the encoding semigroups $RSG_n$ and $NSG_n$}).
Let $RSG_n$ be the semigroup generated by
 the letters of $\A_n$ and the relations $(1)-(10)$.
Everywhere the integer parameters $m,p,q$ will satisfy the inequalities
 $3\leq m\leq n$, $2\leq p\leq \frac{n+1}{2}$, and
 $2\leq q\leq \frac{n}{2}$.

$$\begin{array}{l}
(1) \quad d_0 d_1 d_2 = 1;  \\

(2) \quad b_i d_i = d_i b_i =1; \\

(3) \quad a_i = a_{i+1} d_{i-1}, \quad b_i = a_{i-1} c_{i+1}, \quad
           c_i = b_{i-1} c_{i+1}, \quad d_i = a_{i+1} c_{i-1}; \\

(4) \quad x_{2p-1,i-1} = d_{i-1}^{p-1} (x_{2p-1,i} d_{i+1}) b_{i-1}^{p-2},
    \quad x_{2q,i-1} = d_{i-1}^{q-2} (b_{i+1} x_{2q,i} d_{i+1})
          b_{i-1}^{q-2}; \\

(5) \quad x_{2p-1,i} d_i^{p-1} = a_i (x_{2p-1,i} d_i^{p-1}) c_i,
    \quad b_i^{p-1} x_{2p-1,i} b_i = a_i (b_i^{p-1} x_{2p-1,i} b_i) c_i; \\

(6) \quad d_i x_{2q,i} d_i^{q-1} = a_i (d_i x_{2q,i} d_i^{q-1}) c_i,
    \quad b_i^{q-1} x_{2q,i} b_i = a_i (b_i^{q-1} x_{2q,i} b_i) c_i; \\

(7) \quad (d_{i} c_{i}) w = w (d_{i} c_{i}), \mbox{ where }
     w\in \{c_{i+1}, \; b_{i} d_{i+1} d_{i}, \; x_{m,i+1} \}; \\

(8) \quad uv=vu, \mbox{ where }
 u\in \{a_i b_i, \; b_{i-1} d_i d_{i-1} b_i, \;
        x_{2p-1,i} b_i, \; d_i x_{2q,i} b_i \}, \\ \qquad
 v\in \{ a_{i+1}, \; b_{i+1}, \; c_{i+1}, \; b_{i} d_{i+1} d_{i}, \;
         x_{m,i+1} \}; \\

(9) \quad  (x_{2p-1,i} b_i) D_{p,i}  =  D_{p-1,i} (x_{2p-1,i} b_i),
     \mbox{ where } D_{k,i}=d_i^k d_{i+1}^k d_{i-1}^k (k\geq 1); \\

(10) \quad (d_i x_{2q,i} b_i) D_{q,i}  =  D_{q,i} (d_i x_{2q,i} b_i). 
  \end{array} $$
Let us introduce the semigroup $NSG_n$
 generated by the letters of $\A_n$ and the relations $(1)-(8)$,
 $(9') \; x_{m,i} b_i (d_i^2 d_{i+1}^2 d_{i-1}^2)=x_{m,i} b_i$.
For any set $J=\{j_1,\ldots,j_k\}$ of integers $j_l\geq 3$,
 by $RSG_J$ we denote the semigroup generated by
 the letters $\{a_i,b_i,c_i,d_i,x_{m,i} \vl$ $i\in\Z_3,\; m\in J \}$
 and the relations $(1)-(10)$ containing only these letters.
Let the semigroup $NSG_J$ be generated by the same letters and
 the relations $(1)-(8)$, $(9')$ for $m\in J$.
\ed
\smallskip

The semigroups $RSG_n$ and $NSG_n$ are monoids,
 the empty word $\es$ is the identity element.
In Section~2, we shall give a geometric interpretation of
 the letters of $\A_n$ and the relations (1)--(10).
One of the relations in (2) is superfluous: it can be obtained from
 the relation (1) and the other relations in (2).
Then the total number of the relations (1)--(10) is $3(n^2+7n-2)$.
The semigroups $RSG_2=NSG_2$
 generated by 12 letters $a_i,b_i,c_i,d_i$ ($i\in \Z_3$) and
 48 defining relations $(1)-(3),(7)-(8)$,
 which contain only the letters of $\A_2$,
 coincide with the Dynnikov semigroup $DS$ from \cite{Dy1}.
The semigroups $RSG_3\cong NSG_3$ and $RSG_{\{4\}}\not\cong NSG_{\{4\}}$
 are generated by 15 letters and 84 relations.
By $|J|$ denote the number of elements of a set
 $J=\{j_1,\ldots,j_k\}$, $j_l\geq 3$.
Then the semigroups $RSG_J$ and $NSG_J$ are generated by
 $3(4+|J|)$ letters and $3(16+11|J|+|J|^2)$ relations.
\smallskip

{\bf Definition 5} (\emph{the automorphisms $\rho_n,\e_n$,
 the mirror image}).
Consider the following map of the letters of $\A_n$:
$\rho(a_i)=c_i$, 
$\rho(b_i)=d_i$,
$\rho(c_i)=a_i$,
$\rho(d_i)=b_i$,
$\rho(x_{2p-1,i})=x_{2p-1,i} b_i c_i$, 
$\rho(x_{2q,i})=x_{2q,i}$.
Let $W_n$ be the set of all words in the alphabet $\A_n$.
By the formula $\rho(uv)=\rho(v)\rho(u)$ (where $u,v\in W_n$) 
 the map $\rho$ extends to \emph{the involutive automorphisms}
 $\rho_n: RSG_n\to RSG_n$ and $\e_n: NSG_n\to NSG_n$.
Similarly, one can define the morphisms
 $\rho_J: RSG_J\to RSG_J$ and $\e_J: NSG_J\to NSG_J$.
\emph{The mirror image} of a spatial graph $G\subset\R^3$ is
 the spatial graph $\bar G\subset\R^3$ reflection symmetric to $G$.
\ed

 
\subsection{ Main results}

Theorems 1--2 reduce the isotopy classification of arbitrary spatial graphs
 to a pure algebraic word problem in finitely presented semigroups.
\smallskip

{\bf Theorem~1a.}
Any spatial $n$-graph $G\subset\R^3$
 is encoded by an element $w_G\in RSG_n$.
\smallskip

{\bf Theorem~1b.}
Two arbitrary spatial $n$-graphs $G,H\subset\R^3$ are
 \emph{rigidly isotopic} in $\R^3$
 if and only if the corresponding elements of 
 the semigroup $RSG_n$ are equal: $w_{G}=w_{H}$.
\smallskip

{\bf Theorem~1c.}
An element $w\in RSG_n$ encodes a spatial $n$-graph
 if and only if the element $w$ is central, i.e. 
 $w$ commutes with each element of $RSG_n$.
Moreover, there is an algorithm to decide,
 whether a given element $w\in RSG_n$ is central,
 which is linear in the length of the word $w$.
\smallskip

Theorem~1 means that
 the center of $RSG_n$ encodes uniquely all spatial $n$-graphs
 up to rigid isotopy in $\R^3$.
Proposition~1 in Subsection~5.2 shows
 that the whole semigroup $RSG_n$ describes a wider class of
 \emph{rigid three-page tangles},
 see Definition~14 in Subsection~4.3.
\smallskip

{\bf Theorem~2.}
The center of the semigroup $NSG_n$ encodes
 all spatial $n$-graphs up to \emph{non-rigid isotopy} in $\R^3$.
Moreover, there is an algorithm to decide,
 whether a given element $v\in NSG_n$ is central,
 which is linear in the length of the word $v$.
\smallskip

{\bf Theorem~3.}
Let $\{G\}$ be the set of all non-oriented spatial
 graphs\footnote{
 In particular, the set $\{G\}$ contains all non-oriented knots
  and links.}
 considered up to homeomorphism $f:S^3\to S^3$, which can reverse
 the orientation of $S^3$.
There exists a comlexity function $tp:\{G\}\to\N$
 with the following properties:

(3.1) for any $k\in\N$, there is a finite
 number of spatial graphs $G$ with $tp(G)=k$;

(3.2) there are well-defined operations on spatial graphs:
 the disjoint union $G\sqcup H$, a vertex sum $G*H$,
 an edge sum $G\vee H$ such that
 $tp(G\sqcup H)=tp(G*H)=tp(G)+tp(H)+2$ and
 $tp(G\vee H)=tp(G)+tp(H)+3$.
\smallskip

Theorem 3 was motivated by the notion of Matveev's
 complexity for 3-dimensional manifolds \cite{Mat}.
Theorems~1--2 imply several algebraic and geometric corollaries.
\smallskip

{\bf Corollary 1.}
a) Let a spatial graph $G\subset \R^3$
 be encoded by $w_G\in RSG_n$.
The graph $G$ is \emph{rigid isotopic} to its mirror image $\bar G\subset\R^3$
 if and only if $\rho_n(w_G)=w_G$ in $RSG_n$.
\smallskip

b) Let a spatial $n$-graph $G\subset \R^3$
 be encoded by $v_G\in NSG_n$.
The graph $G$ is \emph{non-rigid} isotopic to its mirror image
 $\bar G\subset\R^3$ if and only if we have $\e_n(v_G)=v_G$ in $NSG_n$.
\smallskip

An element $w$ of a semigroup $S$ is \emph{invertible},
 if $w$ has left and right inverses in $S$.
\smallskip

{\bf Corollary 2.}
a) For any integer $2\leq m<n$, the natural inclusion $RSG_m\to RSG_n$ is
 a monomorphism of semigroups.
The group of the invertible elements of $RSG_n$
 coincides with the Dynnikov group $DG\subset DS=RSG_2$,
 generated by 2 letters and 3 relations:
 $$DG=\langle x,y \vl  [ [x,y], x^2yx^{-2} ] = [ [x,y], y^2xy^{-2} ] =
                       [ [x,y], [x^{-1},y^{-1}] ] = 1 \rangle, \;
   [x,y]=xyx^{-1}y^{-1}.$$
b) For any integer $2\leq m<n$, the natural inclusion $NSG_m\to NSG_n$ is
 a monomorphism of semigroups.
The group of the invertible elements of $NSG_n$ coincides with $DG$.
\smallskip

The commutator subgroup of the group $DG$ is the braid group
 $B_{\infty}$ on infinitely many strings \cite{Dy1}.
The method of three-page embeddings could be applied to any $J$-graph.
\smallskip

{\bf Corollary 3.}
a) The center of the semigroup $RSG_J$ (respectively, $NSG_J$)
 encodes all spatial $J$-graphs up to rigid
 (respectively, non-rigid) isotopy in $\R^3$.
Moreover, there is an algorithm to decide,
 whether a given element $w\in RSG_J$ (respectively, $v\in NSG_J$)
 is central, which is linear in the length of the given word.
\smallskip

b) Let a spatial graph $G\subset \R^3$
 be encoded by $w_G\in RSG_J$ (respectively, by $v_G\in NSG_J$).
The graph $G$ is rigid (respectively, non-rigid) isotopic to
 its mirror image $\bar G\subset\R^3$
 if and only if we have $\rho_J(w_G)=w_G$ in $RSG_n$
 (respectively, $\e_J(v_G)=v_G$ in $NSG_J$).
\smallskip

c) For any subset $K\subset J$, the natural inclusions
 $RSG_K\to RSG_J$ and $NSG_K\to NSG_J$ are monomorphisms of semigroups.
The groups of the invertible elements of the semigroups
 $RSG_J$ and $NSG_J$ coincide with the Dynnikov group $DG$.
\smallskip

The following corollary extends original Brunn's result
 on embeddings of links \cite{Br}.
\smallskip

{\bf Corollary 4.}
Any spatial $J$-graph $G\subset\R^3$ is non-rigid isotopic
 to a spatial graph that could be projected into $\R^2$
 with a unique singular point.
\smallskip

It is well-known that not any finite graph could be topologically
 embedded into $\R^2$.
In what minimal space one can embed all finite graphs?
Theorem~1a implies
\smallskip

{\bf Corollary 5.}
Any finite graph (possibly with hanging edges) could be topologically
 embedded into the book $\Y=T\times I$, where $T$ is the cone on
 three points, $I\ap [0,1]$.


\subsection{ The content of the paper.}
In Section~2, the generators and the relations (1)--(10) of
 the semigroup $RSG_n$ are described in a clear geometric way.
Subsection~2.4 contains a scheme for the proof of Theorem~1.
In Section~3, we give the proofs for Theorem~1a and Corollary~5.
Theorems~1b--1c, 2, and Corollaries~1--4 are proved in Section~5.
The hard part of Theorem~1b is a particular case of
 Proposition~1, which is verified in Subsection~5.2.
In Section~6, we deduce Lemma~3 used in the proof of Proposition~1.
Section~7 is devoted to algebraic and geometric approaches
 to classification of spatial graphs via
 three-page embeddings.
We state several open problems.
Theorem~3 is proved in Subsection~7.3.
\smallskip

{\bf Acknowledgments.}
The author is very grateful to D.~Lines, H.~Morton, L.~Paris,
 and V.~Vershinin for many useful discussions.


\setcounter{section}{1}
\section{ A geometric interpretation of the semigroups $RSG_n$ and $NSG_n$}

\subsection{Reidemeister moves for spatial graphs}

First we formulate an analog of 
 the famous Reidemeister theorem for spatial graphs.
\smallskip

{\bf The generalized Reidemeister theorem} \cite{JM}.
Any spatial graph could be represented by its plane diagram
 up to the generalized Reidemeister moves $R1-R5$ of Fig.~1.
In the case of non-rigid isotopy,
 the move $R5'$ instead of $R5$ is valid.
\qed
\smallskip

Subdivisions of edges and adding of zigzags are omitted.
In Fig.~1, "dots" between two arcs denote any finite number of arcs.
The Reidemeister moves are local, but two-dimensional.
Theorem~1b states that the moves $(1)-(10)$ on words in $\A_n$
 also generate any rigid isotopy of graphs.
Moreover, the moves $(1)-(10)$ are local and 1-dimensional.

\begin{picture}(450,75)(10,0)

\put(0,55){\line(0,1){10}}
\put(0,15){\line(0,1){10}}
\put(15,25){\line(0,1){30}}
\put(0,55){\line(1,-2){15}}
\put(0,25){\line(1,2){5}}
\put(10,45){\line(1,2){5}}
\put(25,40){$\sim$}
\put(40,15){\line(0,1){50}}
\put(50,40){$\sim$}
\put(80,55){\line(0,1){10}}
\put(80,15){\line(0,1){10}}
\put(65,25){\line(0,1){30}}
\put(65,25){\line(1,2){15}}
\put(80,25){\line(-1,2){5}}
\put(70,45){\line(-1,2){5}}
\put(35,0){$R1$}

\put(110,70){\line(1,-2){15}}
\put(125,40){\line(-1,-2){15}}
\put(125,70){\line(-1,-2){5}}
\put(110,40){\line(1,2){5}}
\put(110,40){\line(1,-2){5}}
\put(120,20){\line(1,-2){5}}
\put(135,40){$\sim$}
\put(150,15){\line(0,1){50}}
\put(160,15){\line(0,1){50}}
\put(170,40){$\sim$}
\put(185,40){\line(1,2){15}}
\put(185,40){\line(1,-2){15}}
\put(185,70){\line(1,-2){5}}
\put(195,50){\line(1,-2){5}}
\put(200,40){\line(-1,-2){5}}
\put(185,10){\line(1,2){5}}
\put(150,0){$R2$}

\put(230,70){\line(1,-2){5}}
\put(240,50){\line(1,-2){15}}
\put(240,20){\line(1,2){5}}
\put(250,40){\line(1,2){5}}
\put(260,60){\line(1,2){5}}
\put(230,55){\line(1,0){35}}
\put(275,45){$\sim$}
\put(290,35){\line(1,0){35}}
\put(300,70){\line(1,-2){15}}
\put(320,30){\line(1,-2){5}}
\put(290,20){\line(1,2){5}}
\put(300,40){\line(1,2){5}}
\put(310,60){\line(1,2){5}}

\put(350,70){\line(1,-2){25}}
\put(360,20){\line(1,2){5}}
\put(370,40){\line(1,2){15}}
\put(362,55){\line(1,0){11}}
\put(345,55){\line(1,0){8}}
\put(382,55){\line(1,0){8}}
\put(395,45){$\sim$}
\put(422,35){\line(1,0){11}}
\put(405,35){\line(1,0){8}}
\put(442,35){\line(1,0){8}}
\put(420,70){\line(1,-2){25}}
\put(410,20){\line(1,2){15}}
\put(430,60){\line(1,2){5}}
\put(270,0){$R3$}
\put(390,0){$R3$}
\end{picture}


\begin{picture}(450,85)(20,0)
\put(70,-35){ {\bf Fig.~1.}
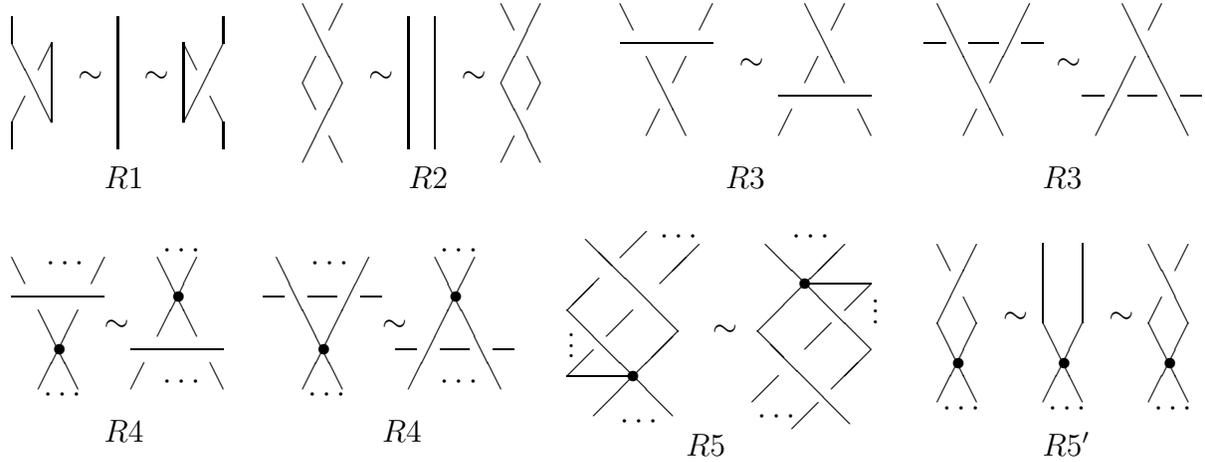
 The Reidemeister moves for spatial graphs in $\R^3$.}

\put(10,60){\line(1,-2){5}}
\put(20,40){\line(1,-2){15}}
\put(20,10){\line(1,2){15}}
\put(40,50){\line(1,2){5}}
\put(10,45){\line(1,0){35}}
\put(28,25){\circle*{4}}
\put(23,55){$\cdots$}
\put(22,5){$\cdots$}
\put(45,32){$\sim$}
\put(55,25){\line(1,0){35}}
\put(65,60){\line(1,-2){15}}
\put(85,20){\line(1,-2){5}}
\put(55,10){\line(1,2){5}}
\put(65,30){\line(1,2){15}}
\put(73,45){\circle*{4}}
\put(67,10){$\cdots$}
\put(67,60){$\cdots$}
\put(45,-10){$R4$}

\put(110,60){\line(1,-2){25}}
\put(120,10){\line(1,2){25}}
\put(122,45){\line(1,0){11}}
\put(105,45){\line(1,0){8}}
\put(142,45){\line(1,0){8}}
\put(128,25){\circle*{4}}
\put(123,55){$\cdots$}
\put(122,5){$\cdots$}
\put(150,32){$\sim$}
\put(172,25){\line(1,0){11}}
\put(155,25){\line(1,0){8}}
\put(192,25){\line(1,0){8}}
\put(170,60){\line(1,-2){25}}
\put(160,10){\line(1,2){25}}
\put(178,45){\circle*{4}}
\put(172,10){$\cdots$}
\put(172,60){$\cdots$}
\put(150,-10){$R4$}

\put(245,15){\circle*{4}}
\put(240,-5){$\cdots$}
\put(230,0){\line(1,1){32}}
\put(245,15){\line(-1,0){25}}
\put(220,22){$\vdots$}
\put(220,15){\line(1,1){10}}
\put(235,30){\line(1,1){10}}
\put(255,50){\line(1,1){15}}
\put(220,40){\line(1,1){12}}
\put(240,60){\line(1,1){10}}
\put(255,65){$\cdots$}
\put(260,0){\line(-1,1){40}}
\put(262,32){\line(-1,1){35}}
\put(275,30){$\sim$}
\put(310,50){\circle*{4}}
\put(305,65){$\cdots$}
\put(310,50){\line(1,0){25}}
\put(335,35){$\vdots$}
\put(335,50){\line(-1,-1){10}}
\put(320,35){\line(-1,-1){10}}
\put(300,15){\line(-1,-1){10}}
\put(292,32){\line(1,1){33}}
\put(292,-3){$\cdots$}
\put(327,-3){\line(-1,1){35}}
\put(335,25){\line(-1,1){40}}
\put(335,25){\line(-1,-1){15}}
\put(315,5){\line(-1,-1){10}}
\put(265,-15){$R5$}

\put(368,20){\circle*{4}}
\put(360,5){\line(1,2){15}}
\put(375,5){\line(-1,2){15}}
\put(360,35){\line(1,2){15}}
\put(375,35){\line(-1,2){5}}
\put(365,55){\line(-1,2){5}}
\put(385,35){$\sim$}
\put(408,20){\circle*{4}}
\put(400,5){\line(1,2){15}}
\put(415,5){\line(-1,2){15}}
\put(400,35){\line(0,1){30}}
\put(415,35){\line(0,1){30}}
\put(425,35){$\sim$}
\put(448,20){\circle*{4}}
\put(440,5){\line(1,2){15}}
\put(455,5){\line(-1,2){15}}
\put(455,35){\line(-1,2){15}}
\put(440,35){\line(1,2){5}}
\put(450,55){\line(1,2){5}}
\put(362,0){$\cdots$}
\put(402,0){$\cdots$}
\put(442,0){$\cdots$}
\put(400,-15){$R5'$}

\end{picture}
\vspace{6mm}

 
\subsection{Geometric interpretation of the alphabet $\A_n$}

The alphabet $\A_n$ was introduced in Definition~3.
Here we describe geometrically the letters of $\A_n$.
\smallskip

{\bf Definition 6} (\emph{the book $\Y$, the pages $P_i$, the axis $\al$}).
\emph{The book with three pages} is the product $\Y=T\times I$, 
 where $T$ is the cone on three points, and $I\ap [0,1]$ is
 the oriented segment.
The interval $I-\bd I$ is homeomorphic to the line $\R$ and is said to be
 \emph{the axis $\al$}.
The set $\Y-\bd\Y$ is the union $P_0\cup P_1\cup P_2$ of
 three half-planes with the common
 oriented boundary $\bd P_0=\bd P_1=\bd P_2=\al$.
The half-planes $P_i$ will be called \emph{the pages} of $\Y$.
\ed
\smallskip

In Fig.~2a and 2b, every letter of $\A_n$ encodes
 a local embedding into the book $\Y$.
In these figures, the page $P_0$ is above the axis $\al$,
 the pages $P_1$, $P_2$ are below $\al$, and $P_2$ is below $P_1$,
  i.e. arcs in $P_2$ are drawn by dashed lines.
In Fig.~2a, every letter with the index
 $i\in\Z_3=\{0,1,2\}$ denotes
 an embedding of two arcs into the disk $P_{i-1}\cup P_{i+1}$.
In Fig.~2b, each letter $x_{m,i}$ encodes an embedding of a neighbourhood
 of an $m$-vertex into the bowed disk $P_{i-1}\cup P_{i+1}$.
More exactly, one of the arcs at each $(2p-1)$-vertex $A_1$ lies in $P_{i-1}$
 and points toward to the positive direction of the axis $\al$.
All the other arcs at $A_1$ lies in $P_{i+1}$.
Exactly $p-1$ of these arcs point toward to
 the positive direction of $\al$, and another ones
 point toward to the negative direction of $\al$.
Similarly, two arcs at any $2q$-vertex $A_2$ lie in $P_{i-1}$,
 one of them points toward to the positive direction of $\al$,
 and another one points toward to the negative direction of $\al$.
Also the other $2q-2$ arcs at $A_2$ lie in $P_{i+1}$, exactly
 $q-1$ of them point toward to the the positive direction of $\al$,
 and the other ones point toward to the negative direction of $\al$.
The disk $P_{i-1}\cup P_{i+1}$ does not lie in a plane.
For rigid spatial graphs,
 one can assume that during any rigid isotopy a neighbourhood
 of each vertex lies in such a bowed disk.
Attaching one local picture of Fig.~2a or 2b to other one
 according to the direction of the axis $\al$,
 one can obtain any word $w$ in $\A_n$.


\begin{picture}(450,90)(0,0)

\put(0,60){\vector(1,0){80}}
\put(30,60){\circle*{3}}
\put(-10,80){\line(1,-2){30}}
\put(50,80){\line(1,-2){30}}
\put(0,60){\line(-1,-2){15}}
\put(60,60){\line(-1,-2){5}}
\put(45,30){\line(1,2){5}}
\put(-15,30){\line(1,0){25}}
\put(20,30){\line(1,0){10}}
\put(35,30){\line(1,0){10}}
\put(-10,80){\line(1,0){60}}
\put(20,20){\line(1,0){60}}
\put(75,65){$\al$}
\put(0,70){$P_0$}
\put(60,25){$P_1$}
\put(-5,35){$P_2$}
\put(30,5){$a_0$}
{\thicklines
\put(30,60){\line(3,-2){30}}
\put(30,60){\line(1,-2){5}}
\put(37,45){\line(1,-2){5}} }

\put(120,60){\vector(1,0){80}}
\put(150,60){\circle*{3}}
\put(110,80){\line(1,-2){30}}
\put(170,80){\line(1,-2){30}}
\put(120,60){\line(-1,-2){15}}
\put(180,60){\line(-1,-2){5}}
\put(165,30){\line(1,2){5}}
\put(105,30){\line(1,0){25}}
\put(140,30){\line(1,0){10}}
\put(155,30){\line(1,0){10}}
\put(110,80){\line(1,0){60}}
\put(140,20){\line(1,0){60}}
\put(195,65){$\al$}
\put(120,70){$P_0$}
\put(180,25){$P_1$}
\put(115,35){$P_2$}
\put(150,5){$b_0$}
{\thicklines
\put(150,60){\line(3,-2){30}}
\put(150,60){\line(-1,-2){5}}
\put(142,45){\line(-1,-2){5}} }

\put(240,60){\vector(1,0){80}}
\put(270,60){\circle*{3}}
\put(230,80){\line(1,-2){30}}
\put(290,80){\line(1,-2){30}}
\put(240,60){\line(-1,-2){15}}
\put(300,60){\line(-1,-2){5}}
\put(285,30){\line(1,2){5}}
\put(225,30){\line(1,0){25}}
\put(260,30){\line(1,0){10}}
\put(275,30){\line(1,0){10}}
\put(230,80){\line(1,0){60}}
\put(260,20){\line(1,0){60}}
\put(315,65){$\al$}
\put(240,70){$P_0$}
\put(300,25){$P_1$}
\put(235,35){$P_2$}
\put(270,5){$c_0$}
{\thicklines
\put(270,60){\line(-2,-1){20}}
\put(270,60){\line(-1,-2){5}}
\put(262,45){\line(-1,-2){5}} }

\put(360,60){\vector(1,0){80}}
\put(390,60){\circle*{3}}
\put(350,80){\line(1,-2){30}}
\put(410,80){\line(1,-2){30}}
\put(360,60){\line(-1,-2){15}}
\put(420,60){\line(-1,-2){5}}
\put(405,30){\line(1,2){5}}
\put(345,30){\line(1,0){25}}
\put(380,30){\line(1,0){10}}
\put(395,30){\line(1,0){10}}
\put(350,80){\line(1,0){60}}
\put(380,20){\line(1,0){60}}
\put(435,65){$\al$}
\put(360,70){$P_0$}
\put(420,25){$P_1$}
\put(355,35){$P_2$}
\put(390,5){$d_0$}
{\thicklines
\put(390,60){\line(-2,-1){20}}
\put(390,60){\line(1,-2){5}}
\put(397,45){\line(1,-2){5}} }

\end{picture}


\begin{picture}(450,100)(0,0)

\put(0,60){\vector(1,0){80}}
\put(30,60){\circle*{3}}
\put(-15,90){\line(1,-2){35}}
\put(45,90){\line(1,-2){35}}
\put(0,60){\line(-1,-2){15}}
\put(60,60){\line(-1,-2){5}}
\put(45,30){\line(1,2){5}}
\put(-15,30){\line(1,0){25}}
\put(20,30){\line(1,0){10}}
\put(35,30){\line(1,0){10}}
\put(-15,90){\line(1,0){60}}
\put(20,20){\line(1,0){60}}
\put(75,65){$\al$}
\put(0,75){$P_0$}
\put(60,25){$P_1$}
\put(-5,35){$P_2$}
\put(30,5){$a_1$}
{\thicklines
\put(30,60){\line(1,2){10}}
\put(30,60){\line(1,-2){5}}
\put(37,45){\line(1,-2){5}} }

\put(120,60){\vector(1,0){80}}
\put(150,60){\circle*{3}}
\put(105,90){\line(1,-2){35}}
\put(165,90){\line(1,-2){35}}
\put(120,60){\line(-1,-2){15}}
\put(180,60){\line(-1,-2){5}}
\put(165,30){\line(1,2){5}}
\put(105,30){\line(1,0){25}}
\put(140,30){\line(1,0){10}}
\put(155,30){\line(1,0){10}}
\put(105,90){\line(1,0){60}}
\put(140,20){\line(1,0){60}}
\put(195,65){$\al$}
\put(120,75){$P_0$}
\put(180,25){$P_1$}
\put(115,35){$P_2$}
\put(150,5){$b_1$}
{\thicklines
\put(140,80){\line(1,-2){15}}
\put(157,45){\line(1,-2){5}} }

\put(240,60){\vector(1,0){80}}
\put(270,60){\circle*{3}}
\put(225,90){\line(1,-2){35}}
\put(285,90){\line(1,-2){35}}
\put(240,60){\line(-1,-2){15}}
\put(300,60){\line(-1,-2){5}}
\put(285,30){\line(1,2){5}}
\put(225,30){\line(1,0){25}}
\put(260,30){\line(1,0){10}}
\put(275,30){\line(1,0){10}}
\put(225,90){\line(1,0){60}}
\put(260,20){\line(1,0){60}}
\put(315,65){$\al$}
\put(240,75){$P_0$}
\put(300,25){$P_1$}
\put(235,35){$P_2$}
\put(270,5){$c_1$}
{\thicklines
\put(270,60){\line(-1,2){10}}
\put(270,60){\line(-1,-2){5}}
\put(262,45){\line(-1,-2){5}} }

\put(360,60){\vector(1,0){80}}
\put(390,60){\circle*{3}}
\put(345,90){\line(1,-2){35}}
\put(405,90){\line(1,-2){35}}
\put(360,60){\line(-1,-2){15}}
\put(420,60){\line(-1,-2){5}}
\put(405,30){\line(1,2){5}}
\put(345,30){\line(1,0){25}}
\put(380,30){\line(1,0){10}}
\put(395,30){\line(1,0){10}}
\put(345,90){\line(1,0){60}}
\put(380,20){\line(1,0){60}}
\put(435,65){$\al$}
\put(360,75){$P_0$}
\put(420,25){$P_1$}
\put(355,35){$P_2$}
\put(390,5){$d_1$}
{\thicklines
\put(400,80){\line(-1,-2){15}}
\put(383,45){\line(-1,-2){5}} }

\end{picture}


\begin{picture}(450,100)(0,0)

\put(0,60){\vector(1,0){80}}
\put(30,60){\circle*{3}}
\put(-15,90){\line(1,-2){35}}
\put(45,90){\line(1,-2){35}}
\put(0,60){\line(-1,-2){15}}
\put(60,60){\line(-1,-2){5}}
\put(45,30){\line(1,2){5}}
\put(-15,30){\line(1,0){25}}
\put(20,30){\line(1,0){10}}
\put(35,30){\line(1,0){10}}
\put(-15,90){\line(1,0){60}}
\put(20,20){\line(1,0){60}}
\put(75,65){$\al$}
\put(0,75){$P_0$}
\put(60,25){$P_1$}
\put(-5,35){$P_2$}
\put(30,5){$a_2$}
{\thicklines
\put(30,60){\line(1,2){10}}
\put(30,60){\line(1,-2){10}} }

\put(120,60){\vector(1,0){80}}
\put(150,60){\circle*{3}}
\put(105,90){\line(1,-2){35}}
\put(165,90){\line(1,-2){35}}
\put(120,60){\line(-1,-2){15}}
\put(180,60){\line(-1,-2){5}}
\put(165,30){\line(1,2){5}}
\put(105,30){\line(1,0){25}}
\put(140,30){\line(1,0){10}}
\put(155,30){\line(1,0){10}}
\put(105,90){\line(1,0){60}}
\put(140,20){\line(1,0){60}}
\put(195,65){$\al$}
\put(120,75){$P_0$}
\put(180,25){$P_1$}
\put(115,35){$P_2$}
\put(150,5){$b_2$}
{\thicklines
\put(160,80){\line(-1,-2){20}} }

\put(240,60){\vector(1,0){80}}
\put(270,60){\circle*{3}}
\put(225,90){\line(1,-2){35}}
\put(285,90){\line(1,-2){35}}
\put(240,60){\line(-1,-2){15}}
\put(300,60){\line(-1,-2){5}}
\put(285,30){\line(1,2){5}}
\put(225,30){\line(1,0){25}}
\put(260,30){\line(1,0){10}}
\put(275,30){\line(1,0){10}}
\put(225,90){\line(1,0){60}}
\put(260,20){\line(1,0){60}}
\put(315,65){$\al$}
\put(240,75){$P_0$}
\put(300,25){$P_1$}
\put(235,35){$P_2$}
\put(270,5){$c_2$}
{\thicklines
\put(270,60){\line(-1,2){10}}
\put(270,60){\line(-1,-2){10}} }

\put(360,60){\vector(1,0){80}}
\put(390,60){\circle*{3}}
\put(345,90){\line(1,-2){35}}
\put(405,90){\line(1,-2){35}}
\put(360,60){\line(-1,-2){15}}
\put(420,60){\line(-1,-2){5}}
\put(405,30){\line(1,2){5}}
\put(345,30){\line(1,0){25}}
\put(380,30){\line(1,0){10}}
\put(395,30){\line(1,0){10}}
\put(345,90){\line(1,0){60}}
\put(380,20){\line(1,0){60}}
\put(435,65){$\al$}
\put(360,75){$P_0$}
\put(420,25){$P_1$}
\put(355,35){$P_2$}
\put(390,5){$d_2$}
{\thicklines
\put(380,80){\line(1,-2){20}} }

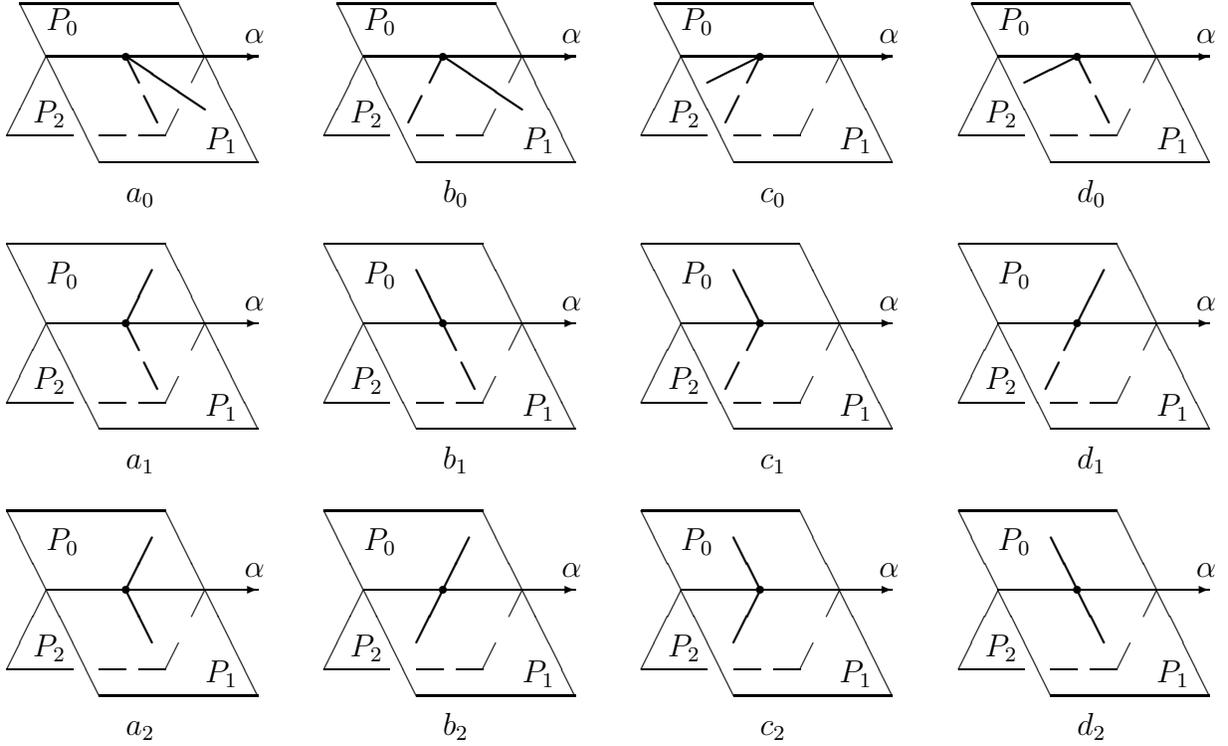
\put(70,-15){ {\bf Fig. 2a.} The Dynnikov letters of the alphabet $\A_2$.}
\end{picture}
\vspace{3mm}


\begin{picture}(450,90)(0,0)

\put(0,60){\vector(1,0){80}}
\put(30,60){\circle*{5}}
\put(-10,80){\line(1,-2){30}}
\put(50,80){\line(1,-2){30}}
\put(0,60){\line(-1,-2){15}}
\put(60,60){\line(-1,-2){5}}
\put(45,30){\line(1,2){5}}
\put(-15,30){\line(1,0){25}}
\put(20,30){\line(1,0){10}}
\put(35,30){\line(1,0){10}}
\put(-10,80){\line(1,0){60}}
\put(20,20){\line(1,0){60}}
\put(75,65){$\al$}
\put(0,70){$P_0$}
\put(60,25){$P_1$}
\put(-5,35){$P_2$}
\put(30,5){$x_{3,0}$}
\put(30,65){$A_1$}
{\thicklines
\put(30,60){\line(3,-2){30}}
\put(30,60){\line(-1,-1){15}}
\put(30,60){\line(1,-2){5}}
\put(37,45){\line(1,-2){5}} }

\put(120,60){\vector(1,0){80}}
\put(150,60){\circle*{5}}
\put(110,80){\line(1,-2){30}}
\put(170,80){\line(1,-2){30}}
\put(120,60){\line(-1,-2){15}}
\put(180,60){\line(-1,-2){5}}
\put(165,30){\line(1,2){5}}
\put(105,30){\line(1,0){25}}
\put(140,30){\line(1,0){10}}
\put(155,30){\line(1,0){10}}
\put(110,80){\line(1,0){60}}
\put(140,20){\line(1,0){60}}
\put(195,65){$\al$}
\put(120,70){$P_0$}
\put(180,25){$P_1$}
\put(115,35){$P_2$}
\put(150,5){$x_{4,0}$}
\put(150,65){$A_2$}
{\thicklines
\put(150,60){\line(3,-2){30}}
\put(150,60){\line(-2,-1){20}}
\put(150,60){\line(1,-2){5}}
\put(157,45){\line(1,-2){5}}
\put(150,60){\line(-1,-2){5}}
\put(142,45){\line(-1,-2){5}} }

\put(240,60){\vector(1,0){80}}
\put(270,60){\circle*{5}}
\put(230,80){\line(1,-2){30}}
\put(290,80){\line(1,-2){30}}
\put(240,60){\line(-1,-2){15}}
\put(300,60){\line(-1,-2){5}}
\put(285,30){\line(1,2){5}}
\put(225,30){\line(1,0){25}}
\put(260,30){\line(1,0){10}}
\put(275,30){\line(1,0){10}}
\put(230,80){\line(1,0){60}}
\put(260,20){\line(1,0){60}}
\put(315,65){$\al$}
\put(240,70){$P_0$}
\put(300,25){$P_1$}
\put(235,35){$P_2$}
\put(270,5){$x_{5,0}$}
\put(270,65){$A_1$}
{\thicklines
\put(270,60){\line(-2,-1){20}}
\put(270,60){\line(-1,-1){15}}
\put(270,60){\line(5,-4){25}}
\put(270,60){\line(2,-1){30}}
\put(270,60){\line(1,-2){5}}
\put(277,45){\line(1,-2){5}} }

\put(360,60){\vector(1,0){80}}
\put(390,60){\circle*{5}}
\put(350,80){\line(1,-2){30}}
\put(410,80){\line(1,-2){30}}
\put(360,60){\line(-1,-2){15}}
\put(420,60){\line(-1,-2){5}}
\put(405,30){\line(1,2){5}}
\put(345,30){\line(1,0){25}}
\put(380,30){\line(1,0){10}}
\put(395,30){\line(1,0){10}}
\put(350,80){\line(1,0){60}}
\put(380,20){\line(1,0){60}}
\put(435,65){$\al$}
\put(360,70){$P_0$}
\put(420,25){$P_1$}
\put(355,35){$P_2$}
\put(390,5){$x_{6,0}$}
\put(390,65){$A_2$}
{\thicklines
\put(390,60){\line(-2,-1){20}}
\put(390,60){\line(2,-1){30}}
\put(390,60){\line(-1,-1){15}}
\put(390,60){\line(5,-4){25}}
\put(390,60){\line(1,-2){5}}
\put(397,45){\line(1,-2){5}}
\put(390,60){\line(-1,-2){5}}
\put(382,45){\line(-1,-2){5}} }

\end{picture}


\begin{picture}(450,100)(0,0)

\put(0,60){\vector(1,0){80}}
\put(30,60){\circle*{5}}
\put(-15,90){\line(1,-2){35}}
\put(45,90){\line(1,-2){35}}
\put(0,60){\line(-1,-2){15}}
\put(60,60){\line(-1,-2){5}}
\put(45,30){\line(1,2){5}}
\put(-15,30){\line(1,0){25}}
\put(20,30){\line(1,0){10}}
\put(35,30){\line(1,0){10}}
\put(-15,90){\line(1,0){60}}
\put(20,20){\line(1,0){60}}
\put(75,65){$\al$}
\put(0,75){$P_0$}
\put(60,25){$P_1$}
\put(-5,35){$P_2$}
\put(30,5){$x_{3,1}$}
{\thicklines
\put(30,60){\line(1,2){10}}
\put(30,60){\line(1,-2){5}}
\put(37,45){\line(1,-2){5}}
\put(30,60){\line(-1,-2){5}}
\put(23,45){\line(-1,-2){5}} }

\put(120,60){\vector(1,0){80}}
\put(150,60){\circle*{5}}
\put(105,90){\line(1,-2){35}}
\put(165,90){\line(1,-2){35}}
\put(120,60){\line(-1,-2){15}}
\put(180,60){\line(-1,-2){5}}
\put(165,30){\line(1,2){5}}
\put(105,30){\line(1,0){25}}
\put(140,30){\line(1,0){10}}
\put(155,30){\line(1,0){10}}
\put(105,90){\line(1,0){60}}
\put(140,20){\line(1,0){60}}
\put(195,65){$\al$}
\put(120,75){$P_0$}
\put(180,25){$P_1$}
\put(115,35){$P_2$}
\put(150,5){$x_{4,1}$}
{\thicklines
\put(150,60){\line(1,2){10}}
\put(150,60){\line(-1,2){10}}
\put(150,60){\line(1,-2){5}}
\put(157,45){\line(1,-2){5}}
\put(150,60){\line(-1,-2){5}}
\put(142,45){\line(-1,-2){5}} }

\put(240,60){\vector(1,0){80}}
\put(270,60){\circle*{5}}
\put(225,90){\line(1,-2){35}}
\put(285,90){\line(1,-2){35}}
\put(240,60){\line(-1,-2){15}}
\put(300,60){\line(-1,-2){5}}
\put(285,30){\line(1,2){5}}
\put(225,30){\line(1,0){25}}
\put(260,30){\line(1,0){10}}
\put(275,30){\line(1,0){10}}
\put(225,90){\line(1,0){60}}
\put(260,20){\line(1,0){60}}
\put(315,65){$\al$}
\put(240,75){$P_0$}
\put(300,25){$P_1$}
\put(270,5){$x_{5,1}$}
{\thicklines
\put(270,60){\line(1,2){10}}
\put(270,60){\line(1,-2){5}}
\put(277,45){\line(1,-2){5}}
\put(270,60){\line(2,-1){10}}
\put(283,53){\line(2,-1){10}}
\put(270,60){\line(-3,-2){15}}
\put(245,45){\line(-2,-1){10}}
\put(270,60){\line(-1,-2){5}}
\put(262,45){\line(-1,-2){5}} }

\put(360,60){\vector(1,0){80}}
\put(390,60){\circle*{5}}
\put(345,90){\line(1,-2){35}}
\put(405,90){\line(1,-2){35}}
\put(360,60){\line(-1,-2){15}}
\put(420,60){\line(-1,-2){5}}
\put(405,30){\line(1,2){5}}
\put(345,30){\line(1,0){25}}
\put(380,30){\line(1,0){10}}
\put(395,30){\line(1,0){10}}
\put(345,90){\line(1,0){60}}
\put(380,20){\line(1,0){60}}
\put(435,65){$\al$}
\put(360,75){$P_0$}
\put(420,25){$P_1$}
\put(390,5){$x_{6,1}$}
{\thicklines
\put(390,60){\line(1,2){10}}
\put(390,60){\line(-1,2){10}}
\put(390,60){\line(1,-2){5}}
\put(397,45){\line(1,-2){5}}
\put(390,60){\line(2,-1){10}}
\put(403,53){\line(2,-1){10}}
\put(390,60){\line(-3,-2){15}}
\put(365,45){\line(-2,-1){10}}
\put(390,60){\line(-1,-2){5}}
\put(382,45){\line(-1,-2){5}} }

\end{picture}


\begin{picture}(450,100)(0,0)

\put(0,60){\vector(1,0){80}}
\put(30,60){\circle*{5}}
\put(-15,90){\line(1,-2){35}}
\put(45,90){\line(1,-2){35}}
\put(0,60){\line(-1,-2){15}}
\put(60,60){\line(-1,-2){5}}
\put(45,30){\line(1,2){5}}
\put(-15,30){\line(1,0){25}}
\put(20,30){\line(1,0){10}}
\put(35,30){\line(1,0){10}}
\put(-15,90){\line(1,0){60}}
\put(20,20){\line(1,0){60}}
\put(75,65){$\al$}
\put(0,75){$P_0$}
\put(60,25){$P_1$}
\put(-5,35){$P_2$}
\put(30,5){$x_{3,2}$}
{\thicklines
\put(30,60){\line(1,2){10}}
\put(30,60){\line(-1,2){10}}
\put(30,60){\line(2,-1){30}} }

\put(120,60){\vector(1,0){80}}
\put(150,60){\circle*{5}}
\put(105,90){\line(1,-2){35}}
\put(165,90){\line(1,-2){35}}
\put(120,60){\line(-1,-2){15}}
\put(180,60){\line(-1,-2){5}}
\put(165,30){\line(1,2){5}}
\put(105,30){\line(1,0){25}}
\put(140,30){\line(1,0){10}}
\put(155,30){\line(1,0){10}}
\put(105,90){\line(1,0){60}}
\put(140,20){\line(1,0){60}}
\put(195,65){$\al$}
\put(120,75){$P_0$}
\put(180,25){$P_1$}
\put(115,35){$P_2$}
\put(150,5){$x_{4,2}$}
{\thicklines
\put(150,60){\line(1,2){10}}
\put(150,60){\line(-1,2){10}}
\put(150,60){\line(2,-1){30}}
\put(150,60){\line(-2,-1){20}} }

\put(240,60){\vector(1,0){80}}
\put(270,60){\circle*{5}}
\put(225,90){\line(1,-2){35}}
\put(285,90){\line(1,-2){35}}
\put(240,60){\line(-1,-2){15}}
\put(300,60){\line(-1,-2){5}}
\put(285,30){\line(1,2){5}}
\put(225,30){\line(1,0){25}}
\put(260,30){\line(1,0){10}}
\put(275,30){\line(1,0){10}}
\put(225,90){\line(1,0){60}}
\put(260,20){\line(1,0){60}}
\put(315,65){$\al$}
\put(235,75){$P_0$}
\put(300,25){$P_1$}
\put(235,35){$P_2$}
\put(270,5){$x_{5,2}$}
{\thicklines
\put(270,60){\line(1,2){10}}
\put(270,60){\line(2,1){20}}
\put(270,60){\line(-2,1){20}}
\put(270,60){\line(-1,2){10}}
\put(270,60){\line(2,-1){30}} }

\put(360,60){\vector(1,0){80}}
\put(390,60){\circle*{5}}
\put(345,90){\line(1,-2){35}}
\put(405,90){\line(1,-2){35}}
\put(360,60){\line(-1,-2){15}}
\put(420,60){\line(-1,-2){5}}
\put(405,30){\line(1,2){5}}
\put(345,30){\line(1,0){25}}
\put(380,30){\line(1,0){10}}
\put(395,30){\line(1,0){10}}
\put(345,90){\line(1,0){60}}
\put(380,20){\line(1,0){60}}
\put(435,65){$\al$}
\put(355,75){$P_0$}
\put(420,25){$P_1$}
\put(355,35){$P_2$}
\put(390,5){$x_{6,2}$}
{\thicklines
\put(390,60){\line(2,1){20}}
\put(390,60){\line(-2,1){20}}
\put(390,60){\line(-1,2){10}}
\put(390,60){\line(1,2){10}}
\put(390,60){\line(2,-1){30}}
\put(390,60){\line(-2,-1){20}} }

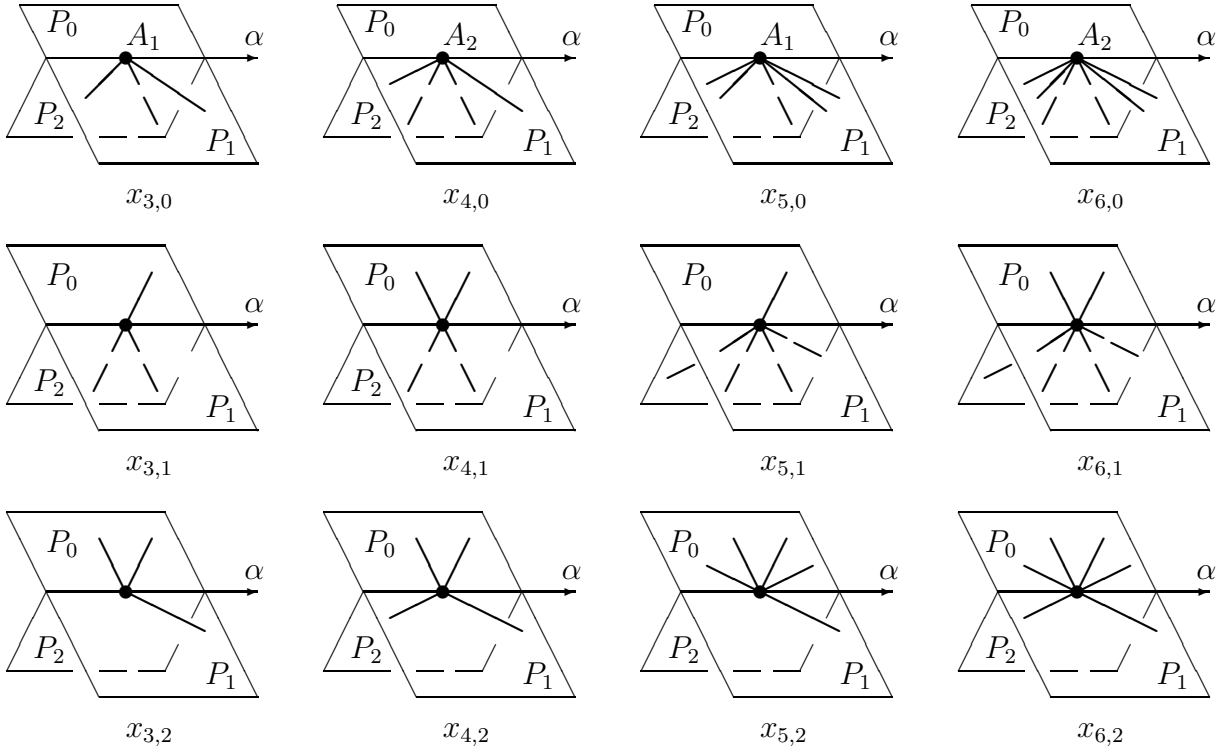
\put(70,-15){ {\bf Fig. 2b.} The letters for vertices of degrees 3, 4, 5, 6.}
\end{picture}
\vspace{3mm}

For example, the words $a_0 c_0$, $a_1 c_1$, $a_2 c_2$ encode the unknot.


\subsection{ Local isotopy moves in the three-page approach}

The relations (1)--(10) could be performed by rigid isotopy,
 which is denoted by "$\s$".


\begin{picture}(450,90)(0,0)

\put(0,60){\vector(1,0){110}}
\put(30,60){\circle*{3}}
\put(55,60){\circle*{3}}
\put(75,60){\circle*{3}}
\put(-10,80){\line(1,-2){30}}
\put(85,80){\line(1,-2){30}}
\put(0,60){\line(-1,-2){15}}
\put(95,60){\line(-1,-2){5}}
\put(-15,30){\line(1,0){25}}
\put(20,30){\line(1,0){15}}
\put(40,30){\line(1,0){10}}
\put(55,30){\line(1,0){10}}
\put(80,30){\line(-1,0){10}}
\put(80,30){\line(1,2){5}}
\put(-10,80){\line(1,0){95}}
\put(20,20){\line(1,0){95}}
\put(105,65){$\al$}
\put(0,70){$P_0$}
\put(-5,35){$P_2$}
\put(50,5){$d_0 d_1 d_2$}
{\thicklines
\put(30,60){\line(-2,-1){20}}
\put(30,60){\line(1,-2){5}}
\put(42,35){\line(-1,2){5}}
\put(42,35){\line(1,2){5}}
\put(50,50){\line(1,2){5}}
\put(65,75){\line(-2,-3){10}}
\put(65,75){\line(2,-3){10}}
\put(75,60){\line(1,-1){30}} }

\put(120,60){$\s$}

\put(150,60){\vector(1,0){55}}
\put(140,80){\line(1,-2){30}}
\put(180,80){\line(1,-2){30}}
\put(150,60){\line(-1,-2){15}}
\put(190,60){\line(-1,-2){5}}
\put(135,30){\line(1,0){25}}
\put(140,80){\line(1,0){40}}
\put(170,20){\line(1,0){40}}
\put(175,30){\line(-1,0){5}}
\put(175,30){\line(1,2){5}}
\put(200,65){$\al$}
\put(150,70){$P_0$}
\put(190,25){$P_1$}
\put(145,35){$P_2$}
\put(170,5){$1$}
{\thicklines
\put(165,45){\line(1,0){25}} }

\put(215,60){$\s$}

\put(245,60){\vector(1,0){75}}
\put(270,60){\circle*{3}}
\put(290,60){\circle*{3}}
\put(235,80){\line(1,-2){30}}
\put(295,80){\line(1,-2){30}}
\put(245,60){\line(-1,-2){15}}
\put(305,60){\line(-1,-2){5}}
\put(230,30){\line(1,0){25}}
\put(235,80){\line(1,0){60}}
\put(265,20){\line(1,0){60}}
\put(265,30){\line(1,0){10}}
\put(290,30){\line(-1,0){10}}
\put(290,30){\line(1,2){5}}
\put(315,65){$\al$}
\put(245,70){$P_0$}
\put(240,35){$P_2$}
\put(260,5){$b_2 d_2$}
{\thicklines
\put(260,40){\line(1,2){10}}
\put(280,75){\line(-2,-3){10}}
\put(280,75){\line(2,-3){10}}
\put(290,60){\line(1,-2){15}} }

\put(330,60){$\s$}

\put(355,60){\vector(1,0){80}}
\put(380,60){\circle*{3}}
\put(405,60){\circle*{3}}
\put(345,80){\line(1,-2){30}}
\put(410,80){\line(1,-2){30}}
\put(355,60){\line(-1,-2){15}}
\put(420,60){\line(-1,-2){5}}
\put(340,30){\line(1,0){25}}
\put(345,80){\line(1,0){65}}
\put(375,20){\line(1,0){65}}
\put(375,30){\line(1,0){10}}
\put(405,30){\line(-1,0){15}}
\put(405,30){\line(1,2){5}}
\put(430,65){$\al$}
\put(355,70){$P_0$}
\put(350,35){$P_2$}
\put(370,5){$d_0 b_0$}
{\thicklines
\put(380,60){\line(-1,-2){10}}
\put(380,60){\line(1,-2){5}}
\put(392,35){\line(-1,2){5}}
\put(392,35){\line(1,2){5}}
\put(405,60){\line(-1,-2){5}}
\put(405,60){\line(1,-2){15}} }

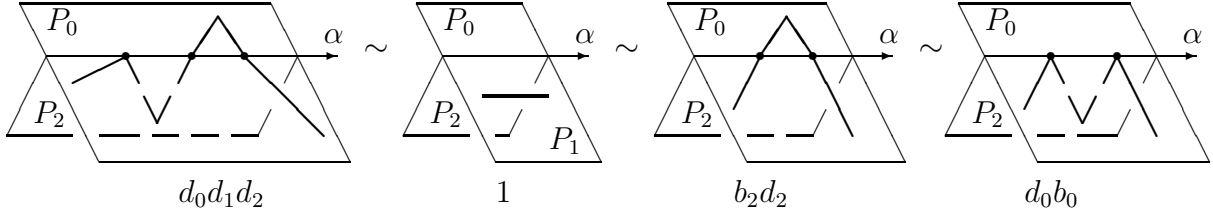
\put(30,-15){ {\bf Fig. 3a.}
The relations $(1)-(2)$ between invertible elements.}
\end{picture}
\vspace{4mm}


\begin{picture}(450,90)(0,0)

\put(0,60){\vector(1,0){70}}
\put(10,60){\circle*{3}}
\put(35,60){\circle*{3}}
\put(-10,80){\line(1,-2){30}}
\put(45,80){\line(1,-2){30}}
\put(0,60){\line(-1,-2){15}}
\put(55,60){\line(-1,-2){5}}
\put(-15,30){\line(1,0){25}}
\put(20,30){\line(1,0){15}}
\put(40,30){\line(-1,0){5}}
\put(40,30){\line(1,2){5}}
\put(-10,80){\line(1,0){55}}
\put(20,20){\line(1,0){55}}
\put(65,65){$\al$}
\put(0,70){$P_0$}
\put(55,25){$P_1$}
\put(-5,35){$P_2$}
\put(30,5){$a_0 d_1$}
{\thicklines
\put(10,60){\line(1,-2){5}}
\put(10,60){\line(3,-2){30}}
\put(22,35){\line(-1,2){5}}
\put(22,35){\line(1,2){5}}
\put(30,50){\line(1,2){10}}
}

\put(75,50){$\s a_2$;}

\put(125,60){\vector(1,0){55}}
\put(135,60){\circle*{3}}
\put(155,60){\circle*{3}}
\put(115,80){\line(1,-2){30}}
\put(155,80){\line(1,-2){30}}
\put(125,60){\line(-1,-2){15}}
\put(165,60){\line(-1,-2){5}}
\put(110,30){\line(1,0){25}}
\put(115,80){\line(1,0){40}}
\put(145,20){\line(1,0){40}}
\put(150,30){\line(-1,0){5}}
\put(150,30){\line(1,2){5}}
\put(175,65){$\al$}
\put(125,70){$P_0$}
\put(165,25){$P_1$}
\put(120,35){$P_2$}
\put(145,5){$a_2 c_1$}
{\thicklines
\put(135,60){\line(6,-5){30}}
\put(135,60){\line(2,3){10}}
\put(155,60){\line(-2,3){10}}
\put(155,60){\line(-1,-2){5}}
\put(147,45){\line(-1,-2){5}}
}

\put(185,50){$\s b_0$;}

\put(235,60){\vector(1,0){60}}
\put(255,60){\circle*{3}}
\put(270,60){\circle*{3}}
\put(225,80){\line(1,-2){30}}
\put(270,80){\line(1,-2){30}}
\put(235,60){\line(-1,-2){15}}
\put(280,60){\line(-1,-2){5}}
\put(220,30){\line(1,0){25}}
\put(225,80){\line(1,0){45}}
\put(255,20){\line(1,0){45}}
\put(265,30){\line(-1,0){10}}
\put(265,30){\line(1,2){5}}
\put(305,65){$\al$}
\put(235,70){$P_0$}
\put(280,25){$P_1$}
\put(230,35){$P_2$}
\put(250,5){$b_2 c_1$}
{\thicklines
\put(245,45){\line(2,3){20}}
\put(257,35){\line(1,2){5}}
\put(270,60){\line(-1,3){5}}
\put(270,60){\line(-1,-2){5}}
}

\put(300,50){$\s c_0$;}

\put(350,60){\vector(1,0){60}}
\put(360,60){\circle*{3}}
\put(385,60){\circle*{3}}
\put(340,80){\line(1,-2){30}}
\put(385,80){\line(1,-2){30}}
\put(350,60){\line(-1,-2){15}}
\put(395,60){\line(-1,-2){5}}
\put(335,30){\line(1,0){25}}
\put(340,80){\line(1,0){45}}
\put(370,20){\line(1,0){45}}
\put(380,30){\line(-1,0){10}}
\put(380,30){\line(1,2){5}}
\put(405,65){$\al$}
\put(350,70){$P_0$}
\put(395,25){$P_1$}
\put(345,35){$P_2$}
\put(365,5){$a_2 c_0$}
{\thicklines
\put(360,60){\line(1,1){15}}
\put(360,60){\line(1,-1){10}}
\put(385,60){\line(-3,-2){15}}
\put(385,60){\line(-1,-2){5}}
\put(372,35){\line(1,2){5}}
}

\put(415,50){$\s d_1$}

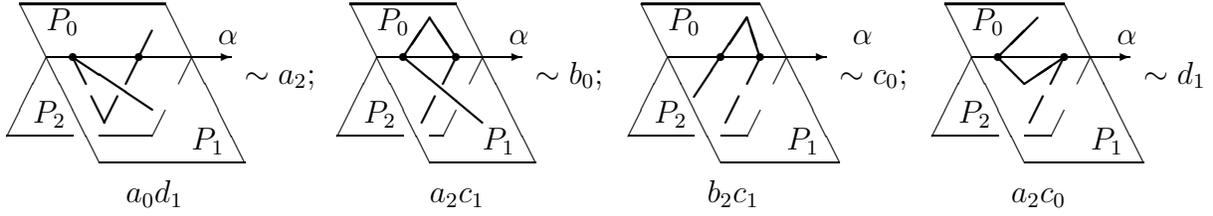
\put(30,-15){ {\bf Fig. 3b.}
The relations (3) are trivial moves at intersection points.}
\end{picture}
\vspace{3mm}
 

\begin{picture}(450,90)(0,5)

\put(0,60){\vector(1,0){85}}
\put(10,60){\circle*{3}}
\put(30,60){\circle*{5}}
\put(50,60){\circle*{3}}
\put(-5,80){\line(1,-4){15}}
\put(70,80){\line(1,-4){15}}
\put(0,60){\line(-1,-3){12}}
\put(75,60){\line(-1,-3){5}}
\put(-12,25){\line(1,0){15}}
\put(15,25){\line(1,0){5}}
\put(25,25){\line(1,0){5}}
\put(35,25){\line(1,0){5}}
\put(45,25){\line(1,0){5}}
\put(65,25){\line(-1,0){10}}
\put(65,25){\line(1,3){5}}
\put(-5,80){\line(1,0){75}}
\put(10,20){\line(1,0){75}}
\put(80,65){$\al$}
\put(20,5){$d_1 (x_{3,2} d_0)$}
{\thicklines
\put(0,40){\line(-1,-2){5}}
\put(10,60){\line(-1,-2){5}}
\put(10,60){\line(2,3){10}}
\put(30,60){\line(1,1){15}}
\put(30,60){\line(-2,3){10}}
\put(30,60){\line(1,-2){10}}
\put(50,60){\line(1,-2){5}}
\put(50,60){\line(-1,-2){10}}
\put(57,45){\line(1,-2){5}}
}

\put(90,45){$\s x_{3,1}$;}

\put(150,60){\vector(1,0){130}}
\put(160,60){\circle*{5}}
\put(180,60){\circle*{3}}
\put(190,60){\circle*{3}}
\put(200,60){\circle*{3}}
\put(210,60){\circle*{3}}
\put(225,60){\circle*{3}}
\put(250,60){\circle*{3}}
\put(260,60){\circle*{3}}
\put(145,80){\line(1,-4){15}}
\put(265,80){\line(1,-4){15}}
\put(150,60){\line(-1,-3){12}}
\put(270,60){\line(-1,-3){5}}
\put(138,25){\line(1,0){15}}
\put(165,25){\line(1,0){5}}
\put(175,25){\line(1,0){5}}
\put(185,25){\line(1,0){5}}
\put(195,25){\line(1,0){5}}
\put(205,25){\line(1,0){5}}
\put(215,25){\line(1,0){5}}
\put(225,25){\line(1,0){5}}
\put(235,25){\line(1,0){5}}
\put(245,25){\line(1,0){5}}
\put(260,25){\line(-1,0){5}}
\put(260,25){\line(1,3){5}}
\put(145,80){\line(1,0){120}}
\put(160,20){\line(1,0){120}}
\put(275,65){$\al$}
\put(170,5){$(x_{3,2} b_2) d_2^2 d_0^2 d_1^2$}
{\thicklines
\put(155,70){\line(1,-2){15}}
\put(160,60){\line(2,3){10}}
\put(170,40){\line(1,2){15}}
\put(170,75){\line(1,0){20}}
\put(185,70){\line(1,-2){15}}
\put(200,60){\line(-2,3){10}}
\put(200,60){\line(1,-2){5}}
\put(200,40){\line(1,0){15}}
\put(210,60){\line(-1,-2){5}}
\put(210,60){\line(1,-2){5}}
\put(220,40){\line(1,-2){5}}
\put(225,60){\line(-1,-2){10}}
\put(225,60){\line(1,-2){5}}
\put(225,30){\line(1,0){7}}
\put(237,35){\line(-1,2){5}}
\put(237,35){\line(1,2){5}}
\put(245,30){\line(-1,0){7}}
\put(245,30){\line(1,2){5}}
\put(245,50){\line(1,2){10}}
\put(255,50){\line(1,2){10}}
}

\put(290,50){$\s$}

\put(315,60){\vector(1,0){130}}
\put(330,60){\circle*{3}}
\put(350,60){\circle*{3}}
\put(380,60){\circle*{3}}
\put(400,60){\circle*{5}}
\put(420,60){\circle*{3}}
\put(310,80){\line(1,-4){15}}
\put(430,80){\line(1,-4){15}}
\put(315,60){\line(-1,-3){12}}
\put(435,60){\line(-1,-3){5}}
\put(303,25){\line(1,0){15}}
\put(330,25){\line(1,0){5}}
\put(340,25){\line(1,0){5}}
\put(350,25){\line(1,0){5}}
\put(360,25){\line(1,0){5}}
\put(370,25){\line(1,0){5}}
\put(380,25){\line(1,0){5}}
\put(390,25){\line(1,0){5}}
\put(400,25){\line(1,0){5}}
\put(410,25){\line(1,0){5}}
\put(425,25){\line(-1,0){5}}
\put(425,25){\line(1,3){5}}
\put(310,80){\line(1,0){120}}
\put(325,20){\line(1,0){120}}
\put(440,65){$\al$}
\put(335,5){$d_2 d_0 d_1 (x_{3,2} b_2)$}
{\thicklines
\put(330,60){\line(-2,3){10}}
\put(330,60){\line(1,-2){10}}
\put(350,60){\line(-1,-2){10}}
\put(350,60){\line(1,-2){5}}
\put(365,30){\line(-1,2){5}}
\put(365,30){\line(1,2){5}}
\put(380,60){\line(-1,-2){5}}
\put(380,60){\line(2,3){10}}
\put(400,60){\line(-2,3){10}}
\put(400,60){\line(1,-2){10}}
\put(400,60){\line(1,1){15}}
\put(410,40){\line(1,2){15}}
}

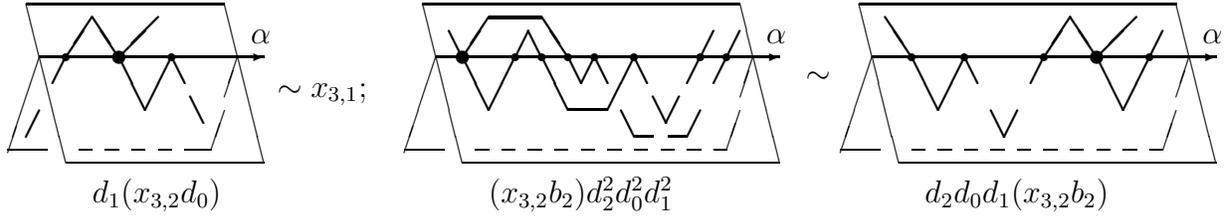
\put(30,-15){ {\bf Fig. 3c.}
The relations (4) and (9) are twistings of arcs at a $(2p-1)$-vertex.}
\end{picture}
\vspace{3mm}


\begin{picture}(450,90)(0,10)

\put(0,60){\vector(1,0){85}}
\put(10,60){\circle*{3}}
\put(30,60){\circle*{5}}
\put(50,60){\circle*{3}}
\put(-5,80){\line(1,-4){15}}
\put(70,80){\line(1,-4){15}}
\put(0,60){\line(-1,-3){12}}
\put(75,60){\line(-1,-3){5}}
\put(-12,25){\line(1,0){15}}
\put(15,25){\line(1,0){5}}
\put(25,25){\line(1,0){5}}
\put(35,25){\line(1,0){5}}
\put(45,25){\line(1,0){5}}
\put(65,25){\line(-1,0){10}}
\put(65,25){\line(1,3){5}}
\put(-5,80){\line(1,0){75}}
\put(10,20){\line(1,0){75}}
\put(80,65){$\al$}
\put(20,5){$b_0 x_{4,2} d_0$}
{\thicklines
\put(0,40){\line(-1,-2){5}}
\put(10,60){\line(-1,-2){5}}
\put(10,60){\line(1,-2){10}}
\put(30,60){\line(-1,1){15}}
\put(30,60){\line(1,1){15}}
\put(30,60){\line(-1,-2){10}}
\put(30,60){\line(1,-2){10}}
\put(50,60){\line(1,-2){5}}
\put(50,60){\line(-1,-2){10}}
\put(57,45){\line(1,-2){5}}
}

\put(90,45){$\s x_{4,1}$;}

\put(150,60){\vector(1,0){130}}
\put(160,60){\circle*{3}}
\put(170,60){\circle*{5}}
\put(180,60){\circle*{3}}
\put(190,60){\circle*{3}}
\put(200,60){\circle*{3}}
\put(210,60){\circle*{3}}
\put(225,60){\circle*{3}}
\put(250,60){\circle*{3}}
\put(260,60){\circle*{3}}
\put(145,80){\line(1,-4){15}}
\put(265,80){\line(1,-4){15}}
\put(150,60){\line(-1,-3){12}}
\put(270,60){\line(-1,-3){5}}
\put(138,25){\line(1,0){15}}
\put(165,25){\line(1,0){5}}
\put(175,25){\line(1,0){5}}
\put(185,25){\line(1,0){5}}
\put(195,25){\line(1,0){5}}
\put(205,25){\line(1,0){5}}
\put(215,25){\line(1,0){5}}
\put(225,25){\line(1,0){5}}
\put(235,25){\line(1,0){5}}
\put(245,25){\line(1,0){5}}
\put(260,25){\line(-1,0){5}}
\put(260,25){\line(1,3){5}}
\put(145,80){\line(1,0){120}}
\put(160,20){\line(1,0){120}}
\put(275,65){$\al$}
\put(160,5){$(d_2 x_{4,2} b_2) d_2^2 d_0^2 d_1^2$}
{\thicklines
\put(155,70){\line(1,-2){10}}
\put(165,70){\line(1,-2){10}}
\put(170,60){\line(-1,-2){5}}
\put(170,60){\line(2,3){10}}
\put(175,50){\line(1,2){10}}
\put(180,75){\line(1,0){10}}
\put(185,70){\line(1,-2){15}}
\put(200,60){\line(-2,3){10}}
\put(200,60){\line(1,-2){5}}
\put(200,40){\line(1,0){15}}
\put(210,60){\line(-1,-2){5}}
\put(210,60){\line(1,-2){5}}
\put(220,40){\line(1,-2){5}}
\put(225,60){\line(-1,-2){10}}
\put(225,60){\line(1,-2){5}}
\put(225,30){\line(1,0){7}}
\put(237,35){\line(-1,2){5}}
\put(237,35){\line(1,2){5}}
\put(245,30){\line(-1,0){7}}
\put(245,30){\line(1,2){5}}
\put(245,50){\line(1,2){10}}
\put(255,50){\line(1,2){10}}
}

\put(290,50){$\s$}

\put(315,60){\vector(1,0){130}}
\put(325,60){\circle*{3}}
\put(335,60){\circle*{3}}
\put(345,60){\circle*{3}}
\put(360,60){\circle*{3}}
\put(385,60){\circle*{3}}
\put(395,60){\circle*{3}}
\put(405,60){\circle*{3}}
\put(415,60){\circle*{5}}
\put(425,60){\circle*{3}}
\put(310,80){\line(1,-4){15}}
\put(430,80){\line(1,-4){15}}
\put(315,60){\line(-1,-3){12}}
\put(435,60){\line(-1,-3){5}}
\put(303,25){\line(1,0){15}}
\put(330,25){\line(1,0){5}}
\put(340,25){\line(1,0){5}}
\put(350,25){\line(1,0){5}}
\put(360,25){\line(1,0){5}}
\put(370,25){\line(1,0){5}}
\put(380,25){\line(1,0){5}}
\put(390,25){\line(1,0){5}}
\put(400,25){\line(1,0){5}}
\put(410,25){\line(1,0){5}}
\put(425,25){\line(-1,0){5}}
\put(425,25){\line(1,3){5}}
\put(310,80){\line(1,0){120}}
\put(325,20){\line(1,0){120}}
\put(440,65){$\al$}
\put(335,5){$d_2^2 d_0^2 d_1^2 (d_2 x_{4,2} b_2)$}
{\thicklines
\put(325,60){\line(-1,3){5}}
\put(325,60){\line(1,-2){10}}
\put(335,60){\line(-1,3){5}}
\put(335,60){\line(1,-2){5}}
\put(335,40){\line(1,0){15}}
\put(345,60){\line(-1,-2){5}}
\put(345,60){\line(1,-2){5}}
\put(360,60){\line(-1,-2){10}}
\put(360,60){\line(1,-2){5}}
\put(360,30){\line(-1,2){5}}
\put(360,30){\line(1,0){7}}
\put(372,35){\line(-1,2){5}}
\put(372,35){\line(1,2){5}}
\put(380,30){\line(-1,0){7}}
\put(380,30){\line(1,2){5}}
\put(385,60){\line(-1,-2){5}}
\put(385,60){\line(1,3){5}}
\put(390,75){\line(1,0){20}}
\put(390,50){\line(1,2){10}}
\put(400,70){\line(1,-2){10}}
\put(410,50){\line(1,2){10}}
\put(415,60){\line(-1,3){5}}
\put(415,60){\line(1,-2){5}}
\put(420,50){\line(1,2){10}}
}

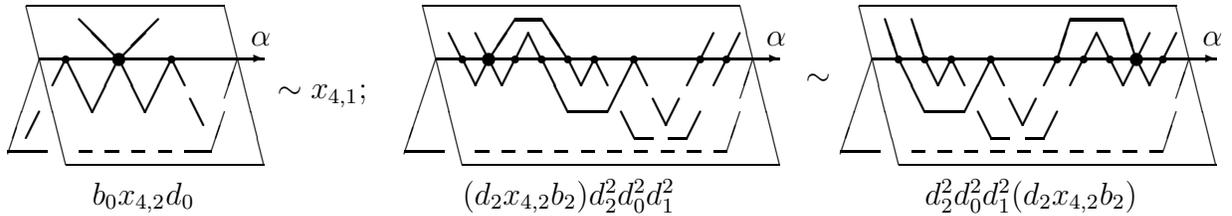
\put(30,-15){ {\bf Fig. 3d.}
The relations (4) and (10) are twistings of arcs at a $2q$-vertex.}
\end{picture}
\vspace{3mm}


\begin{picture}(450,85)(0,15)

\put(0,60){\vector(1,0){85}}
\put(25,60){\circle*{5}}
\put(45,60){\circle*{3}}
\put(-10,80){\line(1,-2){20}}
\put(60,80){\line(1,-2){20}}
\put(0,60){\line(-1,-3){5}}
\put(70,60){\line(-1,-3){5}}
\put(-5,45){\line(1,0){10}}
\put(15,45){\line(1,0){5}}
\put(25,45){\line(1,0){5}}
\put(40,45){\line(1,0){10}}
\put(65,45){\line(-1,0){5}}
\put(65,45){\line(1,3){5}}
\put(-10,80){\line(1,0){70}}
\put(10,40){\line(1,0){70}}
\put(80,65){$\al$}
\put(30,25){$x_{3,2} d_2$}
{\thicklines
\put(25,60){\line(-4,3){20}}
\put(25,60){\line(2,-3){10}}
\put(25,60){\line(2,3){10}}
\put(35,75){\line(2,-3){20}}
}

\put(90,50){$\s$}

\put(110,60){\vector(1,0){100}}
\put(120,60){\circle*{3}}
\put(140,60){\circle*{5}}
\put(160,60){\circle*{3}}
\put(180,60){\circle*{3}}
\put(100,80){\line(1,-2){20}}
\put(185,80){\line(1,-2){20}}
\put(110,60){\line(-1,-3){5}}
\put(195,60){\line(-1,-3){5}}
\put(105,45){\line(1,0){10}}
\put(120,45){\line(1,0){5}}
\put(135,45){\line(1,0){10}}
\put(155,45){\line(1,0){10}}
\put(175,45){\line(1,0){5}}
\put(190,45){\line(-1,0){5}}
\put(190,45){\line(1,3){5}}
\put(100,80){\line(1,0){85}}
\put(120,40){\line(1,0){85}}
\put(205,65){$\al$}
\put(120,25){$a_2 (x_{3,2} d_2) c_2$}
\put(215,50){;}
{\thicklines
\put(120,60){\line(2,-3){10}}
\put(120,60){\line(2,3){10}}
\put(130,75){\line(2,-3){20}}
\put(140,60){\line(2,3){10}}
\put(150,75){\line(2,-3){20}}
\put(180,60){\line(-2,-3){10}}
\put(180,60){\line(-2,3){10}}
}

\put(240,60){\vector(1,0){75}}
\put(255,60){\circle*{3}}
\put(270,60){\circle*{5}}
\put(285,60){\circle*{3}}
\put(230,80){\line(1,-2){20}}
\put(290,80){\line(1,-2){20}}
\put(240,60){\line(-1,-3){5}}
\put(300,60){\line(-1,-3){5}}
\put(235,45){\line(1,0){10}}
\put(250,45){\line(1,0){5}}
\put(267,45){\line(1,0){5}}
\put(275,45){\line(1,0){5}}
\put(295,45){\line(-1,0){10}}
\put(295,45){\line(1,3){5}}
\put(230,80){\line(1,0){60}}
\put(250,40){\line(1,0){60}}
\put(310,65){$\al$}
\put(260,25){$b_2 x_{3,2} b_2$}
{\thicklines
\put(250,50){\line(1,2){10}}
\put(260,70){\line(1,-1){20}}
\put(270,60){\line(2,3){10}}
\put(280,50){\line(1,2){10}}
}

\put(320,50){$\s$}

\put(345,60){\vector(1,0){95}}
\put(355,60){\circle*{3}}
\put(365,60){\circle*{3}}
\put(385,60){\circle*{5}}
\put(405,60){\circle*{3}}
\put(415,60){\circle*{3}}
\put(335,80){\line(1,-2){20}}
\put(415,80){\line(1,-2){20}}
\put(345,60){\line(-1,-3){5}}
\put(340,45){\line(1,0){10}}
\put(355,45){\line(1,0){10}}
\put(375,45){\line(1,0){5}}
\put(385,45){\line(1,0){5}}
\put(395,45){\line(1,0){5}}
\put(420,45){\line(-1,0){10}}
\put(420,45){\line(1,3){5}}
\put(335,80){\line(1,0){80}}
\put(355,40){\line(1,0){80}}
\put(435,65){$\al$}
\put(365,25){$a_2 (b_2 x_{3,2} b_2) c_2$}
{\thicklines
\put(355,60){\line(2,3){10}}
\put(355,60){\line(1,-2){5}}
\put(365,60){\line(-1,-2){5}}
\put(365,60){\line(1,1){10}}
\put(375,70){\line(1,-1){20}}
\put(385,60){\line(1,1){15}}
\put(405,60){\line(-1,-1){10}}
\put(405,60){\line(1,2){5}}
\put(415,60){\line(-1,2){5}}
\put(415,60){\line(-2,-3){10}}
}

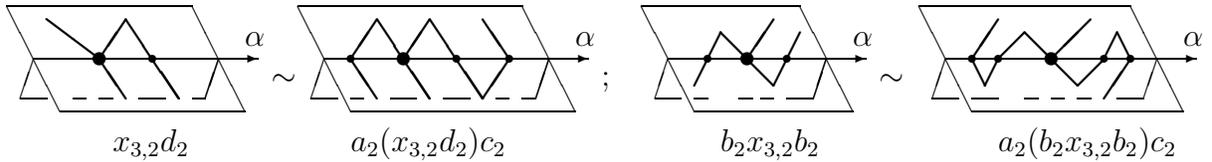
\put(30,5){ {\bf Fig. 3e.}
 Relations (5) mean rotatings of arcs at a $(2p-1)$-vertex.}
\end{picture}
\vspace{-2mm}


\begin{picture}(450,80)(0,25)

\put(0,60){\vector(1,0){85}}
\put(15,60){\circle*{3}}
\put(35,60){\circle*{5}}
\put(55,60){\circle*{3}}
\put(-10,80){\line(1,-2){20}}
\put(60,80){\line(1,-2){20}}
\put(0,60){\line(-1,-3){5}}
\put(70,60){\line(-1,-3){5}}
\put(-5,45){\line(1,0){10}}
\put(10,45){\line(1,0){5}}
\put(20,45){\line(1,0){5}}
\put(30,45){\line(1,0){5}}
\put(40,45){\line(1,0){5}}
\put(65,45){\line(-1,0){10}}
\put(65,45){\line(1,3){5}}
\put(-10,80){\line(1,0){70}}
\put(10,40){\line(1,0){70}}
\put(80,65){$\al$}
\put(20,25){$d_2 x_{4,2} d_2$}
{\thicklines
\put(0,75){\line(1,-1){25}}
\put(20,75){\line(1,-1){30}}
\put(25,50){\line(1,1){20}}
\put(45,70){\line(1,-1){10}}
\put(55,60){\line(1,-2){5}}
}

\put(95,50){$\s$}

\put(120,60){\vector(1,0){100}}
\put(130,60){\circle*{3}}
\put(145,60){\circle*{3}}
\put(165,60){\circle*{5}}
\put(185,60){\circle*{3}}
\put(195,60){\circle*{3}}
\put(110,80){\line(1,-2){20}}
\put(195,80){\line(1,-2){20}}
\put(120,60){\line(-1,-3){5}}
\put(205,60){\line(-1,-3){5}}
\put(115,45){\line(1,0){10}}
\put(130,45){\line(1,0){5}}
\put(145,45){\line(1,0){10}}
\put(160,45){\line(1,0){5}}
\put(170,45){\line(1,0){5}}
\put(185,45){\line(1,0){5}}
\put(200,45){\line(-1,0){5}}
\put(200,45){\line(1,3){5}}
\put(110,80){\line(1,0){85}}
\put(130,40){\line(1,0){85}}
\put(215,65){$\al$}
\put(225,60){;}
\put(130,25){$a_2 (d_2 x_{4,2} d_2) c_2$}
{\thicklines
\put(130,60){\line(2,-3){10}}
\put(130,60){\line(1,2){5}}
\put(135,70){\line(1,-1){20}}
\put(150,75){\line(1,-1){30}}
\put(155,50){\line(1,1){20}}
\put(185,60){\line(-1,1){10}}
\put(185,60){\line(1,-2){5}}
\put(195,60){\line(-1,-2){5}}
\put(195,60){\line(-2,3){10}}
}

\put(250,60){\vector(1,0){75}}
\put(265,60){\circle*{3}}
\put(280,60){\circle*{5}}
\put(295,60){\circle*{3}}
\put(240,80){\line(1,-2){20}}
\put(300,80){\line(1,-2){20}}
\put(250,60){\line(-1,-3){5}}
\put(310,60){\line(-1,-3){5}}
\put(245,45){\line(1,0){10}}
\put(260,45){\line(1,0){5}}
\put(277,45){\line(1,0){5}}
\put(285,45){\line(1,0){5}}
\put(305,45){\line(-1,0){10}}
\put(305,45){\line(1,3){5}}
\put(240,80){\line(1,0){60}}
\put(260,40){\line(1,0){60}}
\put(320,65){$\al$}
\put(270,25){$b_2 x_{4,2} b_2$}
{\thicklines
\put(260,50){\line(1,2){10}}
\put(270,70){\line(1,-1){20}}
\put(270,45){\line(2,3){20}}
\put(290,50){\line(1,2){10}}
}

\put(335,50){$\s$}

\put(360,60){\vector(1,0){95}}
\put(370,60){\circle*{3}}
\put(380,60){\circle*{3}}
\put(400,60){\circle*{5}}
\put(420,60){\circle*{3}}
\put(430,60){\circle*{3}}
\put(350,80){\line(1,-2){20}}
\put(430,80){\line(1,-2){20}}
\put(360,60){\line(-1,-3){5}}
\put(355,45){\line(1,0){10}}
\put(370,45){\line(1,0){10}}
\put(390,45){\line(1,0){5}}
\put(400,45){\line(1,0){5}}
\put(410,45){\line(1,0){5}}
\put(435,45){\line(-1,0){10}}
\put(435,45){\line(1,3){5}}
\put(350,80){\line(1,0){80}}
\put(370,40){\line(1,0){80}}
\put(450,65){$\al$}
\put(380,25){$a_2 (b_2 x_{4,2} b_2) c_2$}
{\thicklines
\put(370,60){\line(2,3){10}}
\put(370,60){\line(1,-2){5}}
\put(380,60){\line(-1,-2){5}}
\put(380,60){\line(1,1){10}}
\put(390,70){\line(1,-1){20}}
\put(385,45){\line(1,1){30}}
\put(420,60){\line(-1,-1){10}}
\put(420,60){\line(1,2){5}}
\put(430,60){\line(-1,2){5}}
\put(430,60){\line(-2,-3){10}}
}

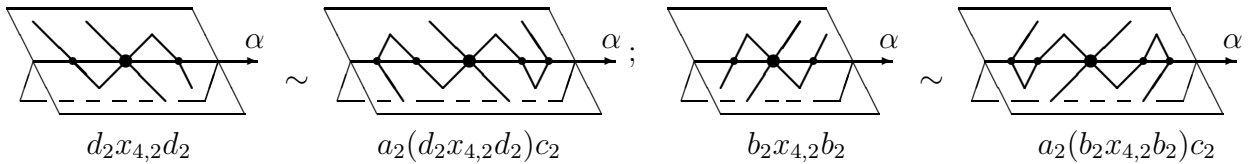
\put(30,5){ {\bf Fig. 3f.}
The relations (6) are rotatings of arcs at a $2q$-vertex.}
\end{picture}
\vspace{0mm}


\begin{picture}(450,90)(0,0)

\put(0,60){\vector(1,0){55}}
\put(10,60){\circle*{3}}
\put(30,60){\circle*{3}}
\put(-10,80){\line(1,-2){20}}
\put(30,80){\line(1,-2){20}}
\put(0,60){\line(-1,-3){5}}
\put(40,60){\line(-1,-3){5}}
\put(-5,45){\line(1,0){10}}
\put(10,45){\line(1,0){5}}
\put(35,45){\line(-1,0){10}}
\put(35,45){\line(1,3){5}}
\put(-10,80){\line(1,0){40}}
\put(10,40){\line(1,0){40}}
\put(40,65){$\al$}
\put(10,25){$d_2 c_2$}
{\thicklines
\put(0,75){\line(2,-3){20}}
\put(30,60){\line(-2,-3){10}}
\put(30,60){\line(-2,3){10}}
}

\put(60,50){,}

\put(80,60){\vector(1,0){70}}
\put(90,60){\circle*{3}}
\put(110,60){\circle*{3}}
\put(70,80){\line(1,-2){20}}
\put(125,80){\line(1,-2){20}}
\put(80,60){\line(-1,-3){5}}
\put(135,60){\line(-1,-3){5}}
\put(75,45){\line(1,0){10}}
\put(90,45){\line(1,0){5}}
\put(105,45){\line(1,0){5}}
\put(115,45){\line(1,0){5}}
\put(130,45){\line(-1,0){5}}
\put(130,45){\line(1,3){5}}
\put(70,80){\line(1,0){55}}
\put(90,40){\line(1,0){55}}
\put(145,65){$\al$}
\put(100,25){$a_2 b_2$}
{\thicklines
\put(90,60){\line(2,3){10}}
\put(90,60){\line(2,-3){10}}
\put(100,45){\line(2,3){20}}
}

\put(150,50){,}

\put(170,70){\vector(1,0){75}}
\put(182,70){\circle*{3}}
\put(195,70){\circle*{3}}
\put(203,70){\circle*{3}}
\put(215,70){\circle*{3}}
\put(165,85){\line(1,-3){15}}
\put(225,85){\line(1,-3){15}}
\put(170,70){\line(-1,-5){5}}
\put(230,70){\line(-1,-5){5}}
\put(165,45){\line(1,0){10}}
\put(180,45){\line(1,0){5}}
\put(190,45){\line(1,0){5}}
\put(200,45){\line(1,0){5}}
\put(210,45){\line(1,0){5}}
\put(225,45){\line(-1,0){5}}
\put(225,45){\line(1,5){5}}
\put(165,85){\line(1,0){60}}
\put(180,40){\line(1,0){60}}
\put(240,75){$\al$}
\put(180,25){$b_1 d_2 d_1 b_2$}
{\thicklines
\put(185,65){\line(-1,2){8}}
\put(192,50){\line(-1,2){5}}
\put(192,50){\line(1,2){5}}
\put(200,65){\line(1,2){8}}
\put(205,50){\line(-1,2){15}}
\put(205,50){\line(1,2){15}}
}

\put(250,60){,}

\put(270,60){\vector(1,0){70}}
\put(280,60){\circle*{5}}
\put(300,60){\circle*{3}}
\put(260,80){\line(1,-2){20}}
\put(315,80){\line(1,-2){20}}
\put(270,60){\line(-1,-3){5}}
\put(325,60){\line(-1,-3){5}}
\put(265,45){\line(1,0){10}}
\put(280,45){\line(1,0){5}}
\put(295,45){\line(1,0){5}}
\put(305,45){\line(1,0){5}}
\put(320,45){\line(-1,0){5}}
\put(320,45){\line(1,3){5}}
\put(260,80){\line(1,0){55}}
\put(280,40){\line(1,0){55}}
\put(335,65){$\al$}
\put(290,25){$x_{3,2} b_2$}
{\thicklines
\put(280,60){\line(2,3){10}}
\put(290,45){\line(-2,3){20}}
\put(290,45){\line(2,3){20}}
}

\put(340,50){,}

\put(360,60){\vector(1,0){90}}
\put(370,60){\circle*{3}}
\put(390,60){\circle*{5}}
\put(410,60){\circle*{3}}
\put(350,80){\line(1,-2){20}}
\put(425,80){\line(1,-2){20}}
\put(360,60){\line(-1,-3){5}}
\put(435,60){\line(-1,-3){5}}
\put(355,45){\line(1,0){10}}
\put(370,45){\line(1,0){5}}
\put(385,45){\line(1,0){10}}
\put(405,45){\line(1,0){5}}
\put(415,45){\line(1,0){5}}
\put(430,45){\line(-1,0){5}}
\put(430,45){\line(1,3){5}}
\put(350,80){\line(1,0){75}}
\put(370,40){\line(1,0){75}}
\put(445,65){$\al$}
\put(380,25){$d_2 x_{4,2} b_2$}
{\thicklines
\put(380,45){\line(-2,3){20}}
\put(380,45){\line(2,3){20}}
\put(400,45){\line(-2,3){20}}
\put(400,45){\line(2,3){20}}
}

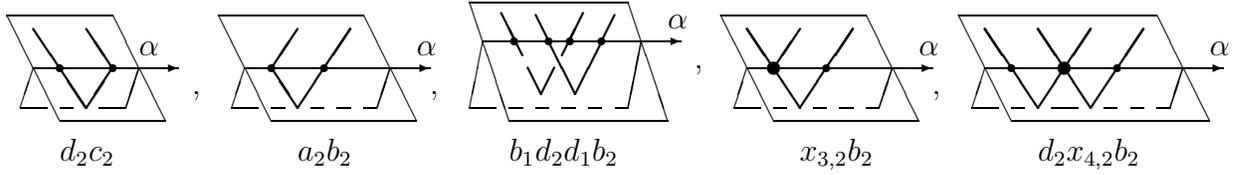
\put(30,5){ {\bf Fig. 3g.}
 These elements commute with
  $a_0, b_0, c_0, b_2 d_0 d_2, x_{m,0}$ in (7)--(8).}
\end{picture}
\vspace{-2mm}

During rigid isotopies from Fig.~3 neighbourhoods of vertices
 lie inside 2 pages.


\subsection{ Scheme for the proof of Theorem~1}
The formal definition of a three-page embedding into the book $\Y$
 is given in Subsection~3.1.
In Subsection~3.2, a three-page embedding $G\subset\Y$
 of a spatial graph $G\subset\R^3$ will be construsted
 from any plane diagram of $G$.
Theorem~1a will be proved in Subsection 3.3 by exploiting an encoding
 of the constructed three-page embedding $G\subset\Y$
 by a word in the alphabet $\A_n$.

The hard part of Theorem~1b is that
 any rigid isotopy of spatial $n$-graphs decomposes on
 the relations $(1)-(10)$.
For the proof of Theorem~1b, the notions of
 graph tangles and three-page tangles
 are introduced in Section~4.
The semigroup $RGT_n$ of all rigid graph tangles will be described
 by generators and the relations (11)--(23)
 in Lemma~1 in Subsection~4.2.
Every three-page embedding of a graph could be transformed to
 an almost balanced tangle, which is
 a "special" three-page tangle.
The semigroup $RBT_n$ of almost balanced tangles
 is turned out to be isomorphic to $RGT_n$, see Lemma~2 in Subsection~5.1.
Under an isomorphism $\ph:RGT_n\to RBT_n$,
 the relations $(11)-(23)$ convert to
 the relations $\ph(11)-\ph(23)$ between words in the alphabet $\A_n$.

Any three-page embedding of a graph could be represented by
 a graph tangle of $RGT_n$.
By Lemma~1 any isotopy between three-page embeddings of graphs
 decomposes on the isotopies $(11)-(23)$ between graph tangles, and
 hence on the isotopies $\ph(11)-\ph(23)$ between almost balanced tangles.
So, it remains to deduce 
 the relations $\ph(11)-\ph(23)$ of $RBT_n$ from
 the relations $(1)-(10)$ of $RSG_n$,
 see Lemma~3 in Subsection~5.2.
Lemma~3 will be checked in Section~6 
 by exploiting technical claims.
Theorem~1c is proved in Subsection~5.3
 by using knot-like three-page tangles.


\setcounter{section}{2}
\section {Three-page embeddings of spatial graphs}

\subsection{ The formal definition of a three-page embedding.}

Let $G$ be a finite graph, $A\in G$ be its point.
Any small segment $\ga\subset G$ with the endpoint $A\in G$
 is said te be \emph{an arc} of $G$.
Hence there are exactly $k$ arcs at each $k$-vertex of $G$.
\smallskip

{\bf Definition 7} (\emph{a three-page embedding}).
Suppose that a spatial graph $G\subset\R^3$ is contained in
 the three-page book $\Y\subset\R^3$.
The embedding $G\subset\Y$ is called \emph{a three-page embedding},
 if the following conditions hold (see Fig.~4a):
\smallskip

(7.1) all vertices of the graph $G$ lie in the axis $\alpha$;

(7.2) the intersection
  $G \cap \alpha = A_1\cup \dots \cup A_k$ is
  a non-empty finite set of points;

(7.3) two arcs with an endpoint $A_l\in G\cap \al$
 that is not a vertex of $G$ lie in different pages $P_i,P_j$ ($i\neq j$);

(7.4) \emph{balance:} neighbourhoods of vertices $A_l\in G\cap \al$
  look like pictures of Fig.~2b.
\smallskip

Since arcs of an embedding $G\subset\Y$, lying in a page $P_i$,
 are not intersected, then by isotopy inside $\Y$
 we may secure the following condition, which will be always assumed:

(7.5) \emph{monotone}: for each $i\in \Z_3$, the restriction of
  the orthogonal projection $\Y\to \alpha \approx \R$ to
  each connected component of the intersection $G\cap P_i$ is
  a monotone function.                   
\ed


\begin{picture}(480,110)(0,5)

\put(10,90){\vector(1,-1){10}}
\put(20,80){\line(1,-1){10}}
\put(50,70){\vector(0,-1){20}}
\put(50,50){\line(0,-1){20}}
\put(70,30){\vector(1,0){20}}
\put(90,30){\line(0,1){20}}
\put(80,60){\line(1,-1){10}}
\put(80,80){\vector(1,1){20}}
\put(100,80){$\alpha$}
\put(60,30){\circle*{5}}

\put(110,30){$D$}
\put(55,15){$B$}
\put(35,25){$L_B$}
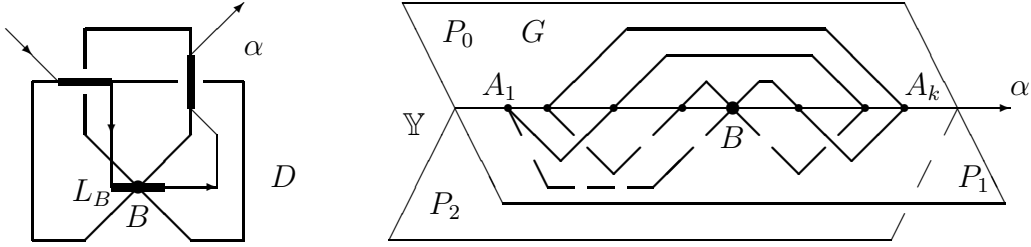
\put(20,-10){ {\bf Fig.~4a.}
 The spatial graph $G$  is encoded 
  by the word $w_G = a_0 a_1 b_2 d_1 x_{4,1} d_2 c_1 c_2$.}
\put(160,50){$\Y$}
\put(190,65){$A_1$}
\put(350,65){$A_k$}
\put(175,85){$P_0$}
\put(205,85){$G$}
\put(280,45){$B$}

\put(180,60){\vector(1,0){210}}
\put(390,65){$\alpha$}
\put(155,10){\line(1,0){190}}
\put(160,100){\line(1,0){190}}
\put(198,24){\line(1,0){190}}
\put(160,100){\line(1,-2){38}}
\put(350,100){\line(1,-2){38}}
\put(155,10){\line(1,2){25}}
\put(370,60){\line(-1,-2){5}}
\put(360,40){\line(-1,-2){5}}
\put(350,20){\line(-1,-2){5}}
\put(170,20){$P_2$}
\put(370,30){$P_1$}

\thicklines

\put(50,70){\line(1,0){25}}
\put(20,70){\line(1,0){55}}
\put(85,70){\line(1,0){15}}
\put(20,10){\line(1,0){20}}
\put(40,10){\line(1,1){40}}
\put(40,90){\line(1,0){40}}
\put(20,10){\line(0,1){60}}
\put(40,50){\line(0,1){15}}
\put(40,50){\line(1,-1){40}}
\put(40,75){\line(0,1){15}}
\put(80,10){\line(1,0){20}}
\put(80,50){\line(0,1){40}}
\put(100,10){\line(0,1){60}}

\put(200,60){\line(1,-1){20}}
\put(200,60){\line(1,-2){5}}
\put(210,40){\line(1,-2){5}}
\put(215,30){\line(1,0){10}}
\put(215,60){\line(1,-1){10}}
\put(215,60){\line(1,1){30}}
\put(220,40){\line(1,1){40}}
\put(230,30){\line(1,0){10}}
\put(240,35){\line(-1,1){10}}
\put(240,35){\line(1,1){10}}
\put(245,90){\line(1,0){75}}
\put(255,50){\line(1,1){20}}
\put(255,30){\line(-1,0){10}}
\put(255,30){\line(1,1){15}}
\put(260,80){\line(1,0){55}}
\put(275,50){\line(1,1){20}}
\put(275,70){\line(1,-1){20}}
\put(295,70){\line(1,0){5}}
\put(300,70){\line(1,-1){30}}
\put(310,35){\line(-1,1){10}}
\put(310,35){\line(1,1){10}}
\put(335,60){\line(-1,-1){10}}
\put(335,60){\line(-1,1){20}}
\put(350,60){\line(-1,-1){20}}
\put(350,60){\line(-1,1){30}}

\linethickness{1mm}

\put(30,70){\line(1,0){20}}
\put(80,60){\line(0,1){20}}
\put(50,30){\line(1,0){20}}

\put(200,60){\circle*{3}}
\put(215,60){\circle*{3}}
\put(240,60){\circle*{3}}
\put(266,60){\circle*{3}}
\put(285,60){\circle*{5}}
\put(310,60){\circle*{3}}
\put(335,60){\circle*{3}}
\put(350,60){\circle*{3}}

\end{picture}
\vspace{3mm}


\subsection{Construction of a three-page embedding from a plane diagram}

Let $D$ be \emph{a plane diagram} of a spatial graph $G$.
Namely, $D\subset\R^2$ is a plane graph with vertices of two types:
 ones correspond to initial vertices of $G\subset\R^3$
 and the others denote usual \emph{crossings} 
 in a planar representation of the spatial graph $G\subset\R^3$.
\smallskip

{\bf Definition~8} (\emph{bridges, upper and lower arcs}).
Let us choose bridges and arcs in $D$.
\smallskip

(a) For each crossing of the plane diagram $D$, let us mark out a small arc
     (\emph{a regular bridge}) in the overcrossing arc.
     See the left pictures of Fig.~4a and 4b.

(b) For any $2q$-vertex $B\in G$,
     let us call two adjacent arcs at $B$ \emph{upper arcs}, 
     and call the other $2q-2$ arcs at $B$ \emph{lower arcs}. 
    Then mark out a small segment (\emph{a singular bridge $L_B$})
     containing $B$ such that the upper arcs are separated from
     the lower arcs by the bridge $L_B$  in some neighbourhood of $B$.
    See the left picture of Fig.~4a.

(c) For every $(2p-1)$-vertex $C\in G$,
     let us call one arc at $C$ \emph{an upper arc},
     and call the other $2p-2$ arcs at $C$ \emph{lower arcs}.
    Then mark out a small segment (\emph{a singular bridge $L_C$}) 
     containing $C$ such that the upper arc is separated from
     the lower arcs by the bridge $L_C$ in some neighbourhood of $C$.
    See the left picture of Fig.~4b.
\ed
\smallskip


\begin{picture}(480,120)(0,-5)

\put(-10,80){\vector(1,0){10}}
\put(0,80){\line(1,0){20}}
\put(40,80){\vector(1,0){10}}
\put(50,80){\line(1,0){20}}
\put(70,60){\vector(-1,0){10}}
\put(60,60){\line(-1,0){10}}
\put(50,40){\vector(1,0){20}}
\put(70,40){\line(1,0){10}}
\put(100,40){\vector(-1,1){20}}
\put(80,60){\vector(1,0){50}}
\put(30,80){\circle*{5}}
\put(90,40){\circle*{5}}
\put(125,65){$\alpha$}

\put(75,45){$L_C$}
\put(90,20){$C$}
\put(50,10){$D$}
\put(160,50){$\Y$}
\put(330,65){$C$}
\put(185,65){$A_1$}
\put(375,65){$A_k$}
\put(230,75){$G$}
\put(170,10){$P_2$}
\put(395,30){$P_1$}
\put(175,85){$P_0$}
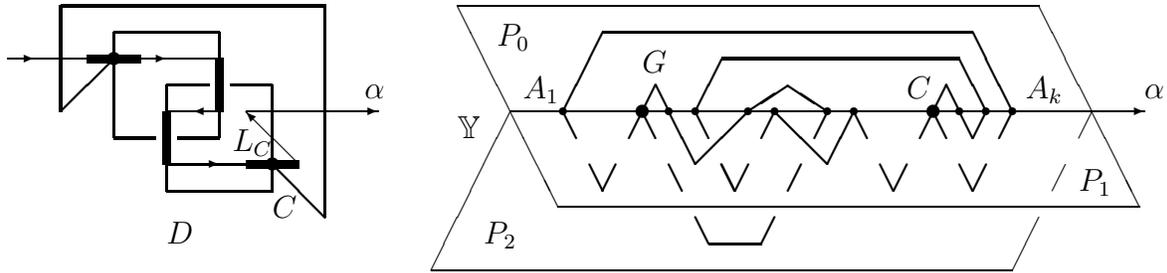
\put(10,-20){ {\bf Fig.~4b.}
 The spatial graph $G$ is encoded 
  by $w_G=a_1 x_{3,1} d_2 a_1 b_2 b_0 c_1 d_0 x_{3,1} d_2 b_2 c_1$}

\put(180,60){\vector(1,0){240}}
\put(420,65){$\alpha$}
\put(150,0){\line(1,0){220}}
\put(160,100){\line(1,0){220}}
\put(198,24){\line(1,0){220}}
\put(160,100){\line(1,-2){38}}
\put(380,100){\line(1,-2){38}}
\put(150,0){\line(1,2){30}}
\put(400,60){\line(-1,-2){5}}
\put(390,40){\line(-1,-2){5}}
\put(380,20){\line(-1,-2){10}}

\thicklines

\put(10,60){\line(0,1){40}}
\put(30,80){\line(0,1){10}}
\put(30,80){\line(-1,-1){20}}
\put(30,90){\line(1,0){40}}
\put(30,50){\line(0,1){30}}
\put(30,50){\line(1,0){16}}
\put(54,50){\line(1,0){16}}
\put(70,50){\line(0,1){10}}
\put(50,60){\line(0,1){10}}
\put(50,70){\line(1,0){16}}
\put(74,70){\line(1,0){16}}
\put(90,40){\line(0,1){30}}
\put(90,40){\line(1,-1){20}}
\put(10,100){\line(1,0){100}}
\put(110,20){\line(0,1){80}}
\put(50,30){\line(1,0){40}}
\put(70,80){\line(0,1){10}}
\put(50,30){\line(0,1){10}}
\put(90,30){\line(0,1){10}}

\put(200,60){\line(1,2){15}}
\put(200,60){\line(1,-2){5}}
\put(210,40){\line(1,-2){5}}
\put(215,90){\line(1,0){140}}
\put(215,30){\line(1,2){5}}
\put(230,60){\line(-1,-2){5}}
\put(230,60){\line(1,-2){5}}
\put(230,60){\line(1,2){5}}
\put(235,70){\line(1,-2){15}}
\put(240,40){\line(1,-2){5}}
\put(250,60){\line(1,2){10}}
\put(250,60){\line(1,-2){5}}
\put(255,10){\line(-1,2){5}}
\put(255,10){\line(1,0){20}}
\put(260,40){\line(1,-2){5}}
\put(260,80){\line(1,0){90}}
\put(265,30){\line(-1,2){5}}
\put(265,30){\line(1,2){5}}
\put(270,60){\line(-1,-1){20}}
\put(270,60){\line(3,2){15}}
\put(275,10){\line(1,2){5}}
\put(280,60){\line(1,-1){20}}
\put(280,60){\line(-1,-2){5}}
\put(290,40){\line(-1,-2){5}}
\put(300,60){\line(-3,2){15}}
\put(300,60){\line(-1,-2){5}}
\put(310,60){\line(-1,-2){10}}
\put(310,60){\line(1,-2){5}}
\put(325,30){\line(-1,2){5}}
\put(325,30){\line(1,2){5}}
\put(340,60){\line(-1,-2){5}}
\put(340,60){\line(1,-2){5}}
\put(340,60){\line(1,2){5}}
\put(345,70){\line(1,-2){10}}
\put(355,30){\line(-1,2){5}}
\put(355,30){\line(1,2){5}}
\put(360,60){\line(-1,-2){5}}
\put(360,60){\line(-1,2){10}}
\put(370,60){\line(-1,-2){5}}
\put(370,60){\line(-1,2){15}}

\linethickness{1mm}
\put(20,80){\line(1,0){20}}
\put(50,40){\line(0,1){20}}
\put(70,60){\line(0,1){20}}
\put(80,40){\line(1,0){20}}

\put(200,60){\circle*{3}}
\put(230,60){\circle*{5}}
\put(240,60){\circle*{3}}
\put(250,60){\circle*{3}}
\put(270,60){\circle*{3}}
\put(280,60){\circle*{3}}
\put(300,60){\circle*{3}}
\put(310,60){\circle*{3}}
\put(340,60){\circle*{5}}
\put(350,60){\circle*{3}}
\put(360,60){\circle*{3}}
\put(370,60){\circle*{3}}

\end{picture}
\vspace{5mm}

Let us take a non-self-intersected oriented path $\alpha$ 
 in the plane of $D$ so that
\smallskip

(1) The ends of the path $\al$ lie far from the diagram $D$.

(2) The path $\alpha$ goes through each bridge only once.

(3) The path $\alpha$ intersects the rest part of $D$
    transversally.

(4) For every vertex $A\in G$,
    the upper arcs at $A$ lie to the left of the oriented path $\al$.
   All the lower arcs at $A$ lie to the right of $\al$. 
   See the left pictures of Fig.~4a and 4b.
 
(5) For any $2q$-vertex $B\in G$,
    an upper arc and exactly $q-1$ lower arcs at $B$
    meet the oriented path $\al$ before $B$.
   Similarly, another lower arc and the other $q-1$ arcs at $B$
    meet $\al$ after $B$.
   See the left picture of Fig.~4a.
 
(6) For each $(2p-1)$-vertex $C$,
    the upper arc and $p-1$ lower arcs at $C$ 
    meet $\al$ after $C$.
   Analogously, the other $p-1$ arcs at $C$ 
    meet $\al$ before $C$.
   See the left picture of Fig.~4b. 
\smallskip

Such a path $\alpha$ could be found as follows.
Let us consider only bridges in $\R^2$, i.e. finitely many arcs.
Pass any path $\alpha$ through these bridges
 according to Conditions (1), (2), and (4).
Condition~(3) will hold, if we shall move $\al$
 in a general position with respect to the diagram $D$.
Let us assume that Condition~(6) does not hold for a $(2p-1)$-vertex $C$.
For instance, let the upper arc $\ga$ at $C$ meet $\al$ before $C$.
Then by using "Reidemeister's move $R2$" (see Fig.~1)
 create two additional transversal intersections
 of $\ga$ and $\alpha$ near the bridge $L_C$.
For example, see the left puicture of Fig.~4b.
Similarly, we deal with Condition~(5).
Deform the plane of $D$ in such a way that
 $\alpha$ becomes\footnote{
 We extend the path $\al$ to an infinite straight line
  $\al\ap\R\subset\R^2$.}
 a straight line and
\smallskip

(7) the restriction of the orthogonal projection
 $\R^2\to \alpha \approx \R$ 
 to each connected component of $D-\alpha$ 
 is a monotone function.
\smallskip

By $P_0$ (respectively, $P_2$) denote the upper part
 (respectively, the lower part) of $\R^2-\alpha$.
See the right pictures of Fig.~4a and 4b.
Let us attach the third page $P_1$ at the axis $\alpha$ and
 push out all the regular bridges into $P_1$ such that 
 each regular bridge becomes \emph{a trivial arc} (a union of two segments).
Conditions~(7.1)--(7.3) of Definition~7 hold by the construction.
Conditions~(4)--(7) of this subsection
 imply Condition~(7.4).


\subsection{ Proof of Theorem 1a.}
Let us take a plane diagram $D$ of a given spatial graph $G$.
Starting with the diagram $D$ one can construct 
 a three-page embedding $G\subset\Y$
 as in Subsection 3.2.
Any three-page embedding is uniquely determined by
 its small part near the axis $\alpha$.
Indeed, in order to reconstruct 
 the whole embedding it is sufficient
 to join all opposite-directed arcs
 in each page starting with interior arcs.
Due to Conditions~(7.3)--(7.5) of Definition~7,
 only the pictures of Fig.~2a and 2b
 could be involved in a three-page embedding near the axis $\alpha$.
For the three-page embedding $G\subset \Y$,
 let us write one by one letters
 of the alphabet $\A_n$, corresponding to
 the intersection points of $G\cap \alpha$.
One gets a word $w_G \in W_n$ (by $W_n$ we denote all words in $\A_n$
 including the empty word $\es$), see the right pictures of Fig.~4a and 4b.
Finally, consider the word $w_G$ as 
 an element of $RSG_n$.
Note that one can rotate any three-page embedding
 around the axis $\al$.
Then any spatial graph $G$ is represented by
 three words obtained from $w_G$ by the index shift $i\mapsto i+1$.
\qed
\smallskip

{\bf Proof of Corollary 5.}
Let us draw a given graph $G$ (possibly with hanging edges)
 in the plane $\R^2$ such that its edges are intersected in
 double points only.
Near each double point push exactly one arc out of $\R^2$.
For the obtained spatial graph $G\subset\R^3$,
 let us construct a three-page embedding $G\subset\Y$ as
 in Subsection~3.2.
\qed


\subsection{ Balanced words in the alphabet $\A_n$.}

By the encoding procedure of the proof of Theorem~1a,
 we may not obtain all words of $W_n$.
A word $w\in W_n$ is called \emph{balanced}, 
 if it encodes some three-page embedding 
 of a spatial graph.
There exists the following geometric criterion for a balanced word:
 "in each page $P_i$ all arcs have to be joined with each other".
Arcs of a non-balanced three-page embedding 
 could go to infinity not meeting each other.
One can rewrite this criterion algebraically
 via the alphabet $\A_n$.
\smallskip

{\bf Definition 9} (\emph{balanced bracket expressions}).
Any expression $\be$, consisting of left and right round brackets,
 is said to be \emph{bracket}.
A bracket expression $\be$ is called \emph{balanced},
 if (by reading $\be$ from left to right) in each place the number of
 the left brackets is not less than 
 the number of the right ones, and 
 their total numbers are equal.
\ed
\smallskip

{\bf Definition 10} (\emph{balanced words}).
Let us consider the following substitution\footnote{
  As usual, we have
  $i\in\Z_3$, $3\leq m\leq n$, $2\leq p\leq \frac{n+1}{2}$, and
  $2\leq q\leq \frac{n}{2}$.}:
 $$\left\lbrace \begin{array}{lll}
  a_i,\; b_i, \; c_i, \; d_i, \; x_{m,i} \to \varnothing,  &
  a_{i\pm 1}, \; b_{i-1}, \; d_{i+1}, \; x_{2p-1,i-1} \to (;  &
  b_{i+1}, \; c_{i\pm 1}, \; d_{i-1} \to ); \\
  x_{2q,i+1} \to )(;  \quad x_{2q,i-1}\to )^{q-1} (^{q-1}; &
  x_{2p-1,i+1} \to )^{p-1} (^{p-1},
  \mbox{ where } & (^j = j \mbox{ brackets } (.
 \end{array} \right.$$
Suppose that after the above substitution, for a word $w\in W_n$,
 we get a balanced bracket expression $\be_i(w)$.
In this case, the word $w$ will be called \emph{$i$-balanced}.
A word $w$ is said to be \emph{balanced},
 if it is $i$-balanced for each $i\in \Z_3$.
\ed
\smallskip

So, a word $w$ is balanced if and only if 
 all three bracket expressions $\be_i(w)$ are balanced.
For example, for the word of Fig.~4a,
 we get the following bracket expressions:
 $\be_0(w)=(( () () ))$,
 $\be_1(w)=() ()$,
 $\be_2(w)=(( )) ()$.
By $W_{n,i}$ let us denote the set of 
 all $i$-balanced words in $\A_n$.
Definition~10 implies that there is an algorithm to decide,
 whether a word $w$ is balanced, which is linear
 in the length of $w$.
Put $BW_{n}=W_{n,0}\cap W_{n,1}\cap W_{n,2}\subset W_n$.
This is the set of all balanced words in $\A_n$.
We shall prove in Lemma~4 that the set $BW_n$ encodes the centers of
 the semigroups $RSG_n$ and $NSG_n$, see Subsection~5.3.


\setcounter{section}{3}
\section{ Graph tangles and three-page tangles}

\subsection{ Graph tangles.}

The initial category of tangles was studied by V.~Turaev \cite{Tu}.
Let us take two horizontal semilines given by coordinates:
 $(r, 0, 0)$ and $(r, 0, 1)$, $r\in \R_+$.
For all $k\in \N$, let us mark out the points $(k,0,0)$, $(k,0,1)$
 in the both semilines.
\smallskip

{\bf Definition 11} (\emph{rigid graph tangles}).
Let $\Ga$ be a non-oriented non-connected infinite graph
 with vertices of degree $\leq n$.
A subset $\Ga\subset\{0\leq z\leq 1\}$, homeomorphic to
 the graph $\Ga$, is called \emph{a graph tangle}, if (see Fig.~5):

(11.1) the set of the 1-vertices of $\Ga$
 coincides with the set of the marked points 
 
 $\{(k,0,0),(k,0,1) \vl k\in \N \}$;

(11.2) the connected components of $\Ga$ lying sufficiently far
 from the origin $0\in \R^3$ are line segments joining
 points $(k,0,0)$ and $(j,0,1)$,
 where the difference $k-j$ is constant for all sufficiently big $j$;

(11.3) a neighbourhood of each vertex $A\in \Ga$
 lies in a plane\footnote{
 Possibly, in a bowed disk as neighbourhoods of $m$-vertices in Fig.~2b.}.
\smallskip

If isotopy of graph tangles inside the layer $\{0<z<1\}$
 keeps Condition~(11.3), then the corresponding isotopy classes
 of graph tangles are called \emph{rigid graph tangles}.
\ed
\smallskip

One can represent graph tangles by 
 their plane diagrams similarly to
 spatial graphs.
\emph{The product} $\Ga_1 \times \Ga_2$ of two graph tangles 
 is the graph tangle obtained by attaching
 the top semiline of $\Ga_2$ to the bottom semiline of $\Ga_1$ and
 then by contracting the new layer $\{0\leq z\leq 2\}$ to 
 the initial one.
So, all rigid isotopy classes of graph tangles
 form a semigroup $RGT_n$.
\emph{The unit graph tangle} $1\in RGT_n$ consists of vertical segments
 joining the points $(k,0,0)$ and $(k,0,1)$, $k\in \N$.
Let us introduce the following rigid graph tangles:
 $\T_n = \{\; \xi_k, \; \eta_k, \; \sm_k, \; \si, \; \la_{m,k} \vl
  k\geq 1, \; 3\leq m\leq n \; \}.$


\begin{picture}(410,75)(-30,-20)

\put(-40,35){\vector(1,0){15}}
\put(-40,35){\vector(0,-1){50}}
\put(-40,40){0}
\put(-35,5){1}
\put(-40,35){\circle*{3}}
\put(-40,5){\circle*{3}}
\put(-25,40){$x$}
\put(-35,-10){$z$}

\thicklines

\put(0,35){\line(1,-1){15}}
\put(15,20){\line(1,1){15}}
\put(-10,5){\line(0,1){30}}
\put(0,5){\line(3,1){30}}
\put(30,15){\line(1,2){10}}
\put(30,5){\line(1,1){10}}
\put(-5,40){\footnotesize $k$}
\put(15,40){\footnotesize $k+1$}
\put(10,-10){$\xi_k$}

\put(0,35){\circle*{3}}
\put(30,35){\circle*{3}}
\put(0,5){\circle*{3}}
\put(30,5){\circle*{3}}
\put(-10,5){\circle*{3}}
\put(-10,35){\circle*{3}}
\put(40,5){\circle*{3}}
\put(40,35){\circle*{3}}

\put(65,5){\line(0,1){30}}
\put(75,5){\line(1,1){15}}
\put(90,20){\line(1,-1){15}}
\put(75,35){\line(3,-1){30}}
\put(105,25){\line(1,-2){10}}
\put(105,35){\line(1,-1){10}}
\put(70,40){\footnotesize $k$}
\put(90,40){\footnotesize $k+1$}
\put(85,-10){$\eta_k$}

\put(65,5){\circle*{3}}
\put(65,35){\circle*{3}}
\put(75,35){\circle*{3}}
\put(105,35){\circle*{3}}
\put(75,5){\circle*{3}}
\put(105,5){\circle*{3}}
\put(115,5){\circle*{3}}
\put(115,35){\circle*{3}}

\put(140,5){\line(0,1){30}}
\put(150,5){\line(1,1){30}}
\put(150,35){\line(1,-1){12}}
\put(168,17){\line(1,-1){12}}
\put(190,5){\line(0,1){30}}
\put(145,40){\footnotesize $k$}
\put(165,40){\footnotesize $k+1$}

\put(140,5){\circle*{3}}
\put(140,35){\circle*{3}}
\put(150,35){\circle*{3}}
\put(180,35){\circle*{3}}
\put(150,5){\circle*{3}}
\put(180,5){\circle*{3}}
\put(190,5){\circle*{3}}
\put(190,35){\circle*{3}}
\put(160,-10){$\sm_k$}

\put(215,5){\line(0,1){30}}
\put(225,5){\line(1,1){12}}
\put(243,23){\line(1,1){12}}
\put(225,35){\line(1,-1){30}}
\put(265,5){\line(0,1){30}}
\put(220,40){\footnotesize $k$}
\put(240,40){\footnotesize $k+1$}

\put(215,5){\circle*{3}}
\put(215,35){\circle*{3}}
\put(225,35){\circle*{3}}
\put(225,5){\circle*{3}}
\put(255,5){\circle*{3}}
\put(255,35){\circle*{3}}
\put(265,5){\circle*{3}}
\put(265,35){\circle*{3}}
\put(235,-10){$\si$}

\put(290,5){\line(0,1){30}}
\put(300,5){\line(1,1){15}}
\put(300,35){\line(1,-1){30}}
\put(340,5){\line(-1,3){10}}
\put(340,35){\line(1,-3){5}}
\put(295,40){\footnotesize $k$}
\put(315,40){\footnotesize $k+1$}

\put(290,5){\circle*{3}}
\put(290,35){\circle*{3}}
\put(300,5){\circle*{3}}
\put(300,35){\circle*{3}}
\put(315,20){\circle*{5}}
\put(330,5){\circle*{3}}
\put(330,35){\circle*{3}}
\put(340,5){\circle*{3}}
\put(340,35){\circle*{3}}
\put(310,-10){$\la_{3,k}$}

\put(365,5){\line(0,1){30}}
\put(375,5){\line(1,1){30}}
\put(375,35){\line(1,-1){30}}
\put(415,5){\line(0,1){30}}
\put(370,40){\footnotesize $k$}
\put(390,40){\footnotesize $k+1$}

\put(365,5){\circle*{3}}
\put(365,35){\circle*{3}}
\put(375,5){\circle*{3}}
\put(375,35){\circle*{3}}
\put(390,20){\circle*{5}}
\put(405,5){\circle*{3}}
\put(405,35){\circle*{3}}
\put(415,5){\circle*{3}}
\put(415,35){\circle*{3}}
\put(385,-10){$\la_{4,k}$}

\put(50,-30){{\bf Fig.~5.} 
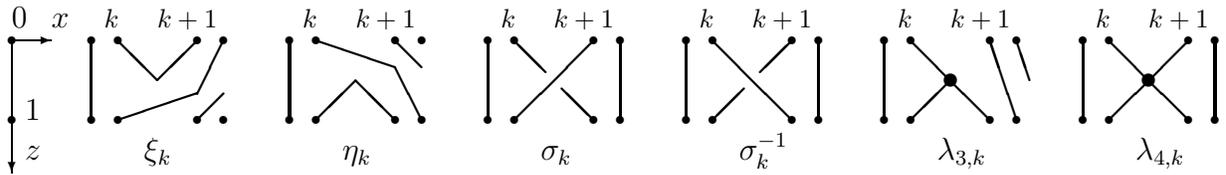
 The graph tangles of the alphabet $\T_4$.}

\end{picture}
\vspace{7mm}

For any $k\in\N$,
 the tangle $\eta_k\xi_k$ is the unknot added to the unit $1\in RGT_n$.


\subsection{ The semigroup $RGT_n$ of rigid graph tangles}

We are working in the PL-category, i.e. graph tangles
 $\Ga\subset\{0\leq z\leq 1\}$ consist of finite polygonal lines.
\smallskip

{\bf Definition~12} (\emph{the graph $\Ga_{xz}$, extremal points,
 singularities}).

(a) By $\Ga_{xz}$ denote the image of a graph tangle
 $\Ga\subset\{0\leq z\leq 1\}$
 under the projection to the $xz$-plane, see Fig.~5.
\emph{Extremal points} of $\Ga_{xz}$
 are images under the $xz$-projection of
 local maxima and minima of the $z$-coordinate on
 the interiors of the edges of $\Ga$.

(b) The images on the $xz$-plane of the vertices of $\Ga$
 (except the 1-vertices),
 the crossings of $\Ga_{xz}$, and the extremal points of $\Ga_{xz}$
 are called \emph{singularities} of $\Ga_{xz}$.
\ed
\smallskip

In Fig.~5, each tangle has exactly one singularity.
The tangles $\xi_k$ and $\eta_k$ have an extremal point.
The tangles $\sm_k$ and $\sm_k^{-1}$ have a crossing,
 $\la_{m,k}$ contains an $m$-vertex.
\smallskip
 
{\bf Definition~13} (\emph{graphs $\Ga_{xz}$ in general position}).
Let us say that a graph $\Ga_{xz}$
 is in \emph{general position} on the $xz$-plane, if
 the following conditions hold (see Fig.~5):
\smallskip

(13.1) the graph $\Ga_{xz}$ has finitely many singularities;

(13.2) each crossing is not an extremal point;

(13.3) for every $(2p-1)$-vertex $C\in\Ga_{xz}$,
 exactly $p-1$ arcs with the endpoint $C$ go up\footnote{
 It means that these $p-1$ arcs point toward to
  the negative direction of the $z$-axis, see Fig.~5.},
 and the other $p$ arcs with the enpoint $C$ go down;

(13.4) for any $2q$-vertex $B\in\Ga_{xz}$,
 exactly $q$ arcs with the endpoint $B$ go up,
 and the other $q$ arcs with the enpoint $B$ go down;

(13.5) no two singularities lie in \emph{a horizontal line},
 which is parallel to the $x$-axis.
\ed
\smallskip

The following lemma extends results of the Turaev work \cite{Tu}
 to rigid graph tangles.

\begin{lemma}
The semigroup $RGT_n$ 
 is generated by the letters of $\T_n$
 and the relations $(11)-(23)$, where $k\geq  1$,
 $3\leq m\leq n$, $2\leq p\leq \frac{n+1}{2}$, 
 $2\leq q\leq \frac{n}{2}$.

$$ \begin{array}{llllll}
(11) &  \xi_k \xi_l = \xi_{l+2} \xi_k, &
        \xi_k \eta_l = \eta_{l+2} \xi_k, &
        \xi_k \sm_l = \sm_{l+2} \xi_k, &
        \xi_k \la_{m,l} = \la_{m,l+2} \xi_k &
        (l\geq k); \\

(12) &  \eta_k \xi_l = \xi_{l-2} \eta_k, &
         \eta_k \eta_l = \eta_{l-2} \eta_k, &
         \eta_k \sm_l = \sm_{l-2} \eta_k, &
         \eta_k \la_{m,l} = \la_{m,l-2} \eta_k &
         (l\geq k+2); \\

(13) & \sm_k \xi_l = \xi_{l} \sm_k, &
        \sm_k \eta_l = \eta_{l} \sm_k, &
        \sm_k \sm_l = \sm_{l} \sm_k, &
        \sm_k \la_{m,l} = \la_{m,l} \sm_k &
        (l\geq k+2);
\end{array}$$
$$(14) \left\{ \begin{array}{lll}
        \la_{2p-1,k} \xi_l =  \xi_{l-1} \la_{2p-1,k}, &
        \la_{2p-1,k} \sm_l =  \sm_{l-1} \la_{2p-1,k}, & \\

        \la_{2p-1,k} \eta_l =  \eta_{l-1} \la_{2p-1,k}, &
        \la_{2p-1,k} \la_{m,l} =  \la_{m,l-1} \la_{2p-1,k} &
        (l\geq k+p), \\

        \la_{2q,k} \xi_l =  \xi_{l} \la_{2q,k}, &
        \la_{2q,k} \sm_l =  \sm_{l} \la_{2q,k}, & \\
        
        \la_{2q,k} \eta_l =  \eta_{l} \la_{2q,k}, &
        \la_{2q,k} \la_{m,l} =  \la_{m,l} \la_{2q,k} & 
        (l\geq k+q);
\end{array} \right.$$
$$\begin{array}{ll}
(15) \;  \eta_{k+1} \xi_k = 1 = \eta_{k} \xi_{k+1}; & \\

(16) \;  \eta_{k+2} \sm_{k+1} \xi_k = \si = \eta_k\sm_{k+1}\xi_{k+2}; &
(19) \; \eta_k \sm_k = \eta_k, \; \sm_k \xi_k = \xi_k; \\

(17) \; \eta_{k+p-1} \la_{2p-1,k+1} \xi_k = \la_{2p-1,k} =
        \eta_k \la_{2p-1,k+1} \xi_{k+p}; &
(20) \; \sm_k \si = 1 = \si \sm_k; \\

(18) \; \eta_{k+q} \la_{2q,k+1} \xi_k = \la_{2q,k} =
        \eta_k \la_{2q,k+1} \xi_{k+q}; &
(21) \; \sm_k \sm_{k+1} \sm_k = \sm_{k+1} \sm_k \sm_{k+1};  \\
\end{array}$$
$$(22) \; \left\{ \begin{array}{ll}
  \la_{2p-1,k+1} \Sm_{k,p} = \Sm_{k,p-1} \la_{2p-1,k}, &
  \la_{2p-1,k} \bar\Sm_{k,p} = \bar\Sm_{k,p-1} \la_{2p-1,k+1}, \\

  \la_{2q,k+1} \Sm_{k,q} = \Sm_{k,q} \la_{2q,k}, &
  \la_{2q,k} \bar\Sm_{k,q} = \bar\Sm_{k,q} \la_{2q,k+1},  \\

  \mbox{ where } \Sm_{k,l} = \sm_k \sm_{k+1} \cdots \sm_{k+l-1}, &
  \bar\Sm_{k,l} = \sm_{k+l-1} \cdots \sm_{k+1} \sm_k \; (l\geq 1);
\end{array} \right.$$
$$(23) \; \left\{ \begin{array}{ll}
 \la_{2p-1,k} \Sm'_{k,p-1} = \Sm'_{k,p-2} \la_{2p-1,k}, &
 \la_{2q,k} \Sm'_{k,q-1} = \Sm'_{k,q-1} \la_{2q,k},  \\

 \mbox{where } \Sm'_{k,0}=1 \mbox{ and, for } l\geq 1, &
 \Sm'_{k,l} = \sm_{k+l-1}^{-1} (\sm_{k+l-2}^{-1} \sm_{k+l-1}^{-1}) \cdots
  (\sm_{k}^{-1} \sm_{k+1}^{-1} \cdots \sm_{k+l-1}^{-1}).
\end{array} \right.$$
\end{lemma}
\begin{proof}
A graph $\Ga_{xz}$ will be in general position on the $xz$-plane after
 a small deformation of the graph $\Ga_{xz}$.
Then the $xz$-plane splits by horizontal lines into
 strips each of that contains exactly one singularity of $\Ga_{xz}$.
Look at all singularities of $\Ga_{xz}$ from the top to the bottom.
Write the corresponding generators from the left to the right.
Then we get a word $u_{\Ga}$ in $\T_n$, see Fig.~5.
The generators $\xi_k,\eta_k$ represent extremal points;
 $\sm_k,\si$ denote overcrossings and undercrossings;
 the letter $\la_{m,k}$ corresponds to some $m$-vertex.
It suffices to prove that any rigid isotopy of graph tangles
 decomposes on the relations $(11)-(23)$.
By generalized Reidemeister's theorem \cite{JM} of Subsection~2.1
 and general position reasonings,
 any rigid isotopy of graph tangles decomposes on
 the following moves:
\smallskip

 1) general position isotopies;

 2) swappings heights of two singularities;

 3) creations or annihilations of two close extremal points;

 4) isotopies of either a crossing or a vertex near an extremal point;

 5) the Reidemeister moves $R1-R5$ of Fig.~1.
\smallskip

The first isotopies preserve the word $u_{\Ga}$ in $\T_n$.
The second isotopies are described by the relations $(11)-(14)$.
The third isotopies provide the relations (15).
It was shown in \cite[proof of Lemma 3.4]{Tu} that 
 all isotopies of a crossing near an extremal point
 decompose geometrically on the relations (16).
Analogously, one can check that all isotopies of a vertex near
 an extremal point decompose on the relations $(17)-(18)$.
Finally, the relations $(19)-(23)$ correspond to the
 Reidemeister moves $R1-R5$, respectively.
\end{proof}


\subsection{Three-page tangles.}

Let us consider three semilines with a common endpoint in
 the horizontal plane $\{z=0\}$.
For example, put:
 $$Y=\{x\geq 0,y=z=0\}\cup \{y\geq 0,x=z=0\}\cup
  \{x\leq 0,y=z=0\}\subset \{z=0\}.$$
Mark out all the integer points in the semilines:
 $\{(j,0,0),(0,k,0),(-l,0,0) \vl j,k,l\in \N \}$.
Let $I\subset \R^3$ be the segment joining
 the points $(0,0,0)$ and $(0,0,1)$.
Put (see Fig.~6):
 $$P_0=\{x>0,y=z=0\}\times I,\quad
   P_1=\{y>0,x=z=0\}\times I,\quad
   P_2=\{x>0,y=z=0\}\times I.$$
Then the product $Y\times I$ is 
 the three-page book with the pages $P_i$, see Definition~7.
\smallskip

{\bf Definition 14} (\emph{rigid three-page tangles}).
Let $\Ga$ be a non-oriented non-connected infinite graph
 with vertices of degree $\leq n$.
A subset $\Ga\subset Y\times I$, homeomorphic to
 the graph $\Ga$, is called \emph{a three-page tangle}, if
 the following conditions hold (see Fig.~6):
\smallskip

(14.1) the set of the 1-vertices of $\Gamma$
  coincides with the set of the marked points
  $$\{(j,0,0),(j,0,1), (0,k,0),(0,k,1), (-l,0,0),(-l,0,1) \vl
   j,k,l\in \N \};$$

(14.2) the vertices of degree $\geq 3$ lie in the segment $I$;

(14.3) \emph{finiteness}: the intersection
  $\Gamma \cap I = A_1\cup \dots \cup A_m$ is a finite set of points;

(14.4) two arcs with an endpoint $A_j\in\Ga\cap I$ that is not a vertex of $\Ga$
  lie in different pages $P_i,P_j$ ($i\neq j$);

(14.5) \emph{balance:} neighbourhoods of
  vertices $A_j\in \Ga\cap \al$
  look like pictures of Fig.~2b;

(14.6) \emph{monotone}: for any $i\in \Z_3$, the restriction of the orthogonal
  projection $Y\times I\to I$ to each connected component of
  the intersection $\Gamma \cap P_i$ is a monotone function.

(14.7) for every $i\in \Z_3$, all connected components of $\Gamma$
 lying in the page $P_i$ 
 sufficiently far from the origin $0\in\R^3$ are 
 line segments that are parallel to each other.
\smallskip

If three-page tangles are considered up to
 rigid (respectively, non-rigid) isotopy in the layer $\{0<z<1\}$,
 then they are called \emph{rigid} (respectively, \emph{non-rigid}).
\ed
\smallskip

The reader can compare the above conditions with
 Definitions~7 and 11, see Fig.~6.
\vspace{5mm}


\begin{picture}(360,75)(-30,0)

\put(-10,70){\vector(1,0){75}}
\put(60,75){$x$}
\put(20,70){\vector(-2,-1){45}}
\put(-25,40){$y$}
\put(55,55){$P_0$}
\put(-25,20){$P_1$}
\put(-20,60){$P_2$}
\put(10,15){$1$}
\put(15,30){\line(1,0){40}}
\put(-5,30){\line(1,0){10}}
\put(20,30){\line(-2,-1){35}}
\put(20,70){\vector(0,-1){55}}
\put(25,15){$z$}
\put(20,75){0}
\put(0,-5){$\ph (\xi_1)=d_2c_2$}

{\thicklines
\put(30,70){\line(-1,-1){20}}
\put(20,40){\line(-1,1){10}}
\put(20,40){\line(2,3){20}}
\put(30,30){\line(1,2){20}}
\put(40,30){\line(1,2){10}}
\put(10,65){\line(-1,-3){5}}
\put(5,50){\line(1,-5){5}}
\put(0,30){\line(0,1){5}}
\put(0,40){\line(0,1){5}}
\put(0,50){\line(0,1){5}}
\put(0,65){\line(0,1){5}}
\put(-10,15){\line(0,1){40}}
}

\put(0,70){\circle*{3}}
\put(0,30){\circle*{3}}
\put(10,66){\circle*{3}}
\put(10,25){\circle*{3}}
\put(-10,55){\circle*{3}}
\put(-10,15){\circle*{3}}
\put(20,60){\circle*{3}}
\put(20,40){\circle*{3}}
\put(30,70){\circle*{3}}
\put(30,30){\circle*{3}}
\put(40,70){\circle*{3}}
\put(40,30){\circle*{3}}
\put(50,70){\circle*{3}}
\put(50,30){\circle*{3}}

\put(90,70){\vector(1,0){75}}
\put(155,55){$P_0$}
\put(75,20){$P_1$}
\put(80,60){$P_2$}
\put(110,15){$1$}
\put(160,75){$x$}
\put(120,70){\vector(-2,-1){45}}
\put(75,40){$y$}
\put(115,30){\line(1,0){40}}
\put(95,30){\line(1,0){10}}
\put(120,30){\line(-2,-1){35}}
\put(120,70){\line(0,-1){40}}
\put(120,70){\vector(0,-1){55}}
\put(125,15){$z$}
\put(120,75){0}
\put(100,-5){$\ph (\eta_1)=a_2 b_2$}

{\thicklines
\put(110,50){\line(1,-1){20}}
\put(110,50){\line(1,1){10}}
\put(120,60){\line(2,-3){20}}
\put(130,70){\line(1,-2){20}}
\put(140,70){\line(1,-2){10}}
\put(110,65){\line(-1,-3){5}}
\put(105,50){\line(1,-5){5}}
\put(100,30){\line(0,1){5}}
\put(100,40){\line(0,1){5}}
\put(100,50){\line(0,1){5}}
\put(100,65){\line(0,1){5}}
\put(90,15){\line(0,1){40}}
}
\put(100,70){\circle*{3}}
\put(100,30){\circle*{3}}
\put(110,66){\circle*{3}}
\put(110,25){\circle*{3}}
\put(90,55){\circle*{3}}
\put(90,15){\circle*{3}}
\put(120,60){\circle*{3}}
\put(120,40){\circle*{3}}
\put(130,70){\circle*{3}}
\put(130,30){\circle*{3}}
\put(140,70){\circle*{3}}
\put(140,30){\circle*{3}}
\put(150,70){\circle*{3}}
\put(150,30){\circle*{3}}

\put(195,70){\vector(1,0){80}}
\put(265,55){$P_0$}
\put(185,10){$P_1$}
\put(190,60){$P_2$}
\put(220,0){$1$}
\put(270,75){$x$}
\put(220,15){\line(1,0){55}}
\put(205,15){\line(1,0){10}}
\put(230,15){\line(-2,-1){35}}
\put(230,70){\vector(-2,-1){45}}
\put(185,40){$y$}
\put(230,70){\vector(0,-1){70}}
\put(235,0){$z$}
\put(230,75){0}
\put(200,-15){$\ph (\sm_1)=b_1 d_2 d_1 b_2$}

\put(230,20){\circle*{3}}
\put(230,35){\circle*{3}}
\put(230,50){\circle*{3}}
\put(230,65){\circle*{3}}
\put(240,70){\circle*{3}}
\put(250,70){\circle*{3}}
\put(260,70){\circle*{3}}
\put(240,15){\circle*{3}}
\put(250,15){\circle*{3}}
\put(260,15){\circle*{3}}
\put(200,0){\circle*{3}}
\put(200,55){\circle*{3}}
\put(205,15){\circle*{3}}
\put(205,70){\circle*{3}}

{\thicklines
\put(250,70){\line(-1,-1){35}}
\put(250,15){\line(-1,1){25}}
\put(210,55){\line(1,-1){10}}
\put(210,55){\line(2,1){10}}
\put(240,70){\line(-2,-1){15}}
\put(230,65){\line(-1,-1){5}}
\put(215,35){\line(1,-1){15}}
\put(215,35){\line(1,1){35}}
\put(240,15){\line(-2,1){10}}
\put(250,15){\line(-1,1){20}}
\put(260,15){\line(0,1){55}}
\put(200,0){\line(0,1){55}}
\put(205,15){\line(0,1){5}}
\put(205,25){\line(0,1){5}}
\put(205,35){\line(0,1){5}}
\put(205,45){\line(0,1){5}}
\put(205,65){\line(0,1){5}}
}

\put(320,70){\vector(1,0){80}}
\put(385,25){$P_0$}
\put(305,20){$P_1$}
\put(310,60){$P_2$}
\put(340,5){$1$}
\put(395,75){$x$}
\put(340,20){\line(1,0){60}}
\put(324,20){\line(1,0){8}}
\put(350,70){\vector(-2,-1){45}}
\put(305,40){$y$}
\put(350,20){\line(-2,-1){35}}
\put(350,70){\vector(0,-1){65}}
\put(355,5){$z$}
\put(350,75){0}
\put(330,-10){$\ph (\la_{3,1})= x_{3,2} b_2$}

\put(350,25){\circle*{3}}
\put(350,45){\circle*{5}}
\put(360,70){\circle*{3}}
\put(375,70){\circle*{3}}
\put(385,70){\circle*{3}}
\put(360,20){\circle*{3}}
\put(375,20){\circle*{3}}
\put(385,20){\circle*{3}}
\put(330,70){\circle*{3}}
\put(330,20){\circle*{3}}
\put(335,64){\circle*{3}}
\put(335,14){\circle*{3}}
\put(320,55){\circle*{3}}
\put(320,5){\circle*{3}}

{\thicklines
\put(350,45){\line(2,5){10}}
\put(350,45){\line(1,-1){25}}
\put(340,35){\line(1,1){10}}
\put(340,35){\line(1,-1){10}}
\put(360,20){\line(-2,1){10}}
\put(385,20){\line(-1,5){10}}
\put(385,70){\line(1,-5){5}}
\put(335,13){\line(0,1){50}}
\put(320,6){\line(0,1){50}}
\put(330,20){\line(0,1){5}}
\put(330,30){\line(0,1){5}}
\put(330,40){\line(0,1){5}}
\put(330,50){\line(0,1){5}}
\put(330,65){\line(0,1){5}}
}
\end{picture}
\vspace{0.5cm}


\begin{picture}(360,100)(0,0)

\put(20,80){\vector(1,0){80}}
\put(95,85){$x$}
\put(40,30){\line(1,0){60}}
\put(24,30){\line(1,0){8}}
\put(50,80){\vector(-2,-1){45}}
\put(5,50){$y$}
\put(50,30){\line(-2,-1){35}}
\put(50,80){\line(0,-1){50}}
\put(50,80){\vector(0,-1){65}}
\put(55,15){$z$}
\put(50,85){0}
\put(30,0){$\phi (\la_{4,1})= d_2 x_{4,2} b_2$}

\put(50,35){\circle*{3}}
\put(50,55){\circle*{5}}
\put(50,75){\circle*{3}}
\put(60,80){\circle*{3}}
\put(75,80){\circle*{3}}
\put(85,80){\circle*{3}}
\put(60,30){\circle*{3}}
\put(75,30){\circle*{3}}
\put(85,30){\circle*{3}}
\put(30,80){\circle*{3}}
\put(30,30){\circle*{3}}
\put(35,74){\circle*{3}}
\put(35,24){\circle*{3}}
\put(20,65){\circle*{3}}
\put(20,15){\circle*{3}}

{\thicklines
\put(50,55){\line(-1,1){10}}
\put(50,55){\line(-1,-1){10}}
\put(50,55){\line(1,-1){25}}
\put(50,55){\line(1,1){25}}
\put(40,65){\line(1,1){10}}
\put(40,45){\line(1,-1){10}}
\put(60,80){\line(-2,-1){10}}
\put(60,30){\line(-2,1){10}}
\put(85,30){\line(0,1){50}}
\put(35,23){\line(0,1){50}}
\put(20,16){\line(0,1){50}}
\put(30,30){\line(0,1){5}}
\put(30,40){\line(0,1){5}}
\put(30,50){\line(0,1){5}}
\put(30,60){\line(0,1){5}}
\put(30,75){\line(0,1){5}}
}

\put(140,80){\vector(1,0){105}}
\put(240,85){$x$}
\put(160,30){\line(1,0){80}}
\put(144,30){\line(1,0){8}}
\put(170,80){\vector(-2,-1){45}}
\put(125,50){$y$}
\put(170,30){\line(-2,-1){35}}
\put(170,80){\vector(0,-1){65}}
\put(175,15){$z$}
\put(170,85){0}
\put(150,0){$\phi (\la_{5,1})= x_{5,2} b_2$}

\put(170,35){\circle*{3}}
\put(170,55){\circle*{5}}
\put(180,80){\circle*{3}}
\put(195,80){\circle*{3}}
\put(220,80){\circle*{3}}
\put(230,80){\circle*{3}}
\put(180,30){\circle*{3}}
\put(220,30){\circle*{3}}
\put(195,30){\circle*{3}}
\put(230,30){\circle*{3}}
\put(150,80){\circle*{3}}
\put(150,30){\circle*{3}}
\put(155,74){\circle*{3}}
\put(155,24){\circle*{3}}
\put(140,65){\circle*{3}}
\put(140,15){\circle*{3}}

{\thicklines
\put(170,55){\line(-1,-1){10}}
\put(170,55){\line(2,5){10}}
\put(170,55){\line(1,1){25}}
\put(170,55){\line(1,-1){25}}
\put(170,55){\line(2,-1){50}}
\put(160,45){\line(1,-1){10}}
\put(180,30){\line(-2,1){10}}
\put(230,30){\line(-1,5){10}}
\put(230,80){\line(1,-5){5}}
\put(155,23){\line(0,1){50}}
\put(140,16){\line(0,1){50}}
\put(150,30){\line(0,1){5}}
\put(150,40){\line(0,1){5}}
\put(150,50){\line(0,1){5}}
\put(150,60){\line(0,1){5}}
\put(150,75){\line(0,1){5}}
}

\put(300,80){\vector(1,0){105}}
\put(400,85){$x$}
\put(320,30){\line(1,0){80}}
\put(304,30){\line(1,0){8}}
\put(330,80){\vector(-2,-1){45}}
\put(285,50){$y$}
\put(330,30){\line(-2,-1){35}}
\put(330,80){\vector(0,-1){65}}
\put(335,15){$z$}
\put(330,85){0}
\put(310,0){$\phi (\la_{6,1})= d_2 x_{6,2} b_2$}

\put(330,35){\circle*{3}}
\put(330,55){\circle*{5}}
\put(330,75){\circle*{3}}
\put(340,80){\circle*{3}}
\put(355,80){\circle*{3}}
\put(380,80){\circle*{3}}
\put(390,80){\circle*{3}}
\put(340,30){\circle*{3}}
\put(355,30){\circle*{3}}
\put(380,30){\circle*{3}}
\put(390,30){\circle*{3}}
\put(310,80){\circle*{3}}
\put(310,30){\circle*{3}}
\put(315,74){\circle*{3}}
\put(315,24){\circle*{3}}
\put(300,65){\circle*{3}}
\put(300,15){\circle*{3}}

{\thicklines
\put(330,55){\line(-1,1){10}}
\put(330,55){\line(-1,-1){10}}
\put(330,55){\line(1,1){25}}
\put(330,55){\line(2,1){50}}
\put(330,55){\line(1,-1){25}}
\put(330,55){\line(2,-1){50}}
\put(320,65){\line(1,1){10}}
\put(320,45){\line(1,-1){10}}
\put(340,80){\line(-2,-1){10}}
\put(340,30){\line(-2,1){10}}
\put(390,30){\line(0,1){50}}
\put(315,23){\line(0,1){50}}
\put(300,16){\line(0,1){50}}
\put(310,30){\line(0,1){5}}
\put(310,40){\line(0,1){5}}
\put(310,50){\line(0,1){5}}
\put(310,60){\line(0,1){5}}
\put(310,75){\line(0,1){5}}
}

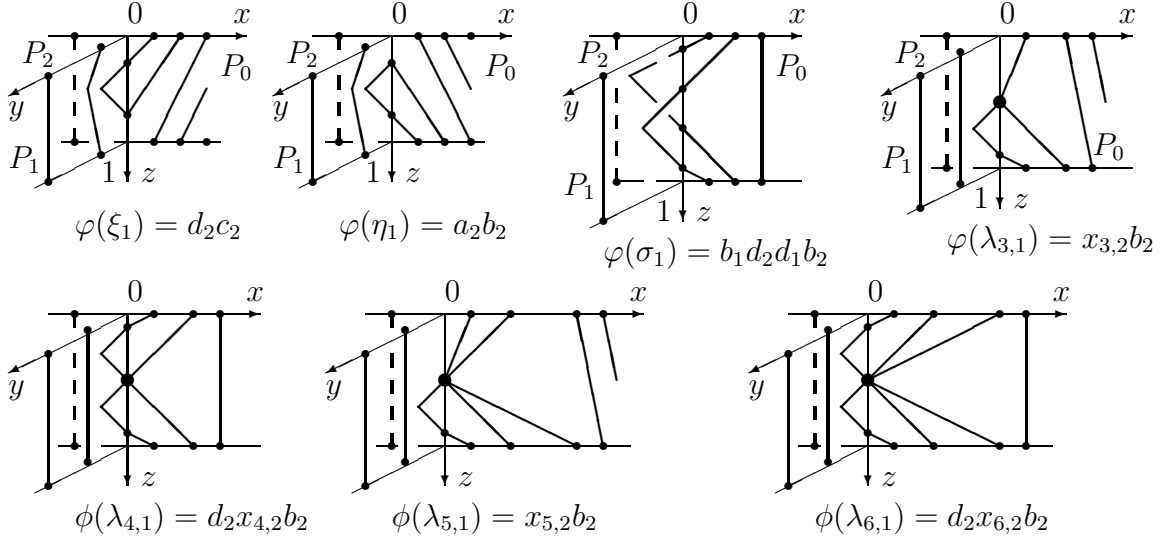
\put(20,-20){ {\bf Fig.~6.} The three-page tangles associated
 with the graph tangles of $\T_6$}
\end{picture}
\vspace{1cm}

All rigid isotopy classes of three-page tangles form a semigroup.
Proposition~1 of Subsection~5.2 shows that
 this semigroup is isomorphic to $RSG_n$.
Any three-page tangle $\Ga\subset\{0\leq z\leq 1\}$
 could be encoded by a word $w_{\Ga}$
 in the alphabet $\A_n$ (see Fig.~2a and 2b) 
 analogously to the proof of Theorem~1a.


\setcounter{section}{4}
\section{ Proofs of Theorems 1b--1c, 2, and Corollaries 1--4}

Theorems~1b and 1c will be proved in Subsections~5.2 and 5.3,
 respectively.
Corollaries~1a and 2a will be checked at the end of
 Subsection~5.3.
Subsection~5.4 is devoted to the proofs of Corollaries~1b, 2b,
 and 3--4.

\subsection{ The semigroup $RBT_n$ of almost balanced tangles}
We are going to select almost balanced tangles among three-page tangles.
The semigroup $RBT_n$ of rigid almost balanced tangles  
 will be isomorphic to the semigroup $RGT_n$ of rigid graph tangles.
\smallskip

{\bf Definition 15} (\emph{almost balanced tangles}).
A three-page tangle $\Ga\subset Y\times I$ is called \emph{almost balanced},
 if the corresponding word $w_{\Ga}$ in $\A_n$ is
 simultaneously 1-balanced and 2-balanced.
Equivalently, one can assume that
 the marked points lying in $P_1,P_2$ are joined
 in any \emph{almost balanced tangle} by vertical segments
 parallel to the $z$-axis.
By $RBT_n$ denote the semigroup forming by
 all rigid isotopy classes of
 almost balanced tangles\footnote{
 \emph{Rigidity} means that under rigid isotopy in
  the layer $\{0<z<1\}$ a neighbourhood of each vertex
  in any almost balanced tangle lies in a plane (possibly, in a
  bowed disk as in Fig.~2b).}.
\ed
\smallskip

Any graph tangle of Definition~11
 could be embedded into $Y\times I$ so that
 its 1-vertices lie in the
 semilines $\{x\geq 0,y=z=0\}$ and $\{x\geq 0,y=0,z=1\}$.
Then we may add two infinite families of vertical segments in $P_1,P_2$
 and get an almost balanced tangle.
Since graph tangles and three-page tangles 
 are defined up to rigid isotopy
 in $\{0<z<1\}$, then a non-canonical monomorphism
 $\ph: RGT_n\to RBT_n$ is well-defined.
\smallskip

{\bf Definition 16} (\emph{the canonical isomorphism $\ph: RGT_n\to RBT_n$}).
Let us take the map $\ph: RGT_n\to RBT_n$ defined on generators as follows
 ($k\in\N$, see Fig.~6):
 $$(24)\quad \left\{ \begin{array}{lll}
   \ph(\xi_k)=d_2^k c_2 b_2^{k-1}, &
   \ph(\sm_k)=d_2^{k-1} b_1 d_2 d_1 b_2^k, &
   \ph(\la_{2p-1,k})=d_2^{k-1} x_{2p-1,2} b_2^k, \\

   \ph(\eta_k)=d_2^{k-1} a_2 b_2^k, &
   \ph(\si)=d_2^k b_1 b_2 d_1 b_2^{k-1}, &
   \ph(\la_{2q,k})=d_2^k x_{2q,2} b_2^k. 
 \end{array} \right. \eqno{\blacksquare} $$

\begin{lemma}
The map $\ph: RGT_n\to RBT_n$ is a well-defined isomorphism of semigroups.
\end{lemma}
\begin{proof}
Let us construct an inverse map $\psi: RBT_n\to RGT_n$.
Associate with each almost balanced tangle $\Ga\in RBT_n$
 a graph tangle $\psi(\Ga)\in RGT_n$ given by the following plane diagram.
By Definition~15 one can suppose that
 the marked points, lying in the pages $P_1,P_2$,
 are joined in $\Ga$ by vertical segments.
Delete from $\Ga$ all these vertical segments.
We get a graph tangle $\psi(\Ga)$ from Definition~11.
The composition $\psi\circ\ph: RGT_n\to RGT_n$ is
 identical on the generators of $RGT_n$.
So, the maps $\ph,\psi$ are the inverse ones.
\end{proof}

By $\ph(11)-\ph(23)$ we denote the images of
 the relations (11)--(23) under the isomorphism $\ph:RGT_n\to RBT_n$,
 i.e. $\ph(11)-\ph(23)$ are relations between words in $\A_n$.


\subsection{ Proof of Theorem~1b}

If two words $u,v\in W_n$ represent
 the same element of the semigroup $RSG_n$, 
 then call them \emph{equivalent} and denote by $u\s v$.
Theorem~1b is a particular case of Proposition~1.
Really, any spatial graph could be represented by a three-page
 tangle encoded by a balanced word.
If two balanced three-page tangles are rigidly isotopic, then
 the corresponding words are equal by Proposition~1 as required.
\qed

\begin{proposition}
The semigroup of all rigid three-page tangles
 (i.e. three-page tangles considered up to rigid isotopy in $\{0<z<1\}$)
 is isomorphic to the semigroup $RSG_n$.
\end{proposition}
\begin{proof}
As it was already mentioned at the end of Subsection~4.3, 
 with any three-page tangle
 one can associate a word $w\in W_n$ 
 and hence an element of the semigroup $RSG_n$.
Conversely, each element of $RSG_n$ encodes
 a three-page embedding.
Indeed, let us draw local three-page embeddings (see Fig.~2a and 2b)
 representing the letters of a given word $w$.
Then extend all arcs until some of them meet with each
 other and another arcs come to the boundary of $Y\times I$.
For any page $P_i$, let us consider the arcs
 coming to the semilines $P_i\cap\{z=0\}$ and $P_i\cap\{z=1\}$.
We may assume that these arcs end at the marked points
 $1,2,\ldots,n_i$ (say) and $1,2,\ldots,m_i$, respectively.
In each page $P_i$, consider the marked points
 $n_i+k$ and $m_i+k$ for all $k\in\N$.
We join the points $n_i+k$ and $m_i+k$ by adding infinitely many
 parallel segments lying in the page $P_i$.
We get a three-page tangle $\Ga(w)$ in the sense of
 Definition~14, see Subsection~4.3.
For example, Fig.~6 shows the three-page tangles
 corresponding to the following elements of $RSG_n$:
 $d_2 c_2$, $a_2 b_2$, $b_1 d_2 d_1 b_2$, $x_{3,2} b_2$,
 $d_2 x_{4,2} b_2$, $x_{5,2} b_2$, $d_2 x_{6,2} b_2$.
The relations $(1)-(10)$ of the semigroup $RSG_n$ could be
 performed by rigid isotopy inside the layer $\{0<z<1\}$, see Fig.~3.

It remains to prove that 
 any rigid isotopy of three-page tangles
 could be decomposed on the relations $(1)-(10)$.
It suffices to do this only for 
 almost balanced tangles.
Really, take the three-page tangle $\Ga$ given by a word $w_{\Ga}$.
Let us consider the marked points, lying in $P_1\cap\{z=0\}$ and
 $P_1\cap\{z=1\}$, that are joined in $\Ga$ with points in $I$.
Let\footnote{
 The numbers $n_i,m_i$ defined here coincide with the same numbers
  from the previous paragraph.}
 $n_1$ and $m_1$ be the maximal marks of the above points
 lying in $P_1\cap \{z=0\}$ and $P_1\cap \{z=1\}$, respectively.
Let $n_2$ and $m_2$ be the maximal marks of the similar marked points
 lying in $P_2\cap \{z=0\}$ and $P_2\cap \{z=1\}$, respectively\footnote{
 For the letters $a_0,b_0$, we get
  $n_1=n_2=0,m_1=m_2=1$ and $n_1=m_2=0,m_1=n_2=1$,
  respectively.}.
For any almost balanced tangle, we have
 $n_1=m_1=n_2=m_2=0$.
Then the word $b_1^{n_2} d_2^{n_1} w b_2^{m_1} d_1^{m_2}$
 is 1-balanced and 2-balanced.
By the relations (2) this map and its inverse
 send equivalent words to equivalent ones.

By Lemma~2 with every almost balanced tangle of $RBT_n$
 one can associate a graph tangle from the semigroup $RGT_n$.
For graph tangles, any rigid isotopy is already decomposed
 on the relations $(11)-(23)$ in Lemma~1.
Then Proposition~1 follows from Lemma~3,
 which will be checked in Subsection~6.3.
\end{proof}

\begin{lemma}
The relations $(1)-(10)$ of $RSG_n$ imply
 the relations $\ph(11)-\ph(23)$ of $RBT_n$.
\end{lemma}


\subsection{Proof of Theorem~1c and Corollaries~1a, 2a.}
Due to Proposition~1 one can identify any element $w\in RSG_n$
 with the corresponding three-page tangle $\Ga(w)$.
\smallskip

{\bf Definition~17} (\emph{knot-like three-page tangles}).
A three-page tangle is called \emph{knot-like}, 
 if it contains a spatial graph added to \emph{the unit tangle}\footnote{
 By definition \emph{the unit tangle} $\Ga(1)$ consists of
  the segments parallel
  to the $z$-axis, such that these segments join marked points with
  the same numbers in every page $P_i$.}
 $\Ga(1)$.
By the definition a knot-like tangle
 is encoded by a balanced word of $BW_n$.
\ed
\smallskip

Then Theorem~1c comes from Lemma~4 stating that
 all balanced words encode all central elements of $RSG_n$.
The algorithm to decide, whether an element of $RSG_n$ is balanced
 (or, equivalently, central), follows from Definition~10 in
 Subsection~3.4.

\begin{lemma}
An element $w\in RSG_n$ encodes a knot-like three-page tangle
 $\Ga(w)$ if and only if
 the element $w$ is central in the semigroup $RSG_n$.
\end{lemma}
\begin{proof}
The part "only if" follows from geometry reasonings:
 a spatial graph could be moved by rigid isotopy to
 any place of a given tangle, i.e. 
 a balanced element commutes with each other by Proposition~1.
Let $w$ be a central element in $RSG_n$.
Then, for each $l\in\N$, we have
 $b_i^l d_i^l w = w b_i^l d_i^l$.
Denote by $m$ (respectively, by $k$)
 the number of arcs in the three-page tangle $\Ga(w)$,
 going out to the left (respectively, to the right) 
 in the page $P_{i-1}$.
For sufficiently large $l$,
 the number of the arcs for the three-page tangle $\Ga(b_i^l d_i^l w)$
 going out to the left in $P_{i-1}$ is $l$,
 and that for $\Ga(w b_i^l d_i^l)$ is $m+l-k$, i.e. $m=k$.
Hence, for $l>m$, the word
 $a_0^l a_1^l w c_1^l c_0^l$ is $i$-balanced.
Since $w$ is central, the word  $w a_0^l a_1^l c_1^l c_0^l$ is also
 $i$-balanced for each $i\in \Z_3$.
Then the tangle $\Ga(w a_0^l a_1^l c_1^l c_0^l)$
 is knot-like as well as $\Ga(w)$.
\end{proof}
\smallskip

\emph{Proof of Corollary~1a.}
If a spatial graph $G$ could be encoded by a word $w_G$,
 then its mirror image $\bar G$ is encoded by $\rho_n(w_G)$. 
So, Theorem~1b implies Corollary~1a.
\qed
\medskip

\emph{Proof of Corollary~2a.}
Every element of $DG\subset DS=RSG_2\subset RSG_n$ is invertible.
Conversely, if an element $w\in RSG_n$ 
 is invertible, then by Proposition~1
 the corresponding three-page tangle $\Ga(w)$ 
 can not contain vertices of degree $\geq 3$.
Indeed, all the relations (1)--(10) preserve the number of the $m$-vertices
 of $\Ga(w)$ or, equivalently, the number of the letters
 $x_{m,0},x_{m,1},x_{m,2}$ in the word $w$.
Hence the element $w$ could contain the letters $a_i,b_i,c_i,d_i$
 only, i.e. $w\in RSG_2=DS$.
But the group of the invertible elements of $DS$
 is the Dynnikov group $DG$ as we know from \cite{Dy1}, i.e. $w\in DG$.
\qed


\subsection{Non-rigid spatial graphs and spatial $J$-graphs}

{\it Proof of Theorem~2.}
Theorem~2 is proved analogously to Theorem~1, if we replace
 the relations~(9)--(10) by
 $$(9')\quad x_{m,i} b_i (d_i^2 d_{i+1}^2 d_{i-1}^2) = x_{m,i} b_i,
 \mbox{ where } 3\leq m\leq n,\; i\in\Z_3.$$
Let us emphasize only key moments.
Non-rigid graph tangles considered up to non-rigid
 isotopy in $\{0<z<1\}$ form a semigroup $NGT_n$.
As in Lemma~1 the semigroup $NGT_n$ is generated by the letters
 of the alphabet $\T_n$ and the relations~(11)-(22),
 $(23')\; \la_{m,k} \sm_k = \la_{m,k}$, where $3\leq m\leq n$, $k\in\N$.
The new relations~$(23')$ correspond to the Reidemeister move $R5'$,
 which switches two arcs at an $m$-vertex, see Fig.~1.
Almost balanced tangles
 considered up to non-rigid isotopy in $\{0<z<1\}$
 form a semigroup $NBT_n$ isomorphic to $NGT_n$.
The canonical isomorphism $\ph: NGT_n\to NBT_n$ is defined also
 by formulas~(24).

All three-page tangles
 considered up to non-rigid isotopy in $\{0<z<1\}$
 form a semigroup isomorphic to $NSG_n$, see Proposition~1 in
 Subsection~5.2.
Really, the relations~$\ph(11)-\ph(22)$ of $NBT_n$ follow from
 the relations (1)--(8) of $NSG_n$, see the proof of Lemma~3 in
 Subsection~6.3.
The new relations~$\ph(23')$ is reduced to~$(9')$ as follows
 (the case of a $2q$-vertex is completely similar to the case of a
 $(2p-1)$-vertex):
$$\begin{array}{l}
 \ph(\la_{2p-1,k} \sm_k^{-1})
  \st{(24)}{=}
 (d_2^{k-1} x_{2p-1,2} b_2^k) (d_2^{k} b_1 b_2 d_1 b_2^{k-1})
  \st{(2),(1')}{\s}
 d_2^{k-1} x_{2p-1,2} (d_2 d_0) (d_0 d_1) d_1 b_2^{k-1}
  \\ \st{(2)}{\s}
 d_2^{k-1} (x_{2p-1,2} b_2) (d_2^2 d_0^2 d_1^2) b_2^{k-1}
  \st{(9')}{\s}
 d_2^{k-1} (x_{2p-1,2} b_2) b_2^{k-1}
  \st{(24)}{=}
 \ph(\la_{2p-1,k}).
\end{array}\;\qed$$
\smallskip

Corollaries~1b and 2b are verified absolutely analogously
 to Corollaries~1a and 2a, respectively.
The proof of Corollary~3 is contained
 in the proofs of the corresponding results for $n$-graphs:
 we should replace the condition $3\leq m\leq n$ by $m\in J$.
Existence of a three-page embedding,
 for any non-rigid $J$-graph (Corollary~3), implies Corollary~4.
\smallskip

\emph{Proof of Corollary~4.}
Take a three-page embedding $G\subset\Y$.
Let $k_i$ be the number of the arcs of $G\cap P_i$.
Split the page $P_i$ into $k_i$ pages.
Let us move the arcs of $G\cap P_i$ to these new pages so that
 each page contains exactly one arc.
Consider a plane $\R^2$ orthogonal to the axis $\al$.
In $\R^3$ slightly deform the above arcs of $G$ so that
 their images under the projection $\R^3\to\R^2$ along $\al$
 are non-intersected loops.
This is a projection required.
\qed


\setcounter{section}{5}
\section{ Proof of Lemma~3}

In Claim 1 we deduce new word equivalences from the
 relations $(1)-(8)$ of the semigroup $RSG_n$.
Claims $1-4$ will imply Lemma~5 on a decomposition
 of any $i$-balanced word, see Subsection~6.2.
The relations $\ph(11)-\ph(23)$ reduce to the relations $(1)-(10)$
 by Lemma~5 and Claim~6 .
Subsection~6.3 finishes the proof of Lemma~3 by exploiting Claims $5-7$.
All relations in this section will be verified formally,
 but they have a clear geometric interpretation
 as well as the relations $(1)-(10)$ in Fig.~3.

\subsection{ New word equivalences in the semigroup $RSG_n$.}

Let $n\geq 2$ be fixed.
The commutativity $uv\s vu$ will be 
 denoted briefly by $u\lra v$.

\begin{claim}
The equivalences $(1)-(8)$ imply the following ones
 (where $i\in \Z_3$ and
 $$w_i\in 
 \B_{n,i}=\{  a_i, \; b_i, \; c_i, \; d_i, \; x_{m,i}, \; 
             b_{i-1} b_i d_{i-1}, \; b_{i-1} d_i d_{i-1} \vl 3\leq m\leq n \; \}):$$
$$\begin{array}{ll}
(25) \; b_i \s d_{i+1} d_{i-1}, \mbox{ or } b_{i+1} \s d_{i-1} d_i, \; b_{i-1} \s d_i d_{i+1}, & 
     \mbox{ or } b_0 \s d_1 d_2, \; b_1\s d_2 d_0, \; b_2 \s d_0 d_1; \\

(26) \; d_i \s b_{i-1} b_{i+1}, \mbox{ or } d_{i-1} \s b_{i+1} b_i, \; d_{i+1} \s b_i b_{i-1},  & 
     \mbox{ or } d_0 \s b_2 b_1, \; d_1\s b_0 b_2, \; d_2 \s b_1 b_0; \\

(27) \; d_{i+1} b_{i-1} \s b_{i-1} d_{i+1} t_i, \; b_{i+1} d_{i-1} \s t_i d_{i-1} b_{i+1},  & 
     \mbox{ where } t_i=b_{i+1} d_{i-1} d_{i+1} b_{i-1}; \\

(28) \; a_i \s a_{i-1} b_{i+1}, \; c_i \s d_{i+1} c_{i-1}; &
(29) \; a_i b_i \s a_{i-1} d_{i-1},  \; d_i c_i \s b_{i-1} c_{i-1}; \\

(30) \; b_i \s a_i b_i c_i, \; d_i \s a_i d_i c_i; &
\end{array} $$
$$\begin{array}{ll}
(31) \; b_i^{p-1} x_{2p-1,i} d_i^{p-1} \s x_{2p-1,i+1} b_{i+1}; &
(32) \; d_i^{p-1} x_{2p-1,i+1} b_{i+1} b_i^{p} \s x_{2p-1,i} b_{i}; \\

(33) \; x_{2q,i+1} \s d_{i-1} b_{i}^{q-2} x_{2q,i} d_{i}^{q-2} b_{i-1}; &
(34) \; b_{i}^{q-1} x_{2q,i} d_{i}^{q-1} \s d_{i+1} x_{2q,i+1} b_{i+1}; \\

(35) \; d_{i} c_{i} \lra w_{i+1}; &
(36) \; b_{i} c_{i} \lra w_{i-1}; \\

(37) \; a_{i} b_{i} \lra w_{i+1}; &
(38) \; a_{i} d_{i} \lra w_{i-1};  \\

(39) \; t_{i}, t_{i}' \lra w_{i}, \mbox{ where } &
     t_i=b_{i+1} d_{i-1} d_{i+1} b_{i-1}, \;
     t_i'=d_{i-1} b_{i+1} b_{i-1} d_{i+1};
\end{array} $$
$$\begin{array}{ll}
(40) \; x_{2p-1,i} b_{i} \lra w_{i+1}; &
(45) \; d_{i+1} b_{i-1} w_i d_{i-1} b_{i+1} \s
         b_{i-1} d_{i+1} w_i b_{i+1} d_{i-1}; \\

(41) \; d_i x_{2q,i} b_i \lra w_{i+1};  &
(46a) \;  b_{i-1}^2 a_i d_{i-1}^2 \s
           ( b_{i-1} a_i d_{i-1} ) d_i^2 ( b_{i-1} b_i d_{i-1} ) b_i; \\

(42) \;  b_i^{p-1} x_{2p-1,i} d_i^{p-1} \lra w_{i-1};  &
(46b) \;  b_{i-1}^2 b_i d_{i-1}^2 \s
           ( b_{i-1} b_i d_{i-1} ) d_i^2 ( b_{i-1} b_i d_{i-1} ) b_i; \\

(43) \; d_{i-1}^{p-1} x_{2p-1,i} b_i b_{i-1}^{p} \lra w_{i};  & 
(46c) \;  b_{i-1}^2 c_i d_{i-1}^2 \s
           d_i ( b_{i-1} d_i d_{i-1} ) b_i^2 ( b_{i-1} c_i d_{i-1} ); \\

(44) \; b_i^{q-1} x_{2q,i} d_i^{q-1} \lra w_{i-1}; &
(46d) \;  b_{i-1}^2 d_i d_{i-1}^2 \s
           d_i ( b_{i-1} d_i d_{i-1} ) b_i^2 ( b_{i-1} d_i d_{i-1} ); 
\end{array}$$
$$(46x) \; \left\{ \begin{array}{l}
         b_{i-1}^2 x_{2p-1,i} d_{i-1}^2 \s
          (b_{i-1} d_i d_{i-1}) d_i x_{2p-1,i} b_i^2 ( b_{i-1} b_i d_{i-1} )
          ( b_{i-1}^2 d_i d_{i-1}^2 ); \\

         b_{i-1}^2 x_{2q,i} d_{i-1}^2 \s
          ( b_{i-1}^2 b_i d_{i-1}^2 ) ( b_{i-1} d_i d_{i-1} )
          d_i^2 x_{2q,i} b_i^2 ( b_{i-1} b_i d_{i-1} )
          ( b_{i-1}^2 d_i d_{i-1}^2 ).
  \end{array} \right.$$
\end{claim}
\emph{Proof.}
The equivalences $(25)-(27)$ follow from $(1)-(2)$.
Due to (2) we have 
 $d_i\s b_i^{-1}$, 
 $b_{i-1} b_i d_{i-1} \s (b_{i-1} d_i d_{i-1})^{-1}$, and 
 $t_i'\s t_i^{-1}$.
Then (37) and $(39)-(41)$ follow from (8).
The other equivalences will be verified step by step exploiting
 the already checked ones.
Since $i\in \Z_3=\{0,1,2\}$, then we have $(i+1)+1=i-1$ and
 $(i-1)-1=i+1$.
$$\begin{array}{ll}
(28): \quad
  a_{i-1} b_{i+1} 
   \st{(3)}{\s} 
  (a_i d_{i+1}) b_{i+1}
   \st{(2)}{\s} a_i, & 

  d_{i+1} c_{i-1} 
   \st{(3)}{\s} 
  d_{i+1} (b_{i+1} c_i)
   \st{(2)}{\s} 
  c_i; \\

(29): \quad
  a_i b_i 
   \st{(28)}{\s} 
  (a_{i-1} b_{i+1}) b_i
   \st{(26)}{\s} 
  a_{i-1} d_{i-1},  & 

  d_i c_i 
   \st{(26)}{\s} 
  (b_{i-1} b_{i+1}) c_i
   \st{(3)}{\s} 
  b_{i-1} c_{i-1}; \\

(30): \quad
  a_i b_i c_i 
   \st{(29)}{\s} 
  a_{i} (d_{i+1} c_{i+1})
   \st{(3)}{\s} 
  a_{i-1} c_{i+1} 
   \st{(3)}{\s} b_{i},  & 

  a_i d_i c_i 
   \st{(29)}{\s} 
  a_{i} (b_{i-1} c_{i-1})
   \st{(28)}{\s} 
  a_{i+1} c_{i-1} 
   \st{(3)}{\s} 
  d_{i}; 
\end{array}$$
$$\begin{array}{l}
(31): \quad
  b_i^{p-1} x_{2p-1,i} d_i^{p-1}
   \st{(4)}{\s}
  b_i^{p-1} (d_i^{p-1} x_{2p-1,i+1} d_{i-1} b_i^{p-2}) d_i^{p-1}
   \st{(2)}{\s} \\ \qquad \st{(2)}{\s}
  x_{2p-1,i+1} (d_{i-1} d_i)
   \st{(25)}{\s}
  x_{2p-1,i+1} b_{i+1}; \\

(32): \quad
  d_i^{p-1} (x_{2p-1,i+1} b_{i+1}) b_i^{p}
   \st{(31)}{\s}
  d_i^{p-1} (b_i^{p-1} x_{2p-1,i} d_i^{p-1}) b_i^{p}
   \st{(2)}{\s}
  x_{2p-1,i} b_{i}; \\

(33): \quad
  d_{i-1} b_{i}^{q-2} x_{2q,i} d_{i}^{q-2} b_{i-1}
   \st{(4)}{\s}
  d_{i-1} b_{i}^{q-2} (d_{i}^{q-2} b_{i-1} x_{2q,i+1} d_{i-1} b_{i}^{q-2})
   d_{i}^{q-2} b_{i-1}
   \st{(2)}{\s}
  x_{2q,i+1}; \\
  
(34): \quad
  b_{i}^{q-1} x_{2q,i} d_{i}^{q-1}
   \st{(4)}{\s}
  b_{i}^{q-1} (d_i^{q-2} b_{i-1} x_{2q,i+1} d_{i-1} b_i^{q-2}) d_{i}^{q-1}
   \st{(2)}{\s}
  (b_i b_{i-1}) x_{2q,i+1} (d_{i-1} d_i)
   \st{(25)}{\s} \\ \qquad \st{(25)}{\s}
  (b_i b_{i-1}) x_{2q,i+1} b_{i+1}
   \st{(26)}{\s}
  d_{i+1} x_{2q,i+1} b_{i+1}.
\end{array}$$

Below in the proof of (35) we commute firstly $b_{i+1}$ with $d_i c_i$ and
 after that we use this equivalence to commute
 $a_{i+1}$ with $d_i c_i$.
$$\begin{array}{l}
(35b): \quad
  b_{i+1} (d_i c_{i}) 
   \st{(30)}{\s} 
  (a_{i+1} b_{i+1} c_{i+1}) (d_i c_i)
   \st{(7)}{\s} 
  a_{i+1} b_{i+1} (d_i c_i) c_{i+1}
   \st{(26)}{\s} \\ \qquad \st{(26)}{\s} 
  a_{i+1} b_{i+1} (b_{i-1} b_{i+1}) c_i c_{i+1}
   \st{(3)}{\s} 
  (a_{i+1} b_{i+1}) b_{i-1} c_{i-1} c_{i+1}
   \st{(8)}{\s} 
  b_{i-1} c_{i-1} (a_{i+1} b_{i+1}) c_{i+1}
   \st{(30)}{\s} \\ \qquad \st{(30)}{\s} 
  b_{i-1} c_{i-1} b_{i+1}
   \st{(3)}{\s} 
  b_{i-1} (b_{i+1} c_i) b_{i+1}
   \st{(26)}{\s}  
  (d_{i} c_i) b_{i+1}; \\ \\

(35a): \quad
  a_{i+1} (d_i c_{i}) 
   \st{(3)}{\s} 
  (a_{i-1} d_{i}) (d_i c_{i})
   \st{(26)}{\s} 
  a_{i-1} (b_{i-1} b_{i+1}) (d_i c_{i})
   \st{(35b)}{\s} 
  a_{i-1} b_{i-1} (d_i c_i) b_{i+1}
   \st{(37)}{\s} \\ \qquad \st{(37)}{\s} 
  (d_{i} c_{i}) (a_{i-1} b_{i-1}) b_{i+1}
   \st{(26)}{\s} 
  (d_{i} c_{i}) (a_{i-1} d_{i})
   \st{(3)}{\s} 
  (d_{i} c_{i}) a_{i+1}. 
\end{array} $$

The remaining equivalences in $(35)$ follow from $(35a),(35b)$, and (7).
The equivalences (36) are easily proved by (29) and (35), 
 as well as (38) by (29) and (37), 
 as well as (42) by (31) and (40), 
 as well as (43) by (32) and (40), 
 as well as (44) by (34) and (41).
The last calculations are straightforward:

$$\begin{array}{l}
(45): \;
  d_{i+1} b_{i-1} w_i d_{i-1} b_{i+1} 
   \st{(27)}{\s}
  (b_{i-1} d_{i+1} t_i) w_i d_{i-1} b_{i+1}
   \st{(39)}{\s} \\ \qquad \st{(39)}{\s} 
  b_{i-1} d_{i+1} (w_i t_i) d_{i-1} b_{i+1}
   \st{(27)}{\s} 
  b_{i-1} d_{i+1} w_i b_{i+1} d_{i-1}; \\ \\

(46a): \;
 b_{i-1}^2 a_i d_{i-1}^2
  \st{(2)}{\s}
 b_{i-1}^2 a_i (d_i b_i) d_{i-1}^2
  \st{(38)}{\s}
 b_{i-1} (a_i d_{i}) (b_{i-1} b_{i}) d_{i-1}^2
  \st{(26)}{\s}
 b_{i-1} a_{i} d_{i} (b_{i-1} b_{i}) d_{i-1} (b_{i+1} b_i)
  \\ \qquad \st{(39)}{\s}
 b_{i-1} a_i b_{i+1} (d_{i} b_{i-1} b_i d_{i-1}) b_i
  \st{(25)}{\s}
 (b_{i-1} a_i d_{i-1}) d_{i}^2 (b_{i-1} b_i d_{i-1}) b_i; \\ \\

(46b): \;
 b_{i-1}^2 b_i d_{i-1}^2
  \st{(2)}{\s}
 b_{i-1} (b_i d_i) b_{i-1} b_i d_{i-1}^2
  \st{(26)}{\s}
 b_{i-1} b_i (d_i b_{i-1} b_i d_{i-1}) (b_{i+1} b_i)
  \st{(39)}{\s} \\ \qquad \st{(39)}{\s}
 b_{i-1} b_i b_{i+1} (d_i b_{i-1} b_i d_{i-1}) b_i
  \st{(25)}{\s}
 (b_{i-1} b_i d_{i-1}) d_i^2 (b_{i-1} b_i d_{i-1}) b_i; 
\end{array}$$

$$\begin{array}{l}
(46c): \;
 b_{i-1}^2 c_i d_{i-1}^2
  \st{(2)}{\s}
 b_{i-1}^2 (d_i b_i) c_i d_{i-1}^2
  \st{(36)}{\s}
 b_{i-1}^2 d_{i} d_{i-1} (b_{i} c_i) d_{i-1}
  \st{(25)}{\s}
 (d_{i} d_{i+1}) (b_{i-1} d_i d_{i-1} b_{i}) c_i d_{i-1}
 \\ \qquad \st{(39)}{\s}
 d_{i} (b_{i-1} d_i d_{i-1} b_{i}) d_{i+1} c_i d_{i-1}
  \st{(26)}{\s}
 d_{i} (b_{i-1} d_i d_{i-1}) b_{i}^2 (b_{i-1} c_i d_{i-1}); \\ \\

(46d): \;
 b_{i-1}^2 d_i d_{i-1}^2
  \st{(2)}{\s}
 b_{i-1}^2 d_i d_{i-1} (b_i d_i) d_{i-1}
  \st{(25)}{\s}
 (d_i d_{i+1}) b_{i-1} d_i d_{i-1} b_i d_{i} d_{i-1}
  \st{(39)}{\s} \\ \qquad \st{(39)}{\s}
 d_i (b_{i-1} d_i d_{i-1} b_i) d_{i+1} d_{i} d_{i-1}
  \st{(26)}{\s}
 d_i (b_{i-1} d_i d_{i-1}) b_i^2 (b_{i-1} d_{i} d_{i-1}); 
\end{array}$$

$$\begin{array}{l}
(46x): \;
 b_{i-1}^2 x_{2p-1,i} d_{i-1}^2
  \st{(2)}{\s}
 b_{i-1}^2 x_{2p-1,i} (b_i b_{i+1}^2 d_{i+1}^2 d_i) d_{i-1}^2
  \st{(40)}{\s} 
 b_{i-1}^2 b_{i+1}^2 (x_{2p-1,i} b_i) d_{i+1}^2 d_i d_{i-1}^2
  \\ \qquad \st{(26)}{\s}
 b_{i-1} (d_i d_{i+1}) b_{i+1}^2 x_{2p-1,i} b_i (b_i b_{i-1})^2 d_i d_{i-1}^2
  \st{(2)}{\s}
 b_{i-1} d_i b_{i+1} x_{2p-1,i} b_i^2 (b_{i-1} b_i) b_{i-1} d_i d_{i-1}^2
 \\ \qquad \st{(25),(2)}{\s}
 (b_{i-1} d_i d_{i-1}) d_i x_{2p-1,i} b_i^2 ( b_{i-1} b_i d_{i-1} )
  ( b_{i-1}^2 d_i d_{i-1}^2 ); 
\end{array}$$

$$\begin{array}{l}
(46x): \;
 b_{i-1}^2 x_{2q,i} d_{i-1}^2
  \st{(2)}{\s}
 b_{i-1}^2 (b_i d_i) x_{2q,i} (b_i b_{i+1}^2 d_{i+1}^2 d_i) d_{i-1}^2
  \st{(41)}{\s}
 b_{i-1}^2 b_i b_{i+1}^2 (d_i x_{2q,i} b_i) d_{i+1}^2 d_i d_{i-1}^2
  \\ \qquad \st{(25)}{\s} 
 b_{i-1}^2 b_i (d_{i-1} d_i)^2 d_i x_{2q,i} b_i d_{i+1}^2 d_i d_{i-1}^2
  \st{(26)}{\s} 
 (b_{i-1}^2 b_i d_{i-1}) d_i d_{i-1} d_i^2 x_{2q,i} b_i (b_i b_{i-1})^2 d_i d_{i-1}^2
  \\ \qquad \st{(2)}{\s}
 ( b_{i-1}^2 b_i d_{i-1}^2 ) ( b_{i-1} d_i d_{i-1} )
   d_i^2 x_{2q,i} b_i^2 ( b_{i-1} b_i d_{i-1} ) ( b_{i-1}^2 d_i d_{i-1}^2 ). \quad\qed
 \end{array} $$
\smallskip

\subsection{ Decomposition of $i$-balanced words}

\begin{claim}
For each $i\in \Z_3$, every $i$-balanced word 
 is equivalent by (1)--(8), (25)--(46) to an $i$-balanced word
 containing only the following letters:
 $a_i$, $b_i$, $c_i$, $d_i$, $x_{m,i}$, $b_{i-1}$, $d_{i-1}$.
\end{claim}
\emph{Proof.}
The rest letters could be eliminated by using the following substitutions:
$$\left\{ \begin{array}{lll}
  x_{2p-1,i+1} \st{(4)}{\s} d_{i+1}^{p-1} x_{2p-1,i-1} d_{i} b_{i+1}^{p-2}, &
  x_{2p-1,i-1} \st{(4)}{\s} d_{i-1}^{p-1} x_{2p-1,i} d_{i} b_{i-1}^{p-2}, &
  a_{i+1} \st{(3)}{\s} a_{i-1} d_i, \\

  x_{2q,i-1} \st{(4)}{\s} d_{i-1}^{q-2} b_{i+1} x_{2q,i} d_{i+1} b_i^{q-2}, &
  x_{2q,i+1} \st{(33)}{\s} d_{i-1} b_i^{q-2} x_{2q,i} d_i^{q-2} b_{i-1}, &
  c_{i+1} \st{(3)}{\s} b_{i} c_{i-1}, \\
 
  a_{i-1} \st{(3)}{\s} a_i d_{i+1}, \;
  c_{i-1} \st{(3)}{\s} b_{i+1} c_i, &
  b_{i+1} \st{(25)}{\s} d_{i-1} d_i, \;
  d_{i+1} \st{(26)}{\s} b_i b_{i-1}. & \qed
\end{array} \right. $$
\smallskip

Fix an index $i\in\Z_3$.
\smallskip

{\bf Definition 18} (\emph{the encoding $\mu(w)$, the depth $d(w)$}).
Let $w$ be an $i$-balanced word in the letters
 $a_i$, $b_i$, $c_i$, $d_i$, $x_{m,i}$, $b_{i-1}$, $d_{i-1}$.
Let us consider the following substitution 
 $\mu: a_i,b_i,c_i,d_i,x_{m,i}\to \bu$;
 $b_{i-1}\to ($; $d_{i-1}\to )$.
By $\mu(w)$ denote the resulting \emph{encoding} consisting of
 brackets and bullets.  
Since $w$ is $i$-balanced,
 then the encoding $\mu(w)$ without bullets
 is a balanced bracket expression, see Definition~9 in
 Subsection~3.4.
For each place $k$, denote by $dif(k)$
 the difference between the number of the left and right brackets in
 the subword of $\mu(w)$ ending at this place.
The maximum of $dif(k)$, for all $k$, is called
 \emph{the depth $d(w)$} of $w$.
For the word $w=b_{i-1}^2 a_i d_{i-1}^2$, we have $\mu(w)=((\bu))$
 and $d(w)=2$.
\ed
\medskip

{\bf Definition 19} (\emph{star decomposable words}).
\emph{A star of the deep $k$} is
 the encoding of the type $(^k \bu )^k$, which has $k$ couples of brackets.
If the encoding $\mu(w)$ decomposes on
 several stars, then
 $w$ is called \emph{star decomposable}.
In this case, the depth $d(w)$ is the maximum among
 the depths of all stars participating in the star decomposition.
\ed
\smallskip

\begin{claim}
Every $i$-balanced word $w$
 is equivalent to a star decomposable word $w'$
 of the same depth $d(w')=d(w)$.
\end{claim}
\begin{proof}
Let us look at the beginning of the encoding $\mu(w)$.
After several initial left brackets
 the encoding $\mu(w)$ contains either
 a right bracket or a bullet.
In the first case, let us delete
 the couple of brackets $()$ by the relation
 $b_{i-1} d_{i-1}\st{(2)}{\s} \es$.
Hence we may assume that
 the next symbol after $k$ left brackets is a bullet.
Since $\mu(w)$ is balanced, then after this bullet
 it may be the sequence of $j$, $0\leq j\leq k$, right brackets.
If $j<k$, then insert into $w$ the subword
 $d_{i-1}^{k-j} b_{i-1}^{k-j}\st{(2)}{\s} \es$
 after the last right bracket.
This operation does not change the depth $d(w)$.
Then, in the resulting word $w_1$,
 the encoding $\mu(w_1)$ contains a star of the depth $k$
 at the beginning.
For instance, starting with the word $w=b_{i-1} a_i^2 d_{i-1}$,
 we get $w_1=b_{i-1} a_i d_{i-1} b_{i-1} a_i d_{i-1}$
 with $\mu(w_1)=(\bu)(\bu)$.
Continuing this process, in finitely many steps, we shall get
 a star decomposable word $w_N$ of the same depth $d(w_N)=d(w)$.
\end{proof}

For any letter $s$, we denote by $s'$ the word $b_{i-1} s d_{i-1}$,
 for example, $a_i'=b_{i-1} a_i d_{i-1}$.

\begin{claim}
Every star decomposable word $w$ is equivalent to a word
 decomposed on the following $i$-balanced subwords\footnote{
 As usual, we suppose that $i\in\Z_3$ and $3\leq m\leq n$,
  the parameters $i,m$ are fixed.}:
 $\{ a_i, \; b_i, \; c_i, \; d_i, \; x_{m,i}, \;
     a_i', \; b_i', \; c_i', \; d_i', \; x_{m,i}'\}$.
\end{claim}
\begin{proof}
Induction on the depth $d(w)$.
Base $d(w)=1$ is trivial.
Suppose the encoding $\mu(w)$
 contains stars of a depth $k\geq 2$.
Apply one of the following moves
 to every such star.

$$\left\{ \begin{array}{l}
 u=b_{i-1}^2 a_i d_{i-1}^2 \st{(46a)}{\s} a_i' d_i^2 b_i' b_i = v, 
 \mbox{ i.e. } \mu(u)=((\bu)) \to \mu(v)=(\bu)\bu\bu(\bu)\bu; \\

 u=b_{i-1}^2 b_i d_{i-1}^2 \st{(46b)}{\s} b_i' d_i^2 b_i' b_i = v,  
 \mbox{ i.e. } \mu(u)=((\bu)) \to \mu(v)=(\bu)\bu\bu(\bu)\bu; \\

 u=b_{i-1}^2 c_i d_{i-1}^2 \st{(46c)}{\s} d_i d_i' b_i^2 c_i' = v, 
 \mbox{ i.e. } \mu(u)=((\bu)) \to \mu(v)=\bu(\bu)\bu\bu(\bu); \\

 u=b_{i-1}^2 d_i d_{i-1}^2 \st{(46d)}{\s} d_i d_i' b_i^2 d_i' = v,
 \mbox{ i.e. } \mu(u)=((\bu)) \to \mu(v)=\bu(\bu)\bu\bu(\bu); 
\end{array} \right.$$

$$\left\{ \begin{array}{l}
  u=b_{i-1}^2 x_{2p-1,i} d_{i-1}^2
   \st{(46x)}{\s}
   (b_{i-1} d_i d_{i-1}) d_i x_{2p-1,i} b_i^2 ( b_{i-1} b_i d_{i-1} )
   ( b_{i-1}^2 d_i d_{i-1}^2 )
   \st{(46d)}{\s} \\ 
  b_i' d_i x_{2p-1,i} b_i^2 b_i' (d_i d_i' b_i^2 d_i') = v,  \mbox{ i.e.  }
   \mu(u)=((\bu)) \to 
   \mu(v)=(\bu) \bu \bu \bu \bu (\bu) \bu (\bu) \bu \bu (\bu); \\

  u=b_{i-1}^2 x_{2q,i} d_{i-1}^2
   \st{(46x)}{\s}
   ( b_{i-1}^2 b_i d_{i-1}^2 ) ( b_{i-1} d_i d_{i-1} )
   d_i^2 x_{2q,i} b_i^2 ( b_{i-1} b_i d_{i-1} ) ( b_{i-1}^2 d_i d_{i-1}^2 )
   \st{(46b,46d)}{\s} \\ 
  \st{(46b,46d)}{\s} 
   (b_i' d_i^2 b_i' b_i) d_i' d_i^2 x_{2q,i} b_i^2 b_i' (d_i d_i' b_i^2 d_i') = v,  
   \mbox{ i.e.  } \\
  \mu(u)=((\bu)) \to 
   \mu(v)=(\bu) \bu \bu (\bu) \bu (\bu) \bu \bu \bu \bu \bu (\bu) \bu (\bu) \bu \bu (\bu).
\end{array} \right.  $$

By using the above moves we get a word $w_1\s w$ of the
 depth $d(w_1)=d(w)-1$.
By Claim~3 the word $w_1$ is equivalent to 
 a star decomposable word $w_2$ of the depth $d(w_2)=d(w_1)=d(w)-1$.
The claim follows from the induction hypothesis for $w_2$.
\end{proof}

\begin{lemma}
For each $i\in \Z_3$, 
 every $i$-balanced word from $W_{n,i}$ is equivalent
 to a word decomposed on $i$-balanced words from the set
 $\B_{n,i}$, see this notation in Claim~1.
\end{lemma}
\emph{Proof.}
By Claims 3 and 4 it remains to eliminate the following words:
$$\begin{array}{l}
   a'_i=b_{i-1} a_i d_{i-1}                
    \st{(25),(2)}{\s}
   (d_i d_{i+1}) a_i (b_i d_i) d_{i-1} 
    \st{(37)}{\s}
   d_i (a_i b_i) d_{i+1} d_i d_{i-1} 
    \st{(26)}{\s}
   d_i a_i b_i^2 (b_{i-1} d_i d_{i-1}); \\

   c'_i=b_{i-1} c_i d_{i-1}                
    \st{(26),(2)}{\s}
   b_{i-1} (b_i d_i c_i) (b_{i+1} b_i) 
    \st{(35)}{\s}
   b_{i-1} b_i b_{i+1} (d_i c_i) b_i 
    \st{(25)}{\s}
   (b_{i-1} b_i d_{i-1}) d_i^2 c_i b_i; 
\end{array}$$
$$\left\{ \begin{array}{l}
   x'_{2p-1,i}=b_{i-1} x_{2p-1,i} d_{i-1}
    \st{(2)}{\s}
   b_{i-1} x_{2p-1,i} (b_i d_i) d_{i-1}
    \st{(25)}{\s}
   (d_i d_{i+1}) (x_{2p-1,i} b_i) d_i d_{i-1}
    \st{(40)}{\s} \\ \st{(40)}{\s}
   d_i (x_{2p-1,i} b_i) d_{i+1} d_i d_{i-1}
    \st{(26)}{\s}
   d_i x_{2p-1,i} b_i^2 (b_{i-1} d_i d_{i-1});
\end{array} \right. $$
$$\left\{ \begin{array}{l}
   x'_{2q,i}=b_{i-1} x_{2q,i} d_{i-1}
    \st{(2)}{\s}
   b_{i-1} (b_i d_i) x_{2q,i} (b_i d_i) d_{i-1}
    \st{(2)}{\s}
   b_{i-1} b_i (b_{i+1} d_{i+1}) (d_i x_{2q,i} b_i) d_i d_{i-1}
    \st{(41)}{\s} \\ \st{(41)}{\s}
   b_{i-1} b_i b_{i+1} (d_i x_{2q,i} b_i) d_{i+1} d_i d_{i-1}
    \st{(25),(26)}{\s}
   (b_{i-1} b_i d_{i-1}) d_i^2 x_{2q,i} b_i^2 (b_{i-1} d_i d_{i-1}). \qed
\end{array} \right. $$
\smallskip

In the equivalences $(35)-(45)$, let us replace
 the condition $w_i\in \B_{n,i}$ by $w_i\in W_{n,i}$.
The obtained relations will be denoted by $(35')-(45')$.

\begin{claim}
The relations $(35')-(45')$ hold for
 arbitrary $i$-balanced words $w_i\in W_{n,i}$.
\end{claim}
\begin{proof}
By Lemma 5 each $i$-balanced word $w\in W_{n,i}$ decomposes on
 the $i$-balanced words of $\B_{n,i}$.
Since the commutative equivalences $(35)-(45)$ hold 
 for the words of $\B_{n,i}$ by Claim~1,
 then they also hold for $w_i\in W_{n,i}$.
\end{proof}


\subsection{ Deduction of the relations $\ph(11)-\ph(23)$ from
 the relations $(1)-(10)$}

The relations $\ph(11)-\ph(23)$ between words in $\A_n$
 were obtained from the relations $(11)-(23)$ of $RGT_n$
 under the isomorphism $\ph: RGT_n\to RBT_n$, see Subsection~5.1.
For each $l\geq 1$, let us denote by $u_l$ any symbol of the set
 $\{\xi_l,\eta_l,\sm_l,\sm_l^{-1},\la_{m,l} \vl 3\leq m\leq n\}$.

For $k\geq 1$, define \emph{the shift maps}
 $\te_k: RGT_n\to RGT_n$ and $\om_k: RBT_n\to RBT_n$ by
 $\te_k(u_l)=u_{k+l}$ and $\om_k(w)=d_2^k w b_2^k$.
Then $\te_k$ is a well-defined homomorphism.
Indeed, any relation of $(11)-(23)$, for each $k>1$,
 is the shift image of the corresponding relation for $k=1$.
For example, the relation $\xi_k \xi_l= \xi_{l+2} \xi_k$
 is obtained from $\xi_1 \xi_{l-k+1}= \xi_{l-k+3} \xi_1$
 under the shift $\te_{k-1}$.
By the relations (2) the shift $\om_k$ sends
 equivalent words to equivalent ones, i.e. 
 $\om_k$ is also a homomorphism.
Moreover, the following diagram is commutative.
$$\begin{CD}
RGT_n      @>{\te_k}>> RGT_n \\
@V{\ph}VV              @VV{\ph}V \\
RBT_n      @>{\om_k}>> RBT_n
\end{CD}$$

\begin{claim}
All the relations $\ph(11)-\ph(23)$ reduce to the relations
 $\ph(11)-\ph(23)$ with $k=1$ by using the equivalence
 $b_2d_2\s 1\s d_2b_2$ of (2).
\end{claim}
\emph{Proof} follows from the commutativity
 of the above diagram. 
For instance, we have
$$\begin{array}{l}
 \ph(\xi_k\xi_l)=
 \ph\circ\te_{k-1}(\xi_1\xi_{l-k+1}) =
 \om_{k-1}\circ\ph(\xi_1\xi_{l-k+1}) =
 d_2^{k-1} \ph(\xi_1\xi_{l-k+1}) b_2^{k-1} \st{\ph(11),k=1}{\s} \\

 d_2^{k-1} \ph(\xi_{l-k+3}\xi_1) b_2^{k-1} =
 \om_{k-1}\circ\ph(\xi_{l-k+3}\xi_1) =
 \ph\circ\te_{k-1}(\xi_{l-k+3}\xi_1) =
 \ph(\xi_{l+2}\xi_k).
\end{array}\eqno{\qed}$$
\medskip

\begin{claim}
Under the map $\ph: RGT_n\to RBT_n$
 the relations $(1)-(10), (25)-(34), (35')-(45')$
 imply the following ones (see $\Sm$-notations in Lemma~1 and 
 $D$-notations in (9)):
$$\begin{array}{lll}
(47) \; \ph(\Sm_{1,l})     \s b_1 d_2^l d_1 b_2^{l}; &
(48) \; \ph(\bar\Sm_{k,l}) \s d_2^{k-1} b_1^l d_2 d_1^l b_2^k; &
(49) \; \ph(\Sm'_{1,l})    \s D_{l+1,2}.
\end{array}$$
\end{claim}
\emph{Proof.}
The following calculations are straightforward:

$$\begin{array}{l}
 (47): \;
  \ph(\Sm_{1,l})
   =
  \ph(\sm_1) \ldots \ph(\sm_{l})
   \st{(24)}{\s}
  (b_1 d_2 d_1 b_2) (d_2 b_1 d_2 d_1 b_2^{2}) \ldots
   (d_2^{l-1} b_1 d_2 d_1 b_2^{l}) \\ \qquad \st{(2)}{\s} 
  (b_1 d_2 d_1) (b_1 d_2 d_1) \ldots (b_1 d_2 d_1 b_2^{l})
  \st{(2)}{\s} 
   b_1 d_2^l d_1 b_2^{l};
 \end{array} $$

$$\begin{array}{l}
 (48_1): \;
 \ph(\bar\Sm_{1,l})
  =
 \ph(\sm_{l}) \ldots \ph(\sm_{2}) \ph(\sm_1)
  \st{(24)}{\s}
 (d_2^{l-1} b_1 d_2 d_1 b_2^{l}) (d_2^{l-2} b_1 d_2 d_1 b_2^{l-1})
   \ldots (b_1 d_2 d_1 b_2)
  \\ \qquad \st{(2)}{\s}  
 d_2^{l-1} (b_1 d_2 d_1 b_2^2)^{l-2} (b_1 d_2 d_1 b_2) (b_2 b_1) d_2 d_1 b_2
  \st{(39')}{\s} \\ \qquad
 d_2^{l-1} (b_1 d_2 d_1 b_2^2)^{l-2} (b_2 b_1) (b_1 d_2 d_1 b_2) d_2 d_1 b_2
  \st{(2)}{\s} 
 d_2^{l-1} (b_1 d_2 d_1 b_2^2)^{l-3} (b_1 d_2 d_1 b_2) (b_2^2 b_1^2) d_2 d_1^2 b_2
  \\ \qquad \st{(39'),(2)}{\s} \cdots \st{(39'),(2)}{\s}
  d_2^{l-1} (b_2^{l-1} b_1^l) d_2 d_1^l b_2
  \st{(2)}{\s} 
 b_1^l d_2 d_1^l b_2;
 \end{array}$$

$$(48_k): \quad
 \ph(\bar\Sm_{k,l})
  =
 \ph( \te_{k-1}( \bar\Sm_{1,l} ) )
  =
 \om_{k-1}( \ph( \bar\Sm_{1,l} ) )
  =
 d_2^{k-1} \ph( \bar\Sm_{1,l} ) b_2^{k-1}
  \st{(48_1)}{\s}
 d_2^{k-1} b_1^l d_2 d_1^l b_2^k; $$

$$\begin{array}{l}
 (49): \;
 \ph(\Sm'_{1,l})
  = 
 \ph(\sm_{l}^{-1}) \ldots \ph(\sm_{2}^{-1} \cdots \sm_{l}^{-1})
  \ph(\sm_{1}^{-1} \cdots \sm_{l}^{-1})
  =
 \ph(\bar\Sm_{l,1}^{-1}) \ldots \ph(\bar\Sm_{2,l-1}^{-1})
  \ph(\bar\Sm_{1,l}^{-1})
  \\ \qquad \st{(48)}{\s} 
 (d_2^{l} b_1 b_2 d_1 b_2^{l-1}) (d_2^{l-1} b_1 b_2 d_1 b_2^{l-2}) \ldots
  (d_2^{2} b_1^{l-1} b_2 d_1^{l-1} b_2) (d_2 b_1^l b_2 d_1^l)
  \st{(2)}{\s} \\ \qquad  \st{(2)}{\s}
 d_2^{l} b_1 (b_2 b_1)^{l-1} b_2 d_1^l
  \st{(26)}{\s}
 d_2^{l} b_1 d_0^{l-1} b_2 d_1^l
  \st{(25)}{\s}
 d_2^{l+1} d_0^{l+1} d_1^{l+1}=D_{l+1,2}. \quad\qed 
 \end{array}$$

{\bf Proof of Lemma 3 from Subsection 5.2.}
Here we deduce the relations $\ph(11)-\ph(23)$
 from the equivalences $(1)-(10), (25)-(34), (35')-(45'), (47)-(49)$.
By $\star$ we denote the following images under the map
 $\ph: RGT_n\to RBT_n$ (see (24) in Subsection~5.1):
 $$(24_1)\quad \left\{ \begin{array}{lll}
   \ph(\xi_1)=d_2 c_2, &
   \ph(\sm_1)=b_1 d_2 d_1 b_2, &
   \ph(\la_{2p-1,1})=x_{2p-1,2} b_2, \\

   \ph(\eta_1)=a_2 b_2, &
   \ph(\sm_1^{-1})=d_2 b_1 b_2 d_1, &
   \ph(\la_{2q,1})=d_2 x_{2q,2} b_2.
 \end{array} \right. $$
The words $\ph(u_l)=d_2^{l-1} \star b_2^{l-1}$
 are 1-balanced (see Fig.~6), i.e. 
 $d_2^{l} \star b_2^{l}\in W_{n,1}$ for each $l\in\N$.
Then the relations $\ph(11)-\ph(14)$ could be proved
 by the following patterns:


$$\begin{array}{l}
(11): \quad
  \ph(\xi_1 u_l)
   \st{(2)}{\s}
  d_2^2 (b_2 c_2) (d_2^{l-1} \star b_2^{l-1})
   \st{(36')}{\s}
  d_2^2 (d_2^{l-1} \star b_2^{l-1}) (b_2 c_2)
   \st{(2)}{\s}
  \ph(u_{l+2} \xi_1); \\

(12): \quad
  \ph(\eta_1 u_l)
   \st{(2)}{\s}
  (a_2 d_2) (d_2^{l-3} \star b_2^{l-3}) b_2^2
   \st{(38')}{\s}
  (d_2^{l-3} \star b_2^{l-3}) (a_2 d_2) b_2^2
   \st{(2)}{\s}
  \ph(u_{l-2} \eta_1); \\

(13): \quad
  \ph(\sm_1 u_l)
   \st{(2)}{\s}
  (b_1 d_2 d_1) (d_2^{l-2} \star b_2^{l-1})
   \st{(25),(2)}{\s}
  d_2^2 (b_2 d_0 d_2 b_0) (d_2^{l-3} \star b_2^{l-3}) b_2^2
   \st{(39')}{\s} \\ \qquad \st{(39')}{\s}
  d_2^2 (d_2^{l-3} \star b_2^{l-3}) (b_2 d_0 d_2 b_0) b_2^2
   \st{(2),(26)}{\s}
  (d_2^{l-1} \star b_2^{l-2}) (b_2 b_1) d_2 d_1 b_2
   \st{(24)}{=}
  \ph(u_{l} \sm_1);
\end{array}$$

$$\begin{array}{l}
 (14): \quad
  \ph(\la_{2p-1,1} u_l)
   \st{(2)}{\s}
  d_2^{p-1} (b_2^{p-1} x_{2p-1,2} d_2^{p-1})
   (d_2^{l-p-1} \star b_2^{l-p-1}) b_2^{p}
   \\ \qquad \st{(42')}{\s}
  d_2^{p-1} (d_2^{l-p-2} \star b_2^{l-p-2})
   (b_2^{p-1} x_{2p-1,2} d_2^{p-1}) b_2^{p}
   \st{(2)}{\s}
  \ph(u_{l-1} \la_{2p-1,1}); \\

 (14): \quad
  \ph(\la_{2q,1} u_l)
   \st{(2)}{\s}
  d_2^{q} (b_2^{q-1} x_{2q,2} d_2^{q-1})
   (d_2^{l-q-1} \star b_2^{l-q-1}) b_2^{q}
  \\ \qquad \st{(44')}{\s}
  d_2^{q} (d_2^{l-q-1} \star b_2^{l-q-1})
   (b_2^{q-1} x_{2q,2} d_2^{q-1}) b_2^{q}
   \st{(2)}{\s}
  \ph(u_{l} \la_{2q,1}).
\end{array}$$
\smallskip

The remaining calculations are straightforward:


$$(15): \;
  \ph(\eta_{2} \xi_1)
   \st{(24)}{=}
  (d_2 a_2 b_2^{2}) (d_2 c_2)
   \st{(2)}{\s}
  d_2 (a_2 b_2 c_2)
   \st{(30)}{\s}
  d_2 b_2
   \st{(2)}{\s}
  1
   \st{(2),(30)}{\s}
  (a_2 d_2 c_2) b_2
   \st{(2)}{\s}
  \ph(\eta_1 \xi_{2});$$


$$\begin{array}{l}
 (16): \quad
  \ph(\eta_{3} \sm_{2} \xi_1)
   \st{(24)}{=}
  (d_2^{2} a_2 b_2^3) (d_2 b_1 d_2 d_1 b_2^2) (d_2 c_2)
   \st{(2)}{\s}
  d_2^{2} a_2 b_2 (b_2 b_1) d_2 d_1 (b_2 c_2)
   \\ \qquad \st{(26)}{\s}
  d_2^{2} (a_2 b_2) d_0 d_2 d_1 (b_2 c_2)
   \st{(37)}{\s} 
  d_2^{2} d_0 (a_2 b_2) d_2 d_1 (b_2 c_2)
   \st{(36)}{\s} 
  d_2 (d_2 d_0) (a_2 b_2) d_2 (b_2 c_2) d_1
   \\ \qquad \st{(2)}{\s}
  d_2 (d_2 d_0) (a_2 b_2 c_2) d_1
   \st{(25)}{\s} 
  d_2 b_1 (a_2 b_2 c_2) d_1
   \st{(30)}{\s}
  d_2 b_1 b_2 d_1
   \st{(24_1)}{=}
  \ph(\sm_1^{-1}),
 \end{array}$$

$$\begin{array}{l}
 (16): \quad
  \ph(\eta_1 \sm_{2} \xi_{3})
   \st{(24)}{=}
  (a_2 b_2) (d_2 b_1 d_2 d_1 b_2^2) (d_2^3 c_2 b_2^2)
   \st{(2)}{\s}
  a_2 b_1 d_2 (d_1 d_2) c_2 b_2^{2}
   \\ \qquad \st{(25)}{\s}
  a_2 b_1 d_2 (b_0 c_2) b_2^{2}
   \st{(3)}{\s}
  (a_0 d_1) b_1 d_2 c_1 b_2^{2}
   \st{(2)}{\s}
  a_0 d_2 (b_1 d_1) c_1 b_2^{2}
   \st{(35)}{\s} 
  a_0 d_2 b_1 b_2 (d_1 c_1) b_2
   \\ \qquad \st{(39)}{\s}
  (d_2 b_1 b_2 d_1) (a_0 c_1) b_2
   \st{(3)}{\s}
  d_2 b_1 b_2 d_1 (d_2 b_2)
   \st{(2)}{\s}
  d_2 b_1 b_2 d_1
   \st{(24_1)}{=}
  \ph(\sm_1^{-1});
 \end{array}$$


$$\begin{array}{l}
 (17): \quad
  \ph(\eta_{p} \la_{2p-1,2} \xi_1)
   \st{(2)}{\s}
  d_2^{p-1} (a_2 b_2^{p-1} x_{2p-1,2} b_2 c_2)
   \st{(5)}{\s}
  d_2^{p-1} (b_2^{p-1} x_{2p-1,2} b_2)
   \st{(2)}{\s}
  \ph(\la_{2p-1,1}), \\ \qquad

  \ph(\eta_{1} \la_{2p-1,2} \xi_{p+1})
   \st{(2)}{\s}
  (a_2 x_{2p-1,2} d_2^{p-1} c_2) b_2^{p}
   \st{(5)}{\s}
  d_2^{2} (x_{2p-1,2} d_2^{p-1}) b_2^{p}
   \st{(2)}{\s}
  \ph(\la_{2p-1,1});
\end{array}$$


$$\begin{array}{l}
 (18): \quad
  \ph(\eta_{q+1} \la_{2q,2} \xi_1)
   \st{(2)}{\s}
  d_2^{q} (a_2 b_2^{q-1} x_{2q,2} b_2 c_2)
   \st{(6)}{\s}
  d_2^{q} (b_2^{q-1} x_{2q,2} b_2)
   \st{(2)}{\s}
  d_2 x_{2q,2} b_2
   \st{(24_1)}{=}
  \ph(\la_{2q,1}), \\ \qquad

  \ph(\eta_{1} \la_{2q,2} \xi_{q+1})
   \st{(2)}{\s}
  (a_2 d_2 x_{2q,2} d_2^{q-1} c_2) b_2^{q}
   \st{(6)}{\s}
  (d_2 x_{2q,2} d_2^{q-1}) b_2^{q}
   \st{(2)}{\s}
  d_2 x_{2q,2} b_2
   \st{(24_1)}{=}
  \ph(\la_{2q,1});
\end{array}$$


$$\begin{array}{l}
 (19): \quad
  \ph(\sm_1 \xi_1)
   \st{(2)}{\s}
  b_1 d_2 (d_1 c_2)
   \st{(28)}{\s}
  b_1 (d_2 c_0)
   \st{(28)}{\s}
  b_1 c_1
   \st{(25)}{\s}
  (d_2 d_0) c_1
   \st{(28)}{\s}
  a_2 b_2
   \st{(24_1)}{=}
  \ph(\xi_1), \\ \quad

  \ph(\eta_1 \sm_1)
   \st{(2)}{\s}
  a_2 (b_2 b_1) d_2 d_1 b_2
   \st{(26)}{\s}
  (a_2 d_0) d_2 d_1 b_2
   \st{(3)}{\s}
  (a_1 d_2) d_1 b_2
   \st{(3)}{\s}
  (a_0 d_1) b_2
   \st{(3),(24_1)}{\s}
  \ph(\eta_1);
 \end{array}$$

$$(20): \quad
  \ph(\sm_1 \sm_1^{-1})
   \st{(24_1)}{=}
  (b_1 d_2 d_1 b_2) (d_2 b_1 b_2 d_1)
   \st{(2)}{\s}
  1
   \st{(2)}{\s}
  (d_2 b_1 b_2 d_1) (b_1 d_2 d_1 b_2)
   \st{(24_1)}{=}
  \ph(\sm_1^{-1} \sm_1); $$


$$\begin{array}{l}
 (21): \quad
  \ph(\sm_{2} \sm_1 \sm_{2}) 
   \st{(2)}{\s}
  d_2 b_1 d_2 d_1 b_2^2 b_1 d_2^2 d_1 b_2^{2} 
   \st{(26)}{\s} 
  d_2 (b_1 d_2 d_1 b_2) d_0 d_2^2 d_1 b_2^{2} 
   \\ \qquad \st{(39)}{\s} 
  d_2 d_0 (b_1 d_2 d_1 b_2) d_2^2 d_1 b_2^{2} 
   \st{(25),(2)}{\s} 
  b_1^2 d_2 (d_1 d_2) d_1 b_2^{2} 
   \st{(26)}{\s} 
  b_1^2 d_2 (d_1 d_2) d_1 (d_0 d_1) b_2 
   \\ \qquad \st{(25)}{\s}   
  b_1^2 d_2 b_0 d_1 (d_0 d_1) b_2 
   \st{(2)}{\s} 
  b_1^2 (d_2 b_0 d_1 d_0 b_2) d_2 d_1 b_2 
   \st{(45)}{\s} 
  b_1^2 (b_0 d_2 d_1 b_2 d_0) d_2 d_1 b_2 
   \\ \qquad \st{(25)}{\s} 
  b_1^2 (d_1 d_2) d_2 d_1 b_2 d_0 d_2 d_1 b_2 
   \st{(26),(2)}{\s} 
  b_1 d_2^2 d_1 b_2 (b_2 b_1) d_2 d_1 b_2 
   \st{(2)}{\s}
  \ph(\sm_1 \sm_{2} \sm_1);
\end{array}$$


$$\begin{array}{l}
 (22): \quad
 \ph(\la_{2p-1,2} \Sm_{1,p})
   \st{(47)}{\s}
  (d_2 x_{2p-1,2} b_2^{2}) (b_1 d_2^p d_1 b_2^{p})
   \st{(26)}{\s}
  d_2 (x_{2p-1,2} b_2) d_0 d_2^p d_1 b_2^{p}
   \\ \qquad \st{(40)}{\s} 
  d_2 d_0 (x_{2p-1,2} b_2) d_2^p d_1 b_2^{p}
   \st{(2)}{\s}
  (d_2 d_0) (d_2^{p-1} b_2^{p-1}) (x_{2p-1,2} d_2^{p-1}) d_1 b_2^{p}
   \\ \qquad \st{(25)}{\s} 
  b_1 d_2^{p-1} (b_2^{p-1} x_{2p-1,2} d_2^{p-1}) d_1 b_2^{p}
   \st{(42)}{\s} 
  b_1 d_2^{p-1} d_1 (b_2^{p-1} x_{2p-1,2} d_2^{p-1}) b_2^{p}
   \\ \qquad \st{(2)}{\s}
  (b_1 d_2^{p-1} d_1 b_2^{p-1}) (x_{2p-1,2} b_2)
   \st{(47)}{\s}
  \ph(\Sm_{1,p-1} \la_{2p-1,1}), 
\end{array}$$

$$\begin{array}{l}
 (22): \quad
  \ph(\la_{2p-1,1} \bar\Sm_{1,p})
   \st{(48)}{\s}
  (x_{2p-1,2} b_2) (b_1^p d_2 d_1^p b_2)
   \st{(2)}{\s} 
  b_1^{p-1} (d_1^{p-1} x_{2p-1,2} b_2 b_1^p) d_2 d_1^p b_2
   \st{(43)}{\s} \\ \qquad
  b_1^{p-1} d_2 (d_1^{p-1} x_{2p-1,2} b_2 b_1^p) d_1^p b_2
   \st{(2)}{\s} 
  (b_1^{p-1} d_2 d_1^{p-1} b_2) (d_2 x_{2p-1,2} b_2^{2})
   \st{(48_1),(24)}{\s}
  \ph(\bar\Sm_{1,p-1} \la_{2p-1,2});
\end{array}$$

$$ \begin{array}{l}
 (22): \quad
 \ph(\la_{2q,2} \Sm_{1,q})
   \st{(47)}{\s}
  (d_2^{2} x_{2q,2} b_2^{2}) (b_1 d_2^q d_1 b_2^{q})
   \st{(2)}{\s}  
  d_2^{2} x_{2q,2} b_2^2 b_1 d_2^q d_1 b_2^{q}
   \\ \qquad \st{(26)}{\s}
  d_2 (d_2 x_{2q,2} b_2) d_0 d_2^q d_1 b_2^{q}
   \st{(41)}{\s} 
  d_2 d_0 (d_2 x_{2q,2} b_2) d_2^q d_1 b_2^{q}
   \st{(2)}{\s} 
  (d_2 d_0) (d_2 x_{2q,2} d_2^{q-1}) d_1 b_2^{q}
   \\ \qquad \st{(25)}{\s}
  b_1 (d_2 x_{2q,2} d_2^{q-1}) d_1 b_2^{q}
   \st{(2)}{\s}  
  b_1 d_2^q (b_2^{q-1} x_{2q,2} d_2^{q-1}) d_1 b_2^{q}
   \st{(44)}{\s}
  b_1 d_2^q d_1 (b_2^{q-1} x_{2q,2} d_2^{q-1}) b_2^{q}
   \\ \qquad \st{(2)}{\s}
  (b_1 d_2^q d_1 b_2^{q}) (d_2 x_{2q,2} b_2)
   \st{(47),(24_1)}{\s}
  \ph(\Sm_{1,q} \la_{2q,1}), 
\end{array}$$

$$\begin{array}{l}
 (22): \quad
  \ph(\la_{2q,1} \bar\Sm_{1,q})
   \st{(48_1)}{\s}
  (d_2 x_{2q,2} b_2) (b_1^q d_2 d_1^q b_2)
   \st{(2)}{\s}  
  d_2 x_{2q,2} b_2 b_1^q d_2 d_1^q b_2
   \\ \qquad \st{(6)}{\s}
  d_2 (d_0 b_1^{q-2} x_{2q,1} d_1^{q-2} b_0) b_2 b_1^q d_2 d_1^q b_2
   \st{(26),(2)}{\s}
  b_1^{p} (d_1 x_{2q,1} b_1) d_2 d_1^q b_2
   \\ \qquad  \st{(41)}{\s}
  b_1^{q} d_2 (d_1 x_{2q,1} b_1)  d_1^q b_2
   \st{(2)}{\s} 
  b_1^{q} d_2 d_1 (d_1^{q-2} b_0 x_{2q,2} d_0 b_1^{q-2}) d_1^{q-1} b_2
   \\ \qquad  \st{(25),(2)}{\s}
  (b_1^{q} d_2 d_1^{q} b_2) (d_2^{2} x_{2q,2} b_2^{2})
   \st{(48_1),(24)}{\s}
  \ph(\bar\Sm_{1,q} \la_{2q,2});
\end{array}$$


$$\begin{array}{l}
 (23): \quad
  \ph(\la_{2p-1,1} \Sm'_{1,p-1})
   \st{(49)}{\s}
  (x_{2p-1,2} b_2) D_{p,2}
   \st{(9)}{\s}
  D_{p-1,2} (x_{2p-1,2} b_2)
   \st{(49)}{\s}
  \ph(\Sm'_{1,p-2} \la_{2p-1,1}), \\ \qquad

  \ph(\la_{2q,1} \Sm'_{1,q-1})
   \st{(49)}{\s}
  (d_2 x_{2q,2} b_2) D_{q,2} 
   \st{(10)}{\s}
  D_{q,2} (d_2 x_{2q,2} b_2)
   \st{(49)}{\s}
  \ph(\Sm'_{1,q-1} \la_{2q,1}). \; \qed
 \end{array}$$


\setcounter{section}{6}
\section{Further approaches to classification of spatial graphs}

Theorems~1--2 reduce the isotopy classification of
 spatial graphs to a word problem in the semigroups
 $RSG_n,NSG_n$.
A solution of the latter problem will provide an algorithmic
 classification of spatial graphs up to ambient isotopy in $\R^3$.
\smallskip

{\bf Problem 1.}
Find an algorithm to decide, whether two central elements of
 the semigroup $RSG_n$ (respectively, $NSG_n$) are equal.
\smallskip

In Subsection~7.1, we shall study the semigroups $RSG_n,NSG_n$
 via representation theory of groups.
In Proposition~3 from Subsection~7.2 we calculate
 a presentation for the fundamental group $\pi_1(S^3-G)$
 of a spatial graph $G$, in terms of three-page embeddings.
We shall need this description to get a lower bound for
 the three-page complexity $tp(G)$.
In Subsection~7.3, we develop a complexity theory for spatial
 graphs by introducing the three-page complexity $tp(G)$ in
 Definition~23.
We obtain a lower bound of $tp(G)$ in terms of the group
 $\pi_1(S^3-G)$, see Proposition~6 in Subsection~7.4.
In Subsection~7.5 we describe all spatial graphs up to complexity~6.


\subsection{Groups associated with the semigroups
 $RSG_n,NSG_n$.}

By Theorems 1--2 any representation of the semigroups
 $RSG_n,NSG_n$, that is non-trivial on their centers, gives
 an isotopy invariant of spatial graphs.
To study representations of groups is more simpler than that
 for semigroups.
So, one wants to map the semigroups $RSG_n,NSG_n$ to some groups, and then
 to study representations of these groups.
\smallskip

{\bf Definition 20} (\emph{the group $\wt{F}$ associated with
 a semigroup $F$}).
Let $F$ be any finitely presented semigroup generated by
 a set $A$ and relations $R$.
Put $A^{-1}=\{a^{-1} \vl a\in A\}$.
Then the \emph{associated group $\wt{F}$} is generated by
 the set $A\cup A^{-1}$ and the relations $R$.
The \emph{natural homomorphism} $F\to \wt{F}$ is given by
 the embedding of generators $A\to A\cup A^{-1}$.
\ed
\smallskip

\begin{lemma}[the associated groups $\wt{RSG_n},\wt{NSG_n}$]
The groups $\wt{RSG_n},\wt{NSG_n}$
 are isomorphic to the free abelian group of rank $n+1$,
 generated by
 $\ti a_i,\ti x_m$, where $i\in\Z_3$, $3\leq m\leq n$.
The natural homomorphisms $RSG_n\to\wt{RSG_n}$
 and $NSG_n\to\wt{NSG_n}$ are given by\footnote{
After the proof of Lemma~6 we shall write elements of
 $\wt{RSG_n},\wt{NSG_n}$ in the abelian form.}
$$\left\{ \begin{array}{l}
 a_i\mapsto \ti a_i, \quad
 b_i\mapsto \ti a_{i-1} - \ti a_{i+1}, \quad
 c_i\mapsto - \ti a_i, \quad
 d_i\mapsto \ti a_{i+1} - \ti a_{i-1}, \quad
 x_{2q,i}\mapsto \ti x_{2q}, \\

 x_{2p-1,0}\mapsto \ti x_{2p-1}, \quad
 x_{2p-1,1}\mapsto \ti x_{2p-1} + \ti a_{2} - \ti a_{0}, \quad
 x_{2p-1,2}\mapsto \ti x_{2p-1} + \ti a_{2} - \ti a_{1}.
\end{array} \right.$$
\end{lemma}
\begin{proof}
Let $\ti a_i, \ti b_i, \ti c_i, \ti d_i$ be the images of $a_i,b_i,c_i,d_i$
 under $RSG_n\to\wt{RSG_n}$ and $NSG_n\to\wt{NSG_n}$.
For the above new symbols, the relations (1)--(10) convert
 to\footnote{
 In the proof of Lemma~6, we write the relations
  $\wt{(1)}-\wt{(10)}$ in the multiplicative form.}
 $\wt{(1)}-\wt{(10)}$.

Firstly, let us consider the case $n=2$, i.e.
 the Dynnikov semigroup $DS=RSG_2=NSG_2$.
The elements $\ti b_i, \ti d_i\in \wt{DS}$ are invertible by~(2).
In order to get a presentation of $\wt{DS}$,
 let us add the symbols $\ti a_i^{-1},\ti c_i^{-1}$
 that are inverses to $\ti a_i, \ti c_i$, respectively.
We have
 $\wt{(3)}\; \ti b_i = \ti a_{i-1} \ti c_{i+1},
        \ti d_i = \ti a_{i+1} \ti c_{i-1}$.
Then $\ti b_i = \ti a_{i-1} \ti c_{i+1}=\ti c_{i-1}^{-1} \ti a_{i+1}^{-1}$
 or $\ti c_{i-1} \ti a_{i-1} \ti c_{i+1} \ti a_{i+1}=1$.
Put $\ti e_i=\ti c_i \ti a_i$, $i\in\Z_3$.
Hence $\ti e_2 \ti e_1=\ti e_0 \ti e_2=\ti e_1 \ti e_0=1$,
 i.e. $\ti e_0=\ti e_1=\ti e_2$ and $\ti e_0^2=1$.
The relation (1) converts to
 $\wt{(1)}\; \ti a_1 \ti c_2 \ti a_2 \ti c_0 \ti a_0 \ti c_1=1$
 or $\ti c_1 \ti a_1 \ti c_2 \ti a_2 \ti c_0 \ti a_0=1$, i.e.
 $\ti e_0^3=1$ and $\ti e_0=\ti e_1=\ti e_2=1$.
So, the elements $\ti a_i,\ti c_i$ are the inverse ones.
Then
 $\wt{(3)}\; \ti a_{i+1} = \ti a_{i-1} \ti d_{i}$ and
       $\ti b_i = \ti a_{i+1} \ti c_{i+1}$
 imply $\ti a_{i+1}=\ti a_{i-1} \ti b_i^{-1}$ and
 $\ti b_i=\ti a_{i-1} \ti a_{i+1}^{-1}$, respectively.
Hence
 $\ti a_{i+1}=\ti a_{i-1} \ti a_{i+1} \ti a_{i-1}^{-1}$, i.e.
 all elements $\ti a_i$ commute with each other.
So, $\wt{DS}$ is the free abelian group of rank~3,
 generated by $\ti a_0,\ti a_1,\ti a_2$.
The other letters are
 $\ti b_i=\ti a_{i-1} \ti a_{i+1}^{-1}$,
 $\ti c_i=\ti a_i^{-1}$, $\ti d_i=\ti a_{i+1} \ti a_{i-1}^{-1}$.

In the general case $n>2$, it suffices to check that
 the images $\ti x_{m,i}$ of the letters $x_{m,i}$
 commute with each other and with all $\ti a_i$.
The relations
 $$\wt{(5)} \; \ti x_{2p-1,i} \ti d_i^{p-1}=
 \ti a_i (\ti x_{2p-1,i} \ti d_i^{p-1}) \ti a_i^{-1},\quad
   \wt{(6)} \; \ti d_i \ti x_{2q,i} \ti d_i^{q-1}=
 \ti a_i (\ti d_i \ti x_{2q,i} \ti d_i^{q-1}) \ti a_i^{-1}$$
 imply $\ti x_{k,i}\lra \ti a_{i}$, where $3\leq k\leq n$.
By~$\wt{(4)}$ we get $\ti a_i\lra \ti x_{k,i\pm 1}$.
By~$\wt{(8)}$ $\ti x_{2p-1,i} \ti b_i\lra \ti \ti x_{m,i+1}$ and
       $\ti d_i \ti x_{2q,i} \ti b_i\lra \ti x_{m,i+1}$
 we have $\ti x_{k,i}\lra \ti x_{m,i+1}$.
By~$\wt{(4)}$ we conclude $\ti x_{m,i+1}\lra \ti x_{k,i\pm 1}$.
The symbols $\ti x_{m,1}, \ti x_{m,2}$ are expressed
 in terms of $\ti x_m=\ti x_{m,0}$ by $\wt{(4)}$, $\wt{(31)}$,
 $\wt{(33)}$.
\end{proof}

For each $i\in\Z_3$, let us introduce the linear functions
 $F_i:\wt{RSG_n},\wt{NSG_n}\to\Z$ by $F_i(\ti a_i)=0$,
 $F(\ti a_{i\pm 1})=1$, $F_i(\ti x_m)=0$ except $i=2$ and $m=2p-1$.
In the latter case, put $F_2(\ti x_{2p-1})=1$.
Denote by $|\be|$ the difference between the number of the left and
 right brackets in a bracket expression $\be$.
The following claim is an easy observation, see Fig.~2.

\begin{claim}
For any word $w\in RSG_n,NSG_n$ and each $i\in\Z_3$,
 we have $F_i(\ti w)=|\be_i(w)|$, see Definition~10 in
 Subsection~3.4.
\qed
\end{claim}

We are ready to describe the images of $RSG_n$ and $NSG_n$ in
 the associated groups.

\begin{proposition}[the centers of $RSG_n,NSG_n$]
Under the natural homomorphisms\\ $RSG_n\to \wt{RSG_n}$ and
 $NSG_n\to \wt{NSG_n}$, the centers of $RSG_n,NSG_n$
 map to the set
 $$\{ \; -z\ti a_0 -z\ti a_1 +z\ti a_2 +k_3\ti x_3 +\cdots +k_n\ti x_n \vl
      \sum_{ 2\leq p\leq \frac{n+1}{2} } k_{2p-1} = 2z, \;
      k_m\geq 0, \; 3\leq m\leq n \; \}.$$
The center of $DS=RSG_2=NSG_2$
 maps to the zero element $0\in \wt{DS}$.
\end{proposition}
\begin{proof}
Let $\ti w$ be the image of either a word $w\in RSG_n$ or
 $w\in NSG_n$ under the natural homomorphism.
By Lemma~6 the word $\ti w\in \wt{RSG_n},\wt{NSG_n}$ has the form
 $\ti w=x\ti a_0 +y\ti a_1 +z\ti a_2+\sum_{m=3}^n k_m\ti x_m$,
 where $x,y,z,k_m\in\Z$.
Since the words $\be_i (w)$ are balanced, then
 $|\be_i(w)|=0$ for each $i\in\Z_3$.
By Claim~8 we have $F_0(\ti w)=y+z=0$,
 $F_1(\ti w)=z+x=0$,
 $F_2(\ti w)=x+y+\sum_p k_{2p-1}=0$, i.e.
 $x=y=-z$ and $\sum_p k_{2p-1}=2z$.
Conversely, any word of this type is the image of the central element
 $w = a_2^z x_{3,0}^{k_3} \cdots x_{n,0}^{k_n} c_0^{z} c_1^{z}\in
      RSG_n,NSG_n.$
\end{proof}
\smallskip

The images of the centers of $RSG_3\cong NSG_3$ are
 $\{ \; -z\ti a_0 -z\ti a_1 +z\ti a_2 +2z\ti x_3 \vl z\geq 0 \; \}$.
For singular knots, the center of $RSG_{ \{4\} }$ maps
 to the subset $\{ k\ti x_4 \vl k\geq 0\}$ of
 the free abelian group of rank 4, generated by the letters $\ti a_i,\ti x_4$.
If the image of a word in $\wt{ RSG_{\{4\}} }$ is
 $\wt{w_G}=k\ti x_4$, then the given singular knot $G$ contains exactly
 $k$ singular points.
More generally, by the image of an encoding word $w_G$ in the associated group
 we may reconstruct only the set of the vertices of $G$.
It means that the algebrac approach to the method of
 three-page embeddings could be effective only via semigroups,
 not via groups.


\subsection{Presentation of the fundamental group $\pi_1(S^3-G)$.}

Further we shall consider spatial graphs $G\subset S^3$ up
 to arbitrary homeomorphism $f:S^3\to S^3$, which can reverse
 the orientation of $S^3$.
We assume that $S^3$ is obtained from $\R^3$ by adding
 a base infinity point $\infty$.
For a knot (or a graph) $K\subset S^3$, the fundamental
 group $\pi_1(S^3-K)$ is said to be \emph{the knot group} $\pi(K)$.
Neuwirth constructed a presentation of the knot
 group by exploiting an arc presentation of a given knot \cite{Ne}.
\emph{The arc presentation} is an embedding of a knot
 into a book with finitely many pages so that
 each page contains exactly one arc.
We modify the Neuwirth construction for
 three-page embedings of spatial graphs.

Adding the infinity point $\infty$ to the axis $\al$ we get
 the circle $\bar\al\subset S^3$.
Let us demonstrate our computations for 
 the trefoil $K$ in Fig.~7.
We are going to choose \emph{the Neuwirth loops}
 lying in $S^3$ near the pages $P_i$,
 the base point is $\infty\in\bar\al$.
For each arc $\ga\subset K\cap P_i$, 
 let us take a loop 
 going around $\ga$ and all the arcs 
 lying in $P_i$ farther from $\al$ than $\ga$.
See the right picture of Fig.~7.
For example, for the arc $A_1A_2\subset P_1$ in Fig.~7,
 we take the loop $r$ near $P_1$.
Similarly, for the arc $A_1A_3\subset P_0$, we take
 the loop $u_0$ near $P_0$.
All the choosen loops will be generators of $\pi(K)$.
Then with every segment $A_j A_{j+1}\subset \al$ one can associate
 two or three loops (at most one near each page) going around nearest arcs.
To get defining \emph{Neuwirth's relations} let us write
 the associated loops from $P_0$ to $P_2$.
For instance, the segment $A_1 A_2\subset \al$ 
 provides the relation $u_0 r=1$.
The segment $A_2 A_3$ gives $u_0 v_1=1$.


\begin{picture}(450,130)(20,0)

\put(150,80){\vector(1,0){300}}
\put(445,85){$\alpha$}
\put(130,80){$\Y$}
\put(150,110){$P_0$}
\put(430,25){$P_1$}
\put(145,35){$P_2$}
\put(90,30){$K$}
\put(140,120){\line(1,0){285}}
\put(165,20){\line(1,0){285}}
\put(140,120){\line(1,-4){25}}
\put(425,120){\line(1,-4){25}}
\put(150,80){\line(-1,-4){13}}
\put(435,80){\line(-1,-4){5}}
\put(422,30){\line(1,4){5}}
\put(137,30){\line(1,0){20}}
\put(165,30){\line(1,0){10}}
\put(180,30){\line(1,0){10}}
\put(195,30){\line(1,0){10}}
\put(210,30){\line(1,0){10}}
\put(225,30){\line(1,0){10}}
\put(240,30){\line(1,0){10}}
\put(255,30){\line(1,0){10}}
\put(270,30){\line(1,0){10}}
\put(285,30){\line(1,0){10}}
\put(300,30){\line(1,0){10}}
\put(315,30){\line(1,0){10}}
\put(330,30){\line(1,0){10}}
\put(345,30){\line(1,0){10}}
\put(360,30){\line(1,0){10}}
\put(375,30){\line(1,0){10}}
\put(390,30){\line(1,0){10}}
\put(410,30){\line(1,0){12}}

{\thicklines
\put(20,25){\line(0,1){35}}
\put(20,25){\line(1,0){60}}
\put(20,60){\line(1,0){30}}
\put(60,60){\line(1,0){30}}
\put(30,35){\line(0,1){20}}
\put(30,35){\line(1,0){40}}
\put(40,45){\line(0,1){10}}
\put(40,45){\line(1,0){20}}
\put(30,65){\line(0,1){30}}
\put(30,95){\line(1,0){60}}
\put(40,65){\line(0,1){20}}
\put(40,85){\line(1,0){40}}
\put(50,60){\line(0,1){15}}
\put(50,75){\line(1,0){20}}
\put(60,45){\line(0,1){15}}
\put(70,35){\line(0,1){20}}
\put(80,25){\line(0,1){30}}
\put(70,65){\line(0,1){10}}
\put(80,65){\line(0,1){20}}
\put(90,60){\line(0,1){35}}

\put(165,80){\circle*{3}}
\put(155,85){\footnotesize $A_1$}
\put(200,80){\circle*{3}}
\put(200,85){\footnotesize $A_2$}
\put(220,80){\circle*{3}}
\put(220,85){\footnotesize $A_3$}
\put(260,80){\circle*{3}}
\put(245,85){\footnotesize $A_4$}
\put(280,80){\circle*{3}}
\put(273,70){\footnotesize $A_5$}
\put(310,80){\circle*{3}}
\put(350,80){\circle*{3}}
\put(370,80){\circle*{3}}
\put(390,80){\circle*{3}}
\put(387,85){\footnotesize $A_{m-1}$}
\put(420,80){\circle*{3}}
\put(418,85){\footnotesize $A_{m}$}

\put(165,80){\line(2,3){10}}
\put(165,80){\line(1,-2){10}}
\put(175,60){\line(1,0){5}}
\put(175,95){\line(1,0){15}}
\put(185,60){\line(1,0){5}}
\put(195,95){\line(1,0){15}}
\put(200,80){\line(-1,-2){10}}
\put(200,80){\line(1,-2){5}}
\put(210,95){\line(2,-3){20}}
\put(210,65){\line(1,-1){10}}
\put(225,50){\line(2,-1){10}}
\put(230,65){\line(1,0){5}}
\put(240,42){\line(4,-1){10}}
\put(240,65){\line(1,0){55}}
\put(255,40){\line(1,0){5}}
\put(260,80){\line(1,1){15}}
\put(260,80){\line(1,-1){10}}
\put(265,40){\line(1,0){5}}
\put(275,40){\line(1,0){5}}
\put(280,60){\line(2,-1){10}}
\put(280,80){\line(1,1){20}}
\put(280,80){\line(3,-2){15}}
\put(280,100){\line(1,1){10}}
\put(285,40){\line(1,0){5}}
\put(290,110){\line(1,0){20}}
\put(295,40){\line(1,0){5}}
\put(295,65){\line(1,1){25}}
\put(300,100){\line(1,0){10}}
\put(303,40){\line(1,0){7}}
\put(305,50){\line(-4,1){10}}
\put(305,50){\line(1,0){5}}
\put(315,40){\line(1,0){5}}
\put(315,50){\line(1,0){5}}
\put(315,60){\line(1,0){5}}
\put(315,60){\line(-2,1){10}}
\put(315,100){\line(1,0){25}}
\put(315,110){\line(1,0){25}}
\put(320,90){\line(1,0){20}}
\put(325,40){\line(1,0){5}}
\put(325,50){\line(1,0){5}}
\put(330,60){\line(-1,0){5}}
\put(330,60){\line(2,1){10}}
\put(333,40){\line(1,0){7}}
\put(335,50){\line(1,0){5}}
\put(345,40){\line(1,0){5}}
\put(345,90){\line(1,0){15}}
\put(345,100){\line(1,0){25}}
\put(345,110){\line(1,0){45}}
\put(350,80){\line(-1,-2){5}}
\put(350,80){\line(2,-3){10}}
\put(355,50){\line(-1,0){10}}
\put(355,50){\line(1,2){5}}
\put(360,65){\line(1,0){30}}
\put(360,40){\line(-1,0){5}}
\put(360,40){\line(2,1){10}}
\put(370,80){\line(-1,1){10}}
\put(370,80){\line(-1,-2){5}}
\put(380,60){\line(-1,-2){5}}
\put(390,80){\line(-1,1){20}}
\put(390,80){\line(-1,-2){5}}
\put(395,65){\line(1,0){10}}
\put(420,80){\line(-1,1){30}}
\put(420,80){\line(-1,-1){15}}
}

\put(175,85){\line(1,0){17}}
\put(175,85){\line(0,1){5}}
\put(192,85){\vector(0,1){30}}
\put(175,100){\line(0,1){10}}
\put(195,105){$u_0$}

\put(182,75){\line(1,0){10}}
\put(182,75){\line(0,-1){30}}
\put(192,75){\line(0,-1){5}}
\put(192,60){\vector(0,-1){35}}
\put(180,35){$r$}

\put(207,67){\line(1,0){10}}
\put(207,67){\line(0,-1){7}}
\put(207,55){\line(0,-1){5}}
\put(217,67){\line(0,-1){5}}
\put(217,55){\line(0,-1){5}}
\put(217,45){\vector(0,-1){10}}
\put(205,40){$v_1$}

\put(237,75){\line(1,0){15}}
\put(237,75){\line(0,-1){50}}
\put(252,75){\line(0,-1){5}}
\put(252,60){\vector(0,-1){35}}
\put(255,50){$s$}

\put(293,60){\line(1,0){12}}
\put(293,60){\line(0,-1){6}}
\put(293,50){\line(0,-1){5}}
\put(293,42){\line(0,-1){7}}
\put(305,60){\line(0,-1){5}}
\put(305,47){\line(0,-1){5}}
\put(305,38){\vector(0,-1){5}}
\put(280,45){$v_2$}

\put(323,75){\line(1,0){12}}
\put(323,75){\line(0,-1){5}}
\put(323,65){\line(0,-1){8}}
\put(323,53){\line(0,-1){6}}
\put(323,43){\line(0,-1){6}}
\put(335,75){\line(0,-1){10}}
\put(335,60){\line(0,-1){5}}
\put(335,48){\line(0,-1){5}}
\put(335,38){\vector(0,-1){5}}
\put(340,55){$v_3$}

\put(393,75){\line(1,0){12}}
\put(393,75){\line(0,-1){30}}
\put(405,75){\line(0,-1){5}}
\put(405,60){\vector(0,-1){35}}
\put(385,35){$t^{-1}$}

\put(270,83){\line(1,0){8}}
\put(270,83){\line(0,1){5}}
\put(270,93){\line(0,1){8}}
\put(278,83){\vector(0,1){30}}
\put(265,105){$u_1$}

\put(300,95){\line(1,0){12}}
\put(300,102){\line(0,1){5}}
\put(300,112){\line(0,1){5}}
\put(312,95){\vector(0,1){20}}
\put(295,85){$u_2$}

\put(325,83){\line(1,0){17}}
\put(325,83){\line(0,1){5}}
\put(325,93){\line(0,1){5}}
\put(325,103){\line(0,1){5}}
\put(325,113){\line(0,1){5}}
\put(342,83){\vector(0,1){30}}
\put(330,113){$u_3$}

\put(20,0){{\bf Fig.~7.} 
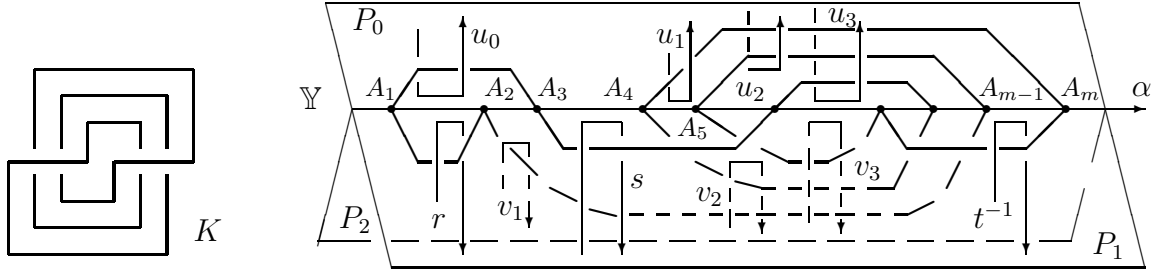
 The trefoil $K$ is encoded by the word
 $w_K=a_2 d_0 d_2 a_1^2 b_2 b_0 c_1^2 c_2$.}

\end{picture}
\vspace{0mm}

In fact, we get the following presentation:
 $$\begin{array}{l}
 \pi(K)=\ab{u_0,u_1,u_2,u_3,r,s,t,v_1,v_2,v_3 \vl \\
 u_0 r = u_0 v_1 = s v_1 = u_1 s v_2 = u_2 s v_3 = u_3 v_3 =
 u_3 t^{-1} v_2 = u_2 t^{-1} v_1 = u_1 t^{-1} = 1}.
 \end{array}$$
We have $r^{-1}=v_1^{-1}=u_0=u_2 t^{-1}=s$,
 $v_2^{-1}=u_1 s=u_3 t^{-1}$, $v_3^{-1}=u_2 s=u_3$, and $u_1=t$.
Then $u_0=s$, $u_1=t$, $u_2=st$, $u_3=tst$, and the relation $u_2s=u_3$
 converts to $sts=tst$.
So, we have the standard presentation
 $\pi(K)=\ab{s,t \vl sts=tst}$ of the trefoil group.

\begin{proposition}
The above modification of the Neuwirth construction provides a presentation
 of the graph group $\pi(G)=\pi_1(S^3-G)$ for any spatial graph $G\subset \Y$.
\end{proposition}
\begin{proof}
For each segment $A_j A_{j+1}\subset \al$, 
 let us choose a small subsegment $I_j\subset A_j A_{j+1}$.
For sufficiently small $\e>0$, 
 let $N(I_j)$ be the $\e$-neighbourhood (a cylinder) of $I_j$, and
 $N(A_1A_m)$ be the $\frac{\e}{2}$-neighbourhood (also a cylinder)
 of the segment $A_1A_m$.
Let us join the two adjacent cylinders $N(I_j),N(I_{j+1})$ by
 an arc $T_j\subset \R^3-\Y$.
Put $X=N(I_1)\cup T_1\cup N(I_2)\cup\ldots\cup T_{m-2}\cup N(I_{m-1})$ and
    $Y=S^3-(G\cup N(A_1 A_m))$.
Then the space $X$ is contractible, i.e. $\pi_1(X)=1$.
The group $\pi_1(Y)$ is freely generated by the Neuwirth loops.
Moreover, the space $X\cup Y$ is homotopically equivalent to $S^3-G$, 
 i.e. $\pi_1(X\cup Y)=\pi(G)$.

By well-known Seifert-Van-Kampen's Theorem \cite[chapter 4]{Ma},
 in order to get a presentation of $\pi_1(X\cup Y)$
 we should add the relations
 corresponded to all generators of $\pi_1(X\cap Y)$.
The intersection $X\cap Y$ consists of 
 the tubes $N(I_j)-N(A_1 A_m)$ joined by the arcs $T_j$.
Hence the group $\pi_1(X\cap Y)$ is generated by
 loops going around the segments
 $I_j\subset A_j A_{j+1}$.
So, all defining relations of $\pi(G)$ are
 the Neuwirth relations.
\end{proof}
\smallskip

{\bf Definition 21} (\emph{the disjoint union, a vertex sum,
 an edge sum, a loop sum of graphs}).

(a) A spatial graph $F\subset S^3$ is called
 \emph{the disjoint union} of spatial graphs $G,H\subset S^3$
 and denoted by $G\sqcup H$, if there is a
 two-sided\footnote{
  It means that we may speak about points inside the sphere $S$
   and about points outside $S$.}
  2-sphere $S\subset S^3$ such that
 $F=G\cup H$, the subgraph $G\subset F$ lies inside the sphere $S$,
 the subgraph $H\subset F$ lies outside $S$.
\smallskip

(b) A spatial graph $F\subset S^3$ is called
 \emph{a vertex sum} of spatial graphs $G,H\subset S^3$
 and denoted by $G*H$, if there is a
 two-sided 2-sphere $S\subset S^3$ such that
 $F=G\cup H$, $F\cap S=v$ is either a vertex or a point inside
  a loop of $G$ and $H$, the subgraph $G-v$ lies
  inside the sphere $S$, $H-v$ lies outside $S$, see Fig.~8.
\smallskip

(c) A spatial graph $F\subset S^3$ is called
 \emph{an edge sum} of spatial graphs $G,H\subset S^3$
 and denoted by $G\vee H$, if there is a
 two-sided 2-sphere $S\subset S^3$ and an edge $e\subset F$ such that
 $F-e=G\cup H$, $F\cap S=e\cap S=$ 1 point,
 $G$ lies inside the sphere $S$, $H$ lies outside $S$.
\smallskip

(d) A spatial graph $F\subset S^3$ is called
 \emph{a loop sum} of spatial graphs $G,H\subset S^3$
 and denoted by $G\circ H$, if there is a
 two-sided 2-sphere $S\subset S^3$ and an arc $I\subset S$ such that
 $F=(G\cup H)-I$, $G\cap H=I$ is inside a loop $e_G\subset G$ and inside
  a loop $e_H\subset H$, the subgraph $G-I\subset F$ lies
  inside the sphere $S$, $H-I$ lies outside $S$, see Fig.~8.
\ed
\medskip

\begin{picture}(420,160)

\thicklines
\qbezier(70,50)(110,70)(80,110)
\qbezier(80,110)(40,170)(30,120)

\qbezier(30,100)(35,5)(100,55)
\qbezier(100,55)(160,100)(90,120)

\qbezier(60,120)(20,120)(-5,70)
\qbezier(-5,70)(-25,25)(45,35)

\put(165,110){\oval(70,70)[b]}
\put(165,110){\oval(70,70)[tl]}
\put(165,135){\oval(70,20)[tr]}
\put(215,90){\oval(70,70)[t]}
\put(215,90){\oval(70,70)[br]}
\put(215,65){\oval(70,20)[bl]}
\put(129,91){\circle*{7}}
\put(135,95){\Large $v$}
\put(210,55){\circle*{7}}
\qbezier(129,91)(150,10)(210,55)
\put(50,70){\large $G$}
\put(160,50){\large $H$}
\put(100,30){\large $G*H$}

\qbezier(360,40)(420,70)(390,110)
\qbezier(390,110)(350,170)(340,120)

\qbezier(340,100)(345,5)(410,55)
\qbezier(410,55)(470,100)(400,120)

\qbezier(370,120)(330,120)(305,70)
\qbezier(305,70)(285,25)(360,40)

\put(361,41){\circle*{7}}
\put(361,41){\line(-2,1){50}}
\put(291,76){\line(2,-1){10}}
\put(291,76){\circle*{7}}
\put(291,96){\circle{40}}
\put(280,130){\large $G\vee H$}
\put(320,65){\Large $e$}
\put(285,90){\large $G$}
\put(360,70){\large $H$}
\end{picture}

\begin{picture}(450,90)

{\thicklines

\put(30,70){\oval(50,50)[b]}
\put(30,70){\oval(50,50)[tl]}
\put(30,85){\oval(50,20)[tr]}
\put(55,65){\circle*{7}}

\put(60,55){\oval(50,50)[t]}
\put(60,55){\oval(50,50)[br]}
\put(60,40){\oval(50,20)[bl]}
\put(35,65){\circle*{7}}
\put(35,65){\line(1,0){20}}

\put(100,60){\huge $\circ$}

\put(150,70){\oval(50,50)[b]}
\put(150,70){\oval(50,50)[tl]}
\put(150,85){\oval(50,20)[tr]}

\put(180,55){\oval(50,50)[t]}
\put(180,55){\oval(50,50)[br]}
\put(180,40){\oval(50,20)[bl]}

\put(220,60){\huge =}

\put(280,80){\oval(50,50)[b]}
\put(280,80){\oval(50,50)[tl]}
\put(280,95){\oval(50,20)[tr]}
\put(305,75){\circle*{7}}

\put(330,65){\oval(90,50)[t]}
\put(330,50){\oval(90,20)[br]}
\put(330,50){\oval(90,20)[bl]}
\put(285,75){\circle*{7}}
\put(285,75){\line(1,0){20}}

\put(380,45){\oval(70,30)[t]}
\put(380,45){\oval(70,30)[br]}
\put(380,35){\oval(70,10)[bl]} }

\put(330,45){\line(0,1){5}}
\put(330,55){\line(0,1){5}}
\put(330,65){\line(0,1){10}}
\put(330,80){\line(0,1){5}}
\put(315,50){\Large $I$}
\put(260,20){\large $G\circ H$}

\put(20,0){{\bf Fig.~8.}
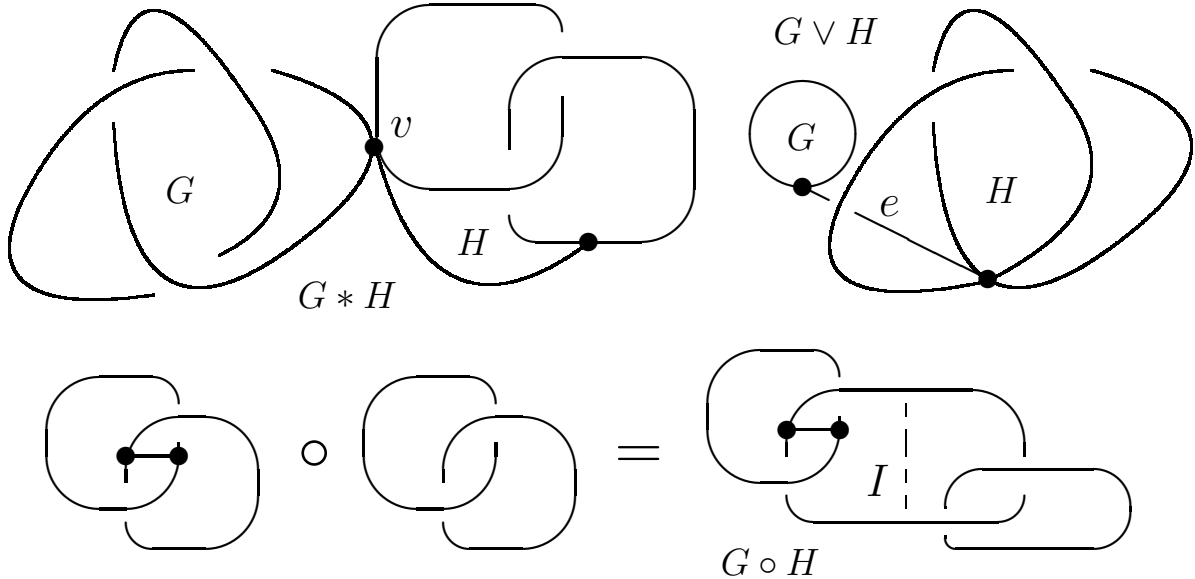
A vertex sum of spatial graphs, an edge sum, a loop sum.}
\end{picture}
\vspace{3mm}

The disjoint union $G\sqcup H$ is determined uniquely up to
   equivalence $f:S^3\to S^3$.
In the case of links, the notion of a loop sum coincides with
  the usual connected sum.
For arbitrary spatial graphs, a loop sum could be
 non-associative and non-commutative since
 the loops $e_G,e_H$ from Definition~21d
 can convert to edges distinct from loops in a loop sum $G\circ H$.
For finitely presented groups $\pi,\pi'$,
 the group $\pi*\pi'$ is called \emph{the free product} of $\pi,\pi'$.
A presentation of $\pi*\pi'$ could be obtained by uniting 
 generators and relations of $\pi,\pi'$ \cite[chapter 3]{Ma}.
Proposition~3 implies the following claim.
\smallskip

\begin{claim}
For any spatial graphs $G,H\subset S^3$, we have
 $\pi(G\sqcup H)\cong\pi(G*H)\cong\pi(G\vee H)\cong\pi(G)*\pi(H)$.
\qed
\end{claim}

Due to Proposition~3 we may also calculate the Alexander
 polynomial of oriented links embedded into $\Y$.
This allows us to classify an infinite family of singular knots
 with arbitrary numbers of singular points and crossing
 \cite[Propositions~2.4--2.5]{Ku2}.
Only singular knots with exactly one singular point and
 at most 6 crossings were classified in \cite{Ge}.


\subsection{Complexity theory for spatial graphs}

In our geometric approach, it is convenient to
 extend the notion of a three-page embedding to a more general one.
\smallskip

{\bf Definition 22} (\emph{general three-page embeddings}).
An embedding $G\subset\Y$ is a \emph{general three-page embedding}, if
 Conditions~(7.1), (7.2), (7.3) of Definition~7 hold and
\smallskip

(7.6) a neighbourhood of each $m$-vertex $A\in G$
       lies in the union of two pages\footnote{
  Possibly, a neighbourhood of an $m$-vertex $A\in G$ lies in
   exactly one page of $\Y$, even for $m=2$.} of $\Y$.
\ed
\medskip

A general three-page embedding of a spatial graph
 could be constructed analogously to
 Subsection~3.2, but we don't need to check Conditions~(4)--(6) of
 Subsection~3.2.
General three-page embeddings also could be encoded by finitely
 many letters.
Moreover, we may prove analogs of Theorems~1--2,
 but in this general case encoding semigroups
 should be generated by much more letters and defining relations.
\smallskip

{\bf Definition 23} (\emph{the arch number $ar(G)$,
 the three-page complexity $tp(G)$}).
(a) An \emph{arch} of a general three-page embedding $G\subset\Y$
 is a connected component of $G-\al$.
\emph{The arch number} $ar(G)$ is the number of the
 arches\footnote{
 In the case of a link $L\subset\Y$, $ar(L)$ equals to the number of the
  letters in the encoded word $w_L$.}
 in a given three-page embedding $G\subset\Y$.
(b) \emph{The three-page complexity $tp(G)$} is the minimal value
 of $ar(G)-2$ over all possible general three-page embeddings
 $G\subset\Y$ of a given spatial graph $G$.
\ed
\smallskip

It is not hard to see that $tp(G)\geq 0$ for all spatial graphs $G$.
Moreover, $tp(G)=0$ if and only if $G$ is the unknot $O_1$.
One can check that $tp(L)=4$ for the Hopf link $L$.
Let $S(p,q)$ be the non-oriented 2-bridge link having
 Shubert's normal form with the parameters $p,q\geq 1$
 \cite[chapter 12.A]{BZ}.
The link $S(p,q)$ could be encoded by the word
 $a_0 a_1^{p-1} b_2 b_1^{q-1} b_0 c_1^{p-q} d_1^{q-1} c_2 c_1^{q-1}$,
 i.e. $tp(S(p,q))\leq 2p+2q-2$.
\medskip

{\bf Conjecture 1.}
The three-page complexity of $S(p,q)$ is $2p+2q-2$ for $p+q\geq 3$.
\smallskip

We have checked Conjecture 1 for the Hopf link $S(2,1)$
 and for the trefoil $S(3,1)$, see Fig.~9 in Subsection~7.5.

\begin{lemma}
For any $k\in\N$, there is a finite number of spatial graphs $G$
 with $tp(G)=k$.
\end{lemma}
\begin{proof}
It suffices to estimate from above the number $TP_k$ of all three-page
 embeddings $G\subset\Y$ of spatial graphs, with $ar(G)=k$.
For such a three-page embedding, the number of the intersection
 points from $G\cap\al$ is not more than $k$.
For an $m$-vertex $A\in G\cap\al$, we may embed a neighbourhood of $A$
 into two pages of $\Y$ by not more than $4^m$ different monotone ways.
Hence we may estimate $TP_k$ very roughly as follows: $TP_k\leq (4^k)^k$.
\end{proof}

{\bf Problem 2.}
Find asymptotics for the number $N_n(k)$ of
 all prime spatial $n$-graphs with the three-page complexity $k$.
\medskip

To get the additivity (3.1) from Theorem~3 we need a geometric inversion.
\smallskip

{\bf Definition 26} (\emph{geometric inversion $f_{a,r}$}).
Let $S_{a,r}$ be a geometric 2-sphere with a center $a\in S^3$ and
 a radius $r>0$.
Then the geometric inversion $f_{a,r}:S^3\to S^3$ is defined
 as follows\footnote{
 Here $x,a\in S^3$ are usual 3-dimensional vectors,
  $|x-a|$ is the length of the vector $x-a$.}:
 $f_{a,r}(x)=a+\frac{r^2}{|x-a|^2}(x-a)$.
In particular, $f_{a,r}(a)=\infty$, $f_{a,r}(\infty)=a$,
 $f_{a,r}(x)=x$ for each $x\in S_{a,r}$.
Moreover, the inversion changes the orientation of $S^3$.
\ed
\smallskip

If the center $a$ of a geometric inversion $f_{a,r}$ lies in the axis $\al$
 of the book $\Y$, then $f_{a,r}(\al)=\al$ and $f_{a,r}(\Y)=\Y$.
Assume that the arches of a spatial graph $G\subset\Y$ are
 \emph{geometric}, i.e. they come to the axis $\al$ as
 perpendicular smooth curves.
Any three-page embedding of a graph is isotopic inside $\Y$ to
  such \emph{a geometric embedding}.

For a geometric three-page embedding $K\subset\Y$ of a singular knot $K$,
 two branches of $K$ are intersected non-transversally
 at a singular point $A\in K$, see the trivial graph $\te_4$ in Fig.~9a.
But we may always select branches at $A$ since any rigid isotopy
 keep a neighbourhood of $A$ in a (non-constant) plane.
Under inversion $f_{a,r}$ with $a\in\al$, a geometric arch
 with endpoints $b,c\in\al$ goes to a similar geometric arch with the
 endpoints $f_{a,r}(b),f_{a,r}(c)\in\al$, in the same page.

\begin{claim}
Let $e_G$ be a loop of a spatial graph $G$.
Let $G\subset\Y$ be a general three-page embedding.
Then there is a geometric inversion $f_{a,r}$
 such that $f_{a,r}(G)\subset\Y$, $ar(f_{a,r}(G))=ar(G)$,
 and the left extreme point from $f_{a,r}(G)\cap\al$
 belongs to the loop $f_{a,r}(e_G)\subset f_{a,r}(G)$.
\end{claim}
\begin{proof}
The loop $e_G$ has two extreme points from $e_G\cap\al$:
 $A_k,A_l$, where $k<l$.
Since the loop $e_G$ contains not more than one vertex of $G$,
 then one of these points (say $A_k$) is not a vertex of $G$.
Since the both geometric arches at $A_k$ come perpendicularly
 to the axis $\al$, we may take a geometric sphere $S_{a,r}$
 with $a\in\al$ and a small radius $r$, such that
 $A_k\in S_{a,r}$ and $G-A_k$ lies outside $S_{a,r}$.
Then $f_{a,r}(A_k)=A_k$, and $f_{a,r}(G-A_k)$ lies inside $S_{a,r}$.
Hence the point $A_k$ is now an extreme point
 from $f_{a,r}(G)\cap\al$, in the axis $\al$.
If $A_k$ is the right extreme point, then take a symmetric
 reflection of $f_{a,r}(G)\subset\Y$ in a plane $\perp\al$.
\end{proof}

\begin{proposition}
For any spatial graphs $G,H\subset S^3$, we have
\smallskip

(a) $tp(G\sqcup H)=tp(G)+tp(H)+2$;$\qquad$
(c) $tp(G\vee H)=tp(G)+tp(H)+3$;

(b) $tp(G*H)=tp(G)+tp(H)+2$;$\qquad$
(d) $tp(G\circ H)\leq tp(G)+tp(H)$.
\end{proposition}
\begin{proof}
(a) Let us take general three-page embeddings $G,H\subset\Y$
 with the minimal numbers of arches, i.e.
 $ar(G)=tp(G)+2$, $ar(H)=tp(H)+2$.
To get a general three-page embedding $G\sqcup H\subset\Y$,
 we attach two copies of $\Y$ along $\al$.
Then $tp(G\sqcup H)\leq ar(G\sqcup H)-2=ar(G)+ar(H)-2=tp(H)+tp(H)+2$.

Conversely, take a general three-page embedding $G\sqcup H\subset\Y$
 such that $ar(G\sqcup H)=tp(G\sqcup H)+2$.
The sphere $S$ from Definition~21a splits the embedding $G\sqcup H\subset\Y$
 into two disjoint parts in such a way that
 $w_{G\sqcup H}=u_1v_1\ldots u_kv_k$,
 where the words $u_1\ldots u_k$ and $v_1\ldots v_k$ encode
 the subgraphs $G,H\subset\Y$.
So, we have
 $ar(G)+ar(H)=ar(G\sqcup H)=tp(G\sqcup H)+2$.
Hence
 $tp(G)+tp(H)\leq ar(G)+ar(H)-4=tp(G\sqcup H)-2$ as required.
\smallskip

(b) Take general three-page embeddings $G,H\subset\Y$
 with $ar(G)=tp(G)+2$, $ar(H)=tp(H)+2$.
By Claim~10 choose another general three-page embeddings
 $G,H\subset\Y$ such that $ar(G)=tp(G)+2$, $ar(H)=tp(H)+2$, and also
 the right extreme point of $G\cap\al$
 (respectively, the left extreme point of $H\cap\al$)
 is the gluing point $v$ from Definition~21b.
Now we may attach the obtained embeddings $G,H\subset\Y$
 to get a general three-page embedding $G*H\subset\Y$
 with $ar(G*H)=ar(G)+ar(H)$ as in the item~(a).

Conversely, take a general three-page embedding $G*H\subset\Y$
 with $ar(G*H)=tp(G\sqcup H)+2$.
The sphere $S$ from Definition~21b splits the embedding $G*H\subset\Y$
 into two parts that are intersected at the point $v$.
These parts form two independent general three-page embeddings
 $G,H\subset\Y$ with $ar(G)+ar(H)=ar(G*H)$.
The proof finishes as in~(a).
\smallskip

The item (c) is analogous to (a) and (b).
Given general three-page embeddings $G,H\subset\Y$ by Claim~10
 we may construct a general three-page embedding $G\sqcup H\subset\Y$
 such that the right extreme point of $G\cap\al$
 (respectively, the left extreme point of $H\cap\al$)
 is an endpoint of the edge $e\subset G\vee H$ from
 Definition~21c.
Also we may assume that neighbourhoods of these endpoints
 lie in two common pages of $\Y$.
Otherwise it suffices to rotate the embedding of $G$ (say)
 to secure the above condition.
Now we are able to add the edge $e\subset G\vee H$ to the embedding
 $G\sqcup H\subset\Y$ and to get a general three-page embedding
 $G\vee H\subset\Y$ with $ar(G\vee H)=ar(G)+ar(H)+1$.
The proof finishes as in the item~(b).
\smallskip

The item (d) is similar to the first part of (c).
Take general three-page embeddings $G,H\subset\Y$
 with $ar(G)=tp(G)+2$, $ar(H)=tp(H)+2$.
By Claim~10 we may construct another general three-page embeddings
 $G,H\subset\Y$ that are intersected at the common
 "vertical" arc $I\perp\al$ from Definition~21d, see Fig.~8.
Here monotone Condition~(7.5) of Definition~7 does not play any role.
Now we may remove the arc $I$ form the union $G\cup H\subset\Y$
 and get a general three-page embedding $G\circ H\subset\Y$
 with $ar(G\circ H)=ar(G)+ar(H)-2=tp(G)+tp(H)+2$.
Then $tp(G\circ H)\leq ar(G\sqcup H)-2=ar(G)+ar(H)-4=tp(H)+tp(H)$.

The reverse of the item~(d) is much harder since the sphere
 $S$ from Definition~21d may intersect an embedding $G\circ H\subset\Y$
 in a terrible way.
\end{proof}

{\bf Conjecture 2.}
The three-page complexity is additive under a loop sum, i.e.\\
 $tp(G\circ H)=tp(G)+tp(H)$ for any spatial graphs $G,H\subset S^3$.
\medskip

If Conjecture 2 is true, then we shall get an hierarchy on the set
 of spatial graphs considered up to equivalence $f:S^3\to S^3$.
Proposition~4a implies that the three-page complexity of the trivial
 $k$-component link $O_k$ is
 $tp(O_k)=k\cdot tp(O_1)+2(k-1)=2k-2$.

Theorem~3 formulated in Subsection~1.3 follows from
 Lemma~7 and Proposition~4.


\subsection{Lower bound of the three-page complexity}

Here we find a lower bound for the three-page complexity $tp(G)$ in terms
 of the group $\pi(G)=\pi_1(S^3-G)$.
\smallskip

{\bf Definition 25} (\emph{the three-letters complexity $tl(\pi)$}).
(a) Let $\pi$ be any finitely presented group.
A presentation of $\pi$ is a \emph{three-letters presentation},
 if it contains $k$ generators and not more than $k-1$ relations,
 each relation consists of 3 generators or their inverses.
\smallskip

(b) \emph{The three-letters complexity} $tl(\pi)$ is the minimal number
 $k$ of generators for all three-letters presentations of $\pi$.
If $\pi$ has no three-letters presentations,
 then set $tl(\pi)=\infty$.
\ed
\smallskip

For instance, $\Z$ is a unique group of three-letters index 1.
All groups of $tl=2$ are $\Z *\Z$ and $\Z_3 *\Z$.
The group $\Z\oplus\Z\cong \ab{a,b,c\vl abc=acb=1}$ has $tl=3$.
All groups with $tl=3$ are listed in \cite[Example 2.11]{Ku2}.
The cyclic groups $\Z_k$ ($k>1$) have $tl=\infty$.
\smallskip

\begin{proposition}
(a) For any $k\in\N$, there is a finite number of groups $\pi$
 with $tl(\pi)=k$.

(b) If groups $\pi_1$ and $\pi_2$ have three-letters presentations,
 then $tl(\pi_1 *\pi_2)=tl(\pi_1)+tl(\pi_2)$.
\end{proposition}
\begin{proof}
(a) It suffices to estimate from above the number $TL_k$
 of all three-letters presentations of complexity $k$.
For such a presentation, there are not more than $3^k$ different
 relations, hence $TL_k\leq (3^k)^{k-1}$.
\smallskip

(b) Since the union of three-letters presentations for
 $\pi_1,\pi_2$ gives a three-letters presentation for $\pi_1 *\pi_2$, then
 $tl(\pi_1 *\pi_2)\leq tl(\pi_1)+tl(\pi_2)$.
If a relation from a three-letters presentation of $\pi_1 *\pi_2$
 contains two generators of $\pi_1$ (say) and a generator $g$
 of $\pi_2$, then $g\in\pi_1$ that is a contradiction.
So, we may split any three-letters presentation of $\pi_1 *\pi_2$
 into two three-letters presentation for $\pi_1$ and $\pi_2$.
Hence $tl(\pi_1 *\pi_2)\geq tl(\pi_1)+tl(\pi_2)$.
\end{proof}
\smallskip

{\bf Definition 26} (\emph{trivial graphs $\te_k$}).
\emph{The $\te_k$-graph} consists of 2 vertices joined by $k\geq 2$ edges.
\emph{The trivial graph $\te_k\subset S^3$} is the $\te_k$ graph
 embedded into $\R^2\subset S^3$.
For example, the trivial graph $\te_2$ is the unknot.
Each trivial graph $\te_k$ has a general three-page embedding
 $\te_k\subset\Y$ such that $ar(\te_k)=k$ and $\te_k\cap\al=2$ points.
\ed

\begin{proposition}
For any spatial graph $G$, distinct from a trivial graph $\te_k$,
 we have $tp(G)\geq tl(\pi(G))$.
For a trivial graph $\te_k$, we get $\pi(\te_k)=F_{k-1}$
 (the free group with $k-1$ generators), $tl(\pi(\te_k))=k-1$,
 $tp(\te_k)=k-2$.
\end{proposition}
\begin{proof}
Let us take a general three-page embedding $G\subset\Y$
 with the minimal number of arches, i.e. $ar(G)=tp(G)+2$.
Proposition~3 from Subsection~7.2 gives a presentation of $\pi(G)$
 with $ar(G)$ generators, all relations contain at most three letters.
Moreover, two Neuwirth's relations corresponding to
 the extreme segments $A_1A_2,A_{l-1}A_l$
 contain exactly two letters.
Hence at least two generators are superfluous, i.e.
 $tl(\pi(G))\leq ar(G)-2=tp(G)$.
The above reason does not hold, when there is exactly one extreme
 segment, i.e. $A_1A_2=A_{l-1}A_l$.
In this case $G$ is a trivial graph $\te_k$.
\end{proof}
\smallskip

{\bf Problem 3.}
Find lower bounds for the three-page complexity in terms of known
 polynomial invariants for links and spatial graphs.


\subsection{Spatial graphs up to complexity 6}

Fig.~9 shows all non-oriented links, spatial 3-graphs, and singular
 knots with three-page complexity $\leq 6$, except disjoint unions.

Fig.~9 contains only two non-trivial links: the Hopf link $4_1$
 and the trefoil $6_1$.
\smallskip

\begin{tabular}{|l|c|c|c|c|c|c|c|}
\hline

three-page complexity     & 0 & 1 & 2 & 3 & 4 & 5 & 6\\
\hline

non-oriented knots    & 1 & 0 & 0 & 0 & 0 & 0 & 1\\
\hline

non-oriented links    & 0 & 0 & 0 & 0 & 1 & 0 & 0\\
\hline

non-oriented spatial 3-graphs      & 0 & 1 & 0 & 1 & 2 & 2 & 2\\
\hline

non-oriented singular knots        & 0 & 0 & 2 & 0 & 2 & 2 & 5\\
\hline
\end{tabular}
\medskip

\begin{picture}(445,70)(10,10)

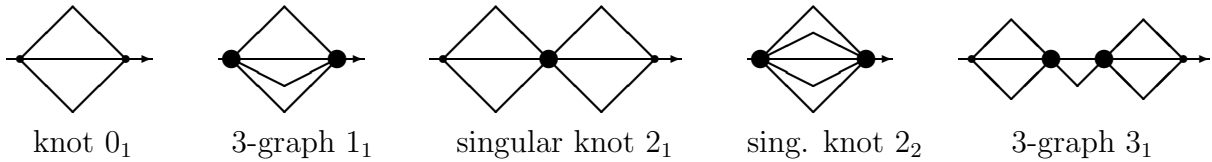
\put(100,-10){{\bf Fig.~9a.} Spatial graphs up to complexity 3.}

\put(10,15){knot $0_1$}
\put(0,50){\vector(1,0){55}}
\put(5,50){\circle*{3}}
\put(45,50){\circle*{3}}

{\thicklines
\put(5,50){\line(1,1){20}}
\put(5,50){\line(1,-1){20}}
\put(45,50){\line(-1,1){20}}
\put(45,50){\line(-1,-1){20}}
}


\put(85,15){3-graph $1_1$}
\put(80,50){\vector(1,0){55}}
\put(85,50){\circle*{7}}
\put(125,50){\circle*{7}}

{\thicklines
\put(85,50){\line(1,1){20}}
\put(85,50){\line(1,-1){20}}
\put(85,50){\line(2,-1){20}}
\put(125,50){\line(-1,1){20}}
\put(125,50){\line(-1,-1){20}}
\put(125,50){\line(-2,-1){20}}
}         


\put(170,15){singular knot $2_1$}
\put(160,50){\vector(1,0){95}}
\put(165,50){\circle*{3}}
\put(205,50){\circle*{7}}
\put(245,50){\circle*{3}}

{\thicklines
\put(165,50){\line(1,1){20}}
\put(165,50){\line(1,-1){20}}
\put(205,50){\line(-1,1){20}}
\put(205,50){\line(-1,-1){20}}
\put(205,50){\line(1,1){20}}
\put(205,50){\line(1,-1){20}}
\put(245,50){\line(-1,1){20}}
\put(245,50){\line(-1,-1){20}}
}


\put(280,15){sing. knot $2_2$}
\put(280,50){\vector(1,0){55}}
\put(285,50){\circle*{7}}
\put(325,50){\circle*{7}}

{\thicklines
\put(285,50){\line(1,1){20}}
\put(285,50){\line(1,-1){20}}
\put(285,50){\line(2,1){20}}
\put(285,50){\line(2,-1){20}}
\put(325,50){\line(-1,1){20}}
\put(325,50){\line(-1,-1){20}}
\put(325,50){\line(-2,1){20}}
\put(325,50){\line(-2,-1){20}}
}


\put(380,15){3-graph $3_1$}
\put(360,50){\vector(1,0){95}}
\put(365,50){\circle*{3}}
\put(395,50){\circle*{7}}
\put(415,50){\circle*{7}}
\put(445,50){\circle*{3}}

{\thicklines
\put(365,50){\line(1,1){15}}
\put(365,50){\line(1,-1){15}}
\put(395,50){\line(-1,1){15}}
\put(395,50){\line(-1,-1){15}}
\put(395,50){\line(1,-1){10}}
\put(415,50){\line(-1,-1){10}}
\put(415,50){\line(1,1){15}}
\put(415,50){\line(1,-1){15}}
\put(445,50){\line(-1,1){15}}
\put(445,50){\line(-1,-1){15}}
}

\end{picture}

\begin{picture}(460,130)(10,10)

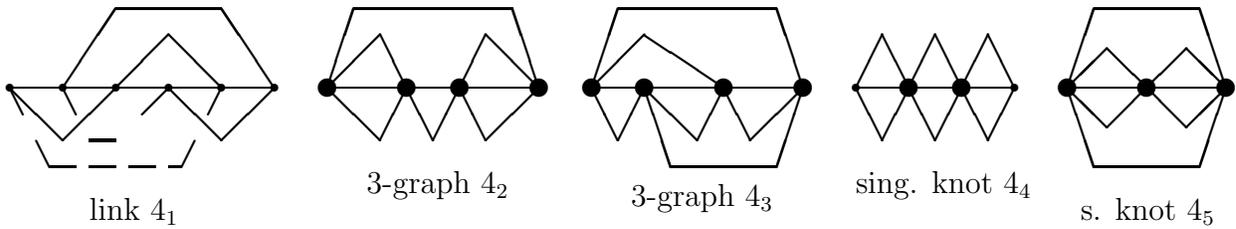
\put(100,0){{\bf Fig.~9b.} Spatial graphs with complexity 4.}

\put(30,20){link $4_1$}
\put(0,70){\line(1,0){100}}
\put(0,70){\circle*{3}}
\put(20,70){\circle*{3}}
\put(40,70){\circle*{3}}
\put(60,70){\circle*{3}}
\put(80,70){\circle*{3}}
\put(100,70){\circle*{3}}

{\thicklines
\put(0,70){\line(1,-1){20}}
\put(0,70){\line(1,-2){5}}
\put(10,50){\line(1,-2){5}}
\put(15,40){\line(1,0){10}}
\put(20,50){\line(1,1){40}}
\put(20,70){\line(1,-2){5}}
\put(20,70){\line(2,3){20}}
\put(30,50){\line(1,0){10}}
\put(30,40){\line(1,0){10}}
\put(40,100){\line(1,0){40}}
\put(45,40){\line(1,0){10}}
\put(60,40){\line(1,0){5}}
\put(60,70){\line(-1,-1){10}}
\put(60,70){\line(1,-1){20}}
\put(60,90){\line(1,-1){20}}
\put(65,40){\line(1,2){5}}
\put(80,70){\line(-1,-2){5}}
\put(100,70){\line(-1,-1){20}}
\put(100,70){\line(-2,3){20}}
}


\put(135,30){3-graph $4_2$}
\put(120,70){\line(1,0){80}}
\put(120,70){\circle*{7}}
\put(150,70){\circle*{7}}
\put(170,70){\circle*{7}}
\put(200,70){\circle*{7}}

{\thicklines
\put(120,70){\line(1,3){10}}
\put(120,70){\line(1,1){20}}
\put(120,70){\line(1,-1){20}}
\put(130,100){\line(1,0){60}}
\put(140,90){\line(1,-2){20}}
\put(140,50){\line(1,2){10}}
\put(160,50){\line(1,2){20}}
\put(170,70){\line(1,-2){10}}
\put(200,70){\line(-1,-1){20}}
\put(200,70){\line(-1,1){20}}
\put(200,70){\line(-1,3){10}}
}


\put(235,25){3-graph $4_3$}
\put(220,70){\line(1,0){80}}
\put(220,70){\circle*{7}}
\put(240,70){\circle*{7}}
\put(270,70){\circle*{7}}
\put(300,70){\circle*{7}}

{\thicklines
\put(220,70){\line(1,3){10}}
\put(220,70){\line(1,1){20}}
\put(220,70){\line(1,-2){10}}
\put(230,50){\line(1,2){10}}
\put(240,70){\line(1,-1){20}}
\put(240,70){\line(1,-3){10}}
\put(230,100){\line(1,0){60}}
\put(250,40){\line(1,0){40}}
\put(270,70){\line(-3,2){30}}
\put(270,70){\line(-1,-2){10}}
\put(270,70){\line(1,-2){10}}
\put(300,70){\line(-1,-1){20}}
\put(300,70){\line(-1,-3){10}}
\put(300,70){\line(-1,3){10}}
}


\put(320,30){sing. knot $4_4$}
\put(320,70){\line(1,0){60}}
\put(320,70){\circle*{3}}
\put(340,70){\circle*{7}}
\put(360,70){\circle*{7}}
\put(380,70){\circle*{3}}

{\thicklines
\put(320,70){\line(1,2){10}}
\put(320,70){\line(1,-2){10}}
\put(330,90){\line(1,-2){20}}
\put(330,50){\line(1,2){20}}
\put(350,90){\line(1,-2){20}}
\put(350,50){\line(1,2){20}}
\put(380,70){\line(-1,-2){10}}
\put(380,70){\line(-1,2){10}}
}


\put(405,20){s. knot $4_5$}
\put(400,70){\line(1,0){60}}
\put(400,70){\circle*{7}}
\put(430,70){\circle*{7}}
\put(460,70){\circle*{7}}

{\thicklines
\put(400,70){\line(1,3){10}}
\put(400,70){\line(1,-3){10}}
\put(400,70){\line(1,1){15}}
\put(400,70){\line(1,-1){15}}
\put(410,40){\line(1,0){40}}
\put(410,100){\line(1,0){40}}
\put(415,55){\line(1,1){30}}
\put(415,85){\line(1,-1){30}}
\put(460,70){\line(-1,3){10}}
\put(460,70){\line(-1,-3){10}}
\put(460,70){\line(-1,1){15}}
\put(460,70){\line(-1,-1){15}}
}
\end{picture}

\begin{picture}(460,120)(20,10)

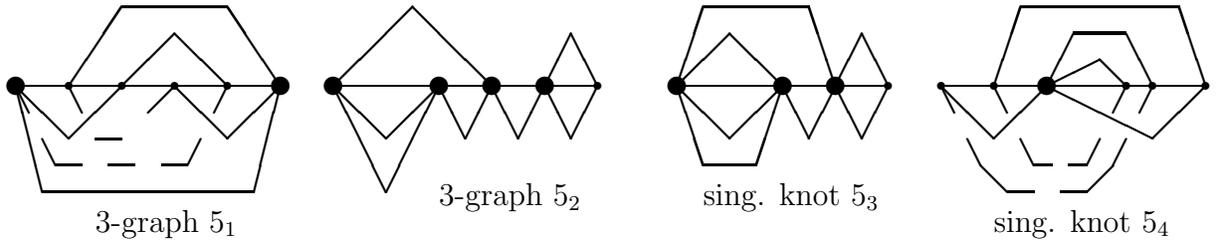
\put(100,0){{\bf Fig.~9c.} Spatial graphs with complexity 5.}

\put(40,15){3-graph $5_1$}
\put(10,70){\line(1,0){100}}
\put(10,70){\circle*{7}}
\put(30,70){\circle*{3}}
\put(50,70){\circle*{3}}
\put(70,70){\circle*{3}}
\put(90,70){\circle*{3}}
\put(110,70){\circle*{7}}

{\thicklines
\put(10,70){\line(1,-4){10}}
\put(10,70){\line(1,-1){20}}
\put(10,70){\line(1,-2){5}}
\put(20,30){\line(1,0){80}}
\put(20,50){\line(1,-2){5}}
\put(25,40){\line(1,0){10}}
\put(35,30){\line(1,0){10}}
\put(30,50){\line(1,1){40}}
\put(30,70){\line(1,-2){5}}
\put(30,70){\line(2,3){20}}
\put(40,50){\line(1,0){10}}
\put(45,40){\line(1,0){10}}
\put(50,100){\line(1,0){40}}
\put(65,40){\line(1,0){10}}
\put(70,70){\line(-1,-1){10}}
\put(70,70){\line(1,-1){20}}
\put(70,90){\line(1,-1){20}}
\put(75,40){\line(1,2){5}}
\put(90,70){\line(-1,-2){5}}
\put(110,70){\line(-1,-1){20}}
\put(110,70){\line(-2,3){20}}
\put(110,70){\line(-1,-4){10}}
}


\put(170,25){3-graph $5_2$}
\put(130,70){\line(1,0){100}}
\put(130,70){\circle*{7}}
\put(170,70){\circle*{7}}
\put(190,70){\circle*{7}}
\put(210,70){\circle*{7}}
\put(230,70){\circle*{3}}

{\thicklines
\put(130,70){\line(1,1){30}}
\put(130,70){\line(1,-1){20}}
\put(130,70){\line(1,-2){20}}
\put(160,100){\line(1,-1){30}}
\put(170,70){\line(-1,-1){20}}
\put(170,70){\line(-1,-2){20}}
\put(170,70){\line(1,-2){10}}
\put(190,70){\line(-1,-2){10}}
\put(190,70){\line(1,-2){10}}
\put(210,70){\line(-1,-2){10}}
\put(210,70){\line(1,-2){10}}
\put(210,70){\line(1,2){10}}
\put(230,70){\line(-1,-2){10}}
\put(230,70){\line(-1,2){10}}
}


\put(270,25){sing. knot $5_3$}
\put(260,70){\line(1,0){80}}
\put(260,70){\circle*{7}}
\put(300,70){\circle*{7}}
\put(320,70){\circle*{7}}
\put(340,70){\circle*{3}}

{\thicklines
\put(260,70){\line(1,3){10}}
\put(260,70){\line(1,-3){10}}
\put(260,70){\line(1,1){20}}
\put(260,70){\line(1,-1){20}}
\put(270,40){\line(1,0){20}}
\put(270,100){\line(1,0){40}}
\put(300,70){\line(-1,1){20}}
\put(300,70){\line(-1,-1){20}}
\put(300,70){\line(-1,-3){10}}
\put(300,70){\line(1,-2){10}}
\put(320,70){\line(-1,3){10}}
\put(320,70){\line(-1,-2){10}}
\put(320,70){\line(1,2){10}}
\put(320,70){\line(1,-2){10}}
\put(340,70){\line(-1,2){10}}
\put(340,70){\line(-1,-2){10}}
}


\put(380,15){sing. knot $5_4$}
\put(360,70){\line(1,0){100}}
\put(360,70){\circle*{3}}
\put(380,70){\circle*{3}}
\put(400,70){\circle*{7}}
\put(430,70){\circle*{3}}
\put(440,70){\circle*{3}}
\put(460,70){\circle*{3}}

{\thicklines
\put(360,70){\line(1,-1){20}}
\put(360,70){\line(1,-2){5}}
\put(370,50){\line(1,-2){5}}
\put(375,40){\line(1,-1){10}}
\put(380,70){\line(1,3){10}}
\put(380,70){\line(1,-2){5}}
\put(385,30){\line(1,0){10}}
\put(390,50){\line(1,-2){5}}
\put(395,40){\line(1,0){7}}
\put(405,30){\line(1,0){10}}
\put(400,70){\line(-1,-1){20}}
\put(400,70){\line(1,2){10}}
\put(400,70){\line(2,1){20}}
\put(400,70){\line(2,-1){40}}
\put(390,100){\line(1,0){60}}
\put(408,40){\line(1,0){7}}
\put(410,90){\line(1,0){20}}
\put(415,30){\line(1,1){10}}
\put(415,40){\line(1,2){5}}
\put(425,40){\line(1,2){5}}
\put(430,70){\line(-1,1){10}}
\put(430,70){\line(-1,-2){5}}
\put(440,70){\line(-1,2){10}}
\put(440,70){\line(-1,-2){5}}
\put(460,70){\line(-1,-1){20}}
\put(460,70){\line(-1,3){10}}
}

\end{picture}

\begin{picture}(460,125)(10,10)

\put(30,15){knot $6_1$}
\put(0,70){\line(1,0){120}}
\put(0,70){\circle*{3}}
\put(20,70){\circle*{3}}
\put(30,70){\circle*{3}}
\put(60,70){\circle*{3}}
\put(80,70){\circle*{3}}
\put(100,70){\circle*{3}}
\put(110,70){\circle*{3}}
\put(120,70){\circle*{3}}

{\thicklines
\put(0,70){\line(1,-2){5}}
\put(0,70){\line(1,-1){20}}
\put(10,50){\line(1,-2){5}}
\put(20,70){\line(1,3){10}}
\put(20,70){\line(1,-2){5}}
\put(20,30){\line(1,0){10}}
\put(30,70){\line(1,1){20}}
\put(30,70){\line(1,-2){5}}
\put(30,50){\line(1,-2){5}}
\put(35,40){\line(1,0){5}}
\put(30,100){\line(1,0){80}}
\put(40,50){\line(1,0){10}}
\put(40,30){\line(1,0){10}}
\put(45,40){\line(1,0){10}}
\put(60,50){\line(1,0){10}}
\put(60,30){\line(1,0){10}}
\put(65,40){\line(1,0){10}}
\put(60,70){\line(-2,-1){40}}
\put(60,70){\line(2,1){20}}
\put(50,90){\line(1,0){40}}
\put(80,70){\line(1,-1){20}}
\put(80,70){\line(-1,-2){5}}
\put(80,30){\line(1,0){10}}
\put(80,40){\line(1,0){5}}
\put(85,40){\line(1,2){5}}
\put(90,30){\line(1,2){5}}
\put(100,70){\line(-2,1){20}}
\put(100,70){\line(-1,-2){5}}
\put(110,70){\line(-1,1){20}}
\put(110,70){\line(-1,-2){5}}
\put(120,70){\line(-1,3){10}}
\put(120,70){\line(-1,-1){20}}
}


\put(165,30){3-graph $6_2$}
\put(140,70){\line(1,0){100}}
\put(140,70){\circle*{3}}
\put(160,70){\circle*{7}}
\put(180,70){\circle*{7}}
\put(200,70){\circle*{7}}
\put(220,70){\circle*{7}}
\put(240,70){\circle*{3}}

{\thicklines
\put(140,70){\line(1,2){10}}
\put(140,70){\line(1,-2){10}}
\put(150,90){\line(1,-2){20}}
\put(160,70){\line(-1,-2){10}}
\put(170,50){\line(1,2){20}}
\put(180,70){\line(1,-2){10}}
\put(190,90){\line(1,-2){20}}
\put(200,70){\line(-1,-2){10}}
\put(210,50){\line(1,2){20}}
\put(220,70){\line(1,-2){10}}
\put(240,70){\line(-1,2){10}}
\put(240,70){\line(-1,-2){10}}
}


\put(290,15){3-graph $6_3$}
\put(260,70){\line(1,0){100}}
\put(260,70){\circle*{7}}
\put(280,70){\circle*{7}}
\put(300,70){\circle*{7}}
\put(320,70){\circle*{3}}
\put(340,70){\circle*{3}}
\put(360,70){\circle*{7}}

{\thicklines
\put(260,70){\line(1,-4){10}}
\put(260,70){\line(1,-1){20}}
\put(260,70){\line(1,-2){5}}
\put(270,30){\line(1,0){80}}
\put(270,50){\line(1,-2){5}}
\put(275,40){\line(1,0){10}}
\put(280,50){\line(1,1){40}}
\put(280,70){\line(1,1){10}}
\put(280,70){\line(1,-2){5}}
\put(280,70){\line(1,3){10}}
\put(290,50){\line(1,0){10}}
\put(295,40){\line(1,0){10}}
\put(290,100){\line(1,0){60}}
\put(300,70){\line(-1,1){10}}
\put(315,40){\line(1,0){10}}
\put(320,70){\line(-1,-1){10}}
\put(320,70){\line(1,-1){20}}
\put(320,90){\line(1,-1){20}}
\put(325,40){\line(1,2){5}}
\put(340,70){\line(-1,-2){5}}
\put(360,70){\line(-1,-1){20}}
\put(360,70){\line(-1,3){10}}
\put(360,70){\line(-1,-4){10}}
}


\put(390,30){sing. knot $6_4$}
\put(380,70){\line(1,0){80}}
\put(380,70){\circle*{3}}
\put(400,70){\circle*{7}}
\put(420,70){\circle*{7}}
\put(440,70){\circle*{7}}
\put(460,70){\circle*{3}}

{\thicklines
\put(380,70){\line(1,2){10}}
\put(380,70){\line(1,-2){10}}
\put(390,50){\line(1,2){20}}
\put(390,90){\line(1,-2){20}}
\put(410,50){\line(1,2){20}}
\put(410,90){\line(1,-2){20}}
\put(430,50){\line(1,2){20}}
\put(430,90){\line(1,-2){20}}
\put(460,70){\line(-1,2){10}}
\put(460,70){\line(-1,-2){10}}
}
\end{picture}


\begin{picture}(460,100)(10,10)

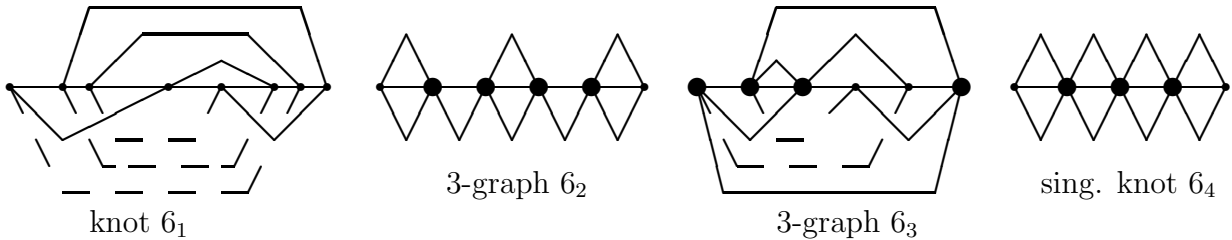
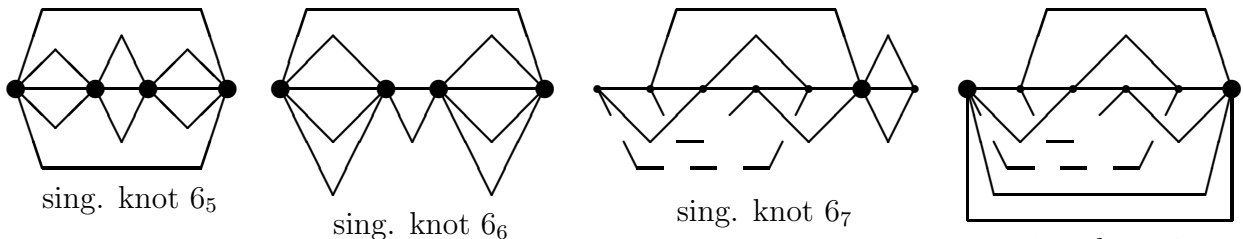
\put(100,-5){{\bf Fig.~9d.} Spatial graphs with complexity 6.}

\put(10,25){sing. knot $6_5$}
\put(0,70){\line(1,0){80}}
\put(0,70){\circle*{7}}
\put(30,70){\circle*{7}}
\put(50,70){\circle*{7}}
\put(80,70){\circle*{7}}

{\thicklines
\put(0,70){\line(1,3){10}}
\put(0,70){\line(1,-3){10}}
\put(0,70){\line(1,1){15}}
\put(0,70){\line(1,-1){15}}
\put(10,100){\line(1,0){60}}
\put(10,40){\line(1,0){60}}
\put(30,70){\line(-1,1){15}}
\put(30,70){\line(-1,-1){15}}
\put(30,70){\line(1,2){10}}
\put(30,70){\line(1,-2){10}}
\put(50,70){\line(-1,2){10}}
\put(50,70){\line(-1,-2){10}}
\put(50,70){\line(1,1){15}}
\put(50,70){\line(1,-1){15}}
\put(80,70){\line(-1,1){15}}
\put(80,70){\line(-1,-1){15}}
\put(80,70){\line(-1,3){10}}
\put(80,70){\line(-1,-3){10}}
}


\put(120,15){sing. knot $6_6$}
\put(100,70){\line(1,0){100}}
\put(100,70){\circle*{7}}
\put(140,70){\circle*{7}}
\put(160,70){\circle*{7}}
\put(200,70){\circle*{7}}

{\thicklines
\put(100,70){\line(1,3){10}}
\put(100,70){\line(1,1){20}}
\put(100,70){\line(1,-1){20}}
\put(100,70){\line(1,-2){20}}
\put(110,100){\line(1,0){80}}
\put(140,70){\line(-1,1){20}}
\put(140,70){\line(-1,-1){20}}
\put(140,70){\line(-1,-2){20}}
\put(140,70){\line(1,-2){10}}
\put(160,70){\line(-1,-2){10}}
\put(160,70){\line(1,1){20}}
\put(160,70){\line(1,-1){20}}
\put(160,70){\line(1,-2){20}}
\put(200,70){\line(-1,1){20}}
\put(200,70){\line(-1,-1){20}}
\put(200,70){\line(-1,-2){20}}
\put(200,70){\line(-1,3){10}}
}

\put(250,20){sing. knot $6_7$}
\put(220,70){\line(1,0){120}}
\put(220,70){\circle*{3}}
\put(240,70){\circle*{3}}
\put(260,70){\circle*{3}}
\put(280,70){\circle*{3}}
\put(300,70){\circle*{3}}
\put(320,70){\circle*{7}}
\put(340,70){\circle*{3}}

{\thicklines
\put(220,70){\line(1,-1){20}}
\put(220,70){\line(1,-2){5}}
\put(230,50){\line(1,-2){5}}
\put(235,40){\line(1,0){10}}
\put(240,70){\line(1,3){10}}
\put(240,70){\line(1,-2){5}}
\put(240,50){\line(1,1){40}}
\put(250,50){\line(1,0){10}}
\put(250,100){\line(1,0){60}}
\put(255,40){\line(1,0){10}}
\put(280,70){\line(-1,-1){10}}
\put(280,70){\line(1,-1){20}}
\put(275,40){\line(1,0){10}}
\put(285,40){\line(1,2){5}}
\put(300,70){\line(-1,1){20}}
\put(300,70){\line(-1,-2){5}}
\put(320,70){\line(-1,3){10}}
\put(320,70){\line(-1,-1){20}}
\put(320,70){\line(1,2){10}}
\put(320,70){\line(1,-2){10}}
\put(340,70){\line(-1,2){10}}
\put(340,70){\line(-1,-2){10}}
}

\put(380,5){sing. knot $6_8$}
\put(360,70){\line(1,0){100}}
\put(360,70){\circle*{7}}
\put(380,70){\circle*{3}}
\put(400,70){\circle*{3}}
\put(420,70){\circle*{3}}
\put(440,70){\circle*{3}}
\put(460,70){\circle*{7}}

{\thicklines
\put(360,70){\line(1,-1){20}}
\put(360,70){\line(1,-4){10}}
\put(360,70){\line(0,-1){50}}
\put(360,70){\line(1,-2){5}}
\put(360,20){\line(1,0){100}}
\put(370,30){\line(1,0){80}}
\put(370,50){\line(1,-2){5}}
\put(375,40){\line(1,0){10}}
\put(380,70){\line(1,3){10}}
\put(380,70){\line(1,-2){5}}
\put(380,50){\line(1,1){40}}
\put(390,50){\line(1,0){10}}
\put(390,100){\line(1,0){60}}
\put(395,40){\line(1,0){10}}
\put(420,70){\line(-1,-1){10}}
\put(420,70){\line(1,-1){20}}
\put(415,40){\line(1,0){10}}
\put(425,40){\line(1,2){5}}
\put(440,70){\line(-1,1){20}}
\put(440,70){\line(-1,-2){5}}
\put(460,70){\line(-1,3){10}}
\put(460,70){\line(-1,-1){20}}
\put(460,70){\line(-1,-4){10}}
\put(460,70){\line(0,-1){50}}
}

\end{picture}



\end{document}